\documentclass[english]{siamltex}
\usepackage[english]{babel}

\usepackage[bookmarks, colorlinks=true, plainpages = false, linkcolor = blue, citecolor = green, urlcolor = blue, filecolor = blue]{hyperref}

\usepackage[T1]{fontenc}
\usepackage[latin1]{inputenc}
\usepackage{textcomp}

\usepackage{amsthm}
\usepackage{amsmath}
\usepackage{amsfonts}
\usepackage{mathrsfs}
\usepackage{dsfont}

\usepackage{amssymb}
\usepackage{enumerate}
\usepackage{graphicx}
\usepackage{mathrsfs}
\usepackage{algorithm}
\usepackage{algorithmicx}
\usepackage{algpseudocode}
\usepackage{array}
\usepackage{caption}
\usepackage{booktabs}
\usepackage{url}

\usepackage{color}

\usepackage[normalem]{ulem}  

\title{Bayesian quadrature and energy minimization\\
 for space-filling design}

\author{Luc Pronzato\footnotemark[2]\ \ and Anatoly Zhigljavsky\footnotemark[1] }

\DeclareMathOperator{\Arg}{Arg}

\DeclareMathOperator{\Id}{Id}

\DeclareMathOperator{\CR}{\mathsf{CR}}
\DeclareMathOperator{\PR}{\mathsf{PR}}

\DeclareMathOperator{\IMSPE}{\mathsf{IMSPE}}

\DeclareMathOperator{\tr}{trace}

\numberwithin{equation}{section}

\newcommand{\bea}{\begin{eqnarray*}}
\newcommand{\eea}{\end{eqnarray*}}
\newcommand{\be}{\begin{eqnarray}}
\newcommand{\ee}{\end{eqnarray}}
\def\ra{\rightarrow}
\def\e1{\mathrm{e}}
\def\mp{\partial}
\def\0b{\boldsymbol{0}}
\def\1b{\boldsymbol{1}}

\def\dd{\mbox{\rm d}}
\def\Supp{\mathrm{Supp}}

\def\SB{{\mathscr B}}

\def\SE{{\mathscr E}}
\def\SH{{\mathcal H}}

\def\SJ{{\mathcal J}}

\def\SM{{\mathscr M}}
\def\SO{{\mathcal O}}
\def\SP{{\mathscr P}}

\def\SSS{\mathbb{S}}
\def\SSp{{\mathscr S}}
\def\SN{{\mathscr N}}

\def\SX{{\mathscr X}}

\def\ma{\alpha}
\def\betah{\hat\beta}
\def\mg{\gamma}
\def\ml{\lambda}
\def\me{\epsilon}

\def\ms{\sigma}
\def\mt{\theta}

\def\mo{\omega}
\def\mob{\boldsymbol{\omega}}
\def\mO{\Omega}

\def\psib{\boldsymbol{\psi}}
\def\Psib{\boldsymbol{\Psi}}
\def\Sigmab{\boldsymbol{\Sigma}}
\def\Lambdab{\boldsymbol{\Lambda}}

\def\betab{\boldsymbol{\beta}}
\def\betabh{\hat{\boldsymbol{\beta}}}

\def\Ex{\mathsf{E}}
\def\Exx{\mathbb{E}}

\def\Ab{\mathbf{A}}

\def\Bb{\mathbf{B}}

\def\eb{\mathbf{e}}

\def\hb{\mathbf{h}}
\def\Hb{\mathbf{H}}
\def\Ib{\mathbf{I}}

\def\kb{\mathbf{k}}
\def\Kb{\mathbf{K}}
\def\Mb{\mathbf{M}}

\def\Pb{\mathbf{P}}
\def\pb{\mathbf{p}}
\def\Qb{\mathbf{Q}}
\def\rb{\mathbf{r}}
\def\Rb{\mathbf{R}}
\def\sb{\mathbf{s}}
\def\Sb{\mathbf{S}}

\def\ub{\mathbf{u}}
\def\Ub{\mathbf{U}}
\def\vb{\mathbf{v}}
\def\Vb{\mathbf{V}}
\def\wb{\mathbf{w}}

\def\xb{\mathbf{x}}
\def\Xb{\mathbf{X}}
\def\yb{\mathbf{y}}

\def\zb{\mathbf{z}}
\def\Zb{\mathbf{Z}}

\theoremstyle{plain}

\newtheorem{theo}{Theorem}[section]

\newtheorem{defi}{Definition}[section]

\newtheorem{lem}{Lemma}[section]

\theoremstyle{definition}

\newtheorem{remark}{Remark}[section]

\newtheorem{example}{Example}[section]

\newcommand{\carre} {\mbox{}~\hfill\rule{2mm}{2mm}}
\newcommand{\fin} {\mbox{}~\hfill{\lower-0.3ex\hbox{$\triangleleft$}}}

\newcommand{\vsp} {\vspace{0.3cm}}

\begin{document}
\maketitle
\renewcommand{\thefootnote}{\fnsymbol{footnote}}
\footnotetext[1]{Cardiff University, UK}
\footnotetext[2]{CNRS, Universit\'e C\^ote d'Azur, I3S, France}
\footnotetext[3]{Luc.Pronzato@cnrs.fr (corresponding author)}

\renewcommand{\thefootnote}{\arabic{footnote}}

\begin{abstract}
A standard objective in computer experiments is to approximate the behaviour of an unknown function on a compact domain from a few evaluations inside the domain. When little is known about the function, space-filling design is advisable: typically, points of evaluation spread out across the available space are obtained by minimizing a geometrical (for instance, covering radius) or a discrepancy criterion measuring distance to uniformity.
The paper investigates connections between design for integration (quadrature design), construction of the (continuous) BLUE for the location model, space-filling design, and minimization of energy (kernel discrepancy) for signed measures. Integrally strictly positive definite kernels define strictly convex energy functionals, with an equivalence between the notions of potential and directional derivative, showing the strong relation between discrepancy minimization and more traditional design of optimal experiments. In particular, kernel herding algorithms, which are special instances of vertex-direction methods used in optimal design, can be applied to the construction of point sequences with suitable space-filling properties.
\end{abstract}



\section{Introduction}\label{S:Intro}

The design of computer experiments, where observations of a real physical phenomenon are replaced by simulations of a complex mathematical model (e.g., based on PDEs), has emerged as a full discipline, central to uncertainty quantification. The final objective of the simulations is often goal-oriented, that is, precisely defined. It may correspond for example to the optimization of the response of a system with respect to its input factors, or to the estimation of the probability that the response will exceed a given threshold when input factors have a given probability distribution. Achieving this objective generally requires sequential learning of the behavior of the response in a particular domain of interest for input factors: the region where the response is close to its optimum, or is close to the given threshold; see, e.g., the references in \cite{Ginsbourger2017}. When simulations are computationally expensive, sequential inference based on the direct use of the mathematical model is unfeasible due to the large number of simulations required and simplified prediction models, approximating the simulated response, have to be used. A most popular approach relies on Gaussian process modelling, where the response (unknown prior to simulation) is considered as the realization of a Gaussian Random Field (RF), with parameterized mean and covariance, and Bayesian inference gives access to the posterior distribution of the RF (after simulation). Typically, in a goal-oriented approach based on stepwise-uncertainty reduction \cite{BectBG2016, BectGLPV2012}, the prediction model is used to select the input factors to be used for the next simulation, the selection being optimal in terms of predicted uncertainty on the target. The construction of a first, possibly crude, prediction model is necessary to initialize the procedure. This amounts at approximating the behaviour of an unknown function $f$ (the model response) on a compact domain $\SX\subset\mathds{R}^d$ (the feasible set for $d$ input factors) from a few evaluations inside the domain. That is the basic design objective we shall keep in mind throughout the paper, although we may use diverted paths where approximation/prediction will be shadowed by other objectives, integration in particular.

In general, little is known about the function {\em a priori}, and it seems intuitively reasonable to spread out points of evaluation across the available space; see \cite{BiedermannD2001}. Such space-filling designs can be obtained by optimizing a geometrical measure of dispersion or a discrepancy criterion measuring distance to uniformity. When using a Gaussian RF model, minimizing the Integrated Mean-Squared Prediction Error (IMSPE) is also a popular approach, although not very much used due to its apparent complexity, see, e.g., \cite{GP-SIAM_2014, GorodetskyK2016}. The paper promotes the use of designs optimized for integration with respect to the uniform measure for their good space-filling properties. It gives a survey of recent results on energy functionals that measure distance to uniformity and places recent approaches proposed for space-filling design, such as \cite{MakJ2017}, in a general framework and perspective encompassing design for integration, construction of the (continuous) Best Linear Unbiased Estimator (BLUE) in a location model with correlated errors, and minimization of energy (kernel discrepancy) for signed measures.

We start by a quick introduction to Bayesian function approximation and integration (Section~\ref{S:RFM}), where the function is considered as the realization of a Gaussian RF with covariance structure defined by some kernel $K$.
Exploiting recent results on the minimization of energy functionals \cite{DamelinHRZ2010, SejdinovicSGF2013, SriperumbudurGFSL2010}, we show in Section~\ref{S:KDEP}  that integrally strictly positive definite kernels define strictly convex energy functionals, with an equivalence between the notions of potential and directional derivative that reveals the strong relation between discrepancy minimization and more traditional design of optimal experiments. We show that Bayesian integration is equivalent to the construction of the BLUE in a model with modified correlation structure, so that the two associated design problems coincide. We also show that the posterior variance in Bayesian integration corresponds to the minimum of a squared kernel discrepancy for signed measures with total mass one and to the minimum of an energy functional for a reduced kernel.
Since the posterior variance criterion in Bayesian integration takes a very simple form, its minimization constitutes an attractive alternative to the minimization of the IMSPE criterion for space-filling design. This is considered in Section~\ref{S:Empirical}. We consider in particular kernel herding algorithms from machine learning, which are special instances of vertex-direction methods used in optimal design and can be used for the construction of point sequences with suitable space-filling properties (any-time designs).
Several auxiliary results are given in appendix. Appendix~\ref{S:CV-properties} provides convergence properties of algorithms presented in Section~\ref{S:Empirical}. Extension to design for the simultaneous estimation of several integrals is considered in Appendix~\ref{S:BQ-several}. A Karhunen-Lo\`eve expansion of the RF model is considered in Appendix~\ref{S:KL}, that yields a Bayesian linear model for which minimization of the posterior variance in Bayesian integration corresponds to a c-optimal design problem.

\section{Random-field models for function approximation and integration}\label{S:RFM}

\subsection{Space-filling design and kernel choice for function approximation}
\label{S:GP regression}

Let $K(\cdot,\cdot)$ denote a symmetric positive definite kernel on $\SX\times\SX$, with associated Reproducing Kernel Hilbert Space (RKHS) $\SH_K$. Denote $K_\xb(\cdot)=K(\xb,\cdot)$ and $\langle \cdot,\cdot \rangle_K$ the scalar product in $\SH_K$, so that the reproducing property gives $\langle f,K_\xb \rangle_K=f(\xb)$ for any $f\in\SH_K$.

Consider first the common framework where the function $f$ to be approximated is supposed to belong to $\SH_K$. Let $\eta_n(\xb)=\sum_{i=1}^n w_i f(\xb_i)=\wb_n^T\yb_n$ be a linear predictor of $f(\xb)$ based on evaluations of $f$ at the $n$-point design $\Xb_n=\{\xb_1,\ldots,\xb_n\}$, with $\xb_i\in\SX$ for all $i$. Throughout the paper we denote $\wb_n=(w_1,\ldots,w_n)^T$, $\yb_n=[f(\xb_1),\ldots,f(\xb_n)]^T$, $\kb_n(\cdot)=[K_{\xb_1}(\cdot),\ldots,K_{\xb_n}(\cdot)]^T$ and $\{\Kb_n\}_{i,j}=K(\xb_i,\xb_j)$, $i,j=1,\ldots,n$.
Cauchy-Schwarz inequality gives the classical result
\begin{eqnarray*}
|f(\xb)-\eta_n(\xb)| = \left| f(\xb)- \sum_{i=1}^n w_i f(\xb_i)\right| &=& \left| \langle f, K_\xb-\sum_{i=1}^n w_i K_{\xb_i} \rangle_K \right| \\
&\leq& \|f\|_{\SH_K} \, \left\|K_\xb-\sum_{i=1}^n w_i K_{\xb_i}\right\|_{\SH_K} \,,
\end{eqnarray*}
where $\|f\|_{\SH_K}$ depends on $f$  but not on $\Xb_n$, and $\rho_n(\xb,\wb)=\left\|K_\xb-\sum_{i=1}^n w_i K_{\xb_i}\right\|_{\SH_K}$ depends on $\Xb_n$ (and $\wb_n$) but not on $f$. Suppose that $\Kb_n$ has full rank. For a given $\Xb_n$, the Best Linear Predictor (BLP) minimizes
$\rho_n(\xb,\wb)$ and corresponds to $\eta_n^*(\xb) = (\wb_n^*)^T\yb_n$,
with $\wb_n^*=\wb_n^*(\xb)=\Kb_n^{-1}\kb_n(\xb)$, which gives ${\rho_n^*}^2(\xb)= \rho_n^2(\xb,\wb_n^*) = K(\xb,\xb)-\kb_n^T(\xb)\Kb_n^{-1} \kb_n(\xb)$.

A less restrictive assumption on $f$ is to suppose that it corresponds to a realization of a RF $Z_x$, with zero mean ($\Exx\{Z_x\}=0$) and covariance $\Exx\{Z_x Z_{x'}\}=\ms^2\, K(\xb,\xb')$ for all $\xb$, $\xb'$ in~$\SX$, $\ms>0$. Then, straightforward calculation shows that $\eta_n^*(\xb)$ is still the BLP (the posterior mean if $Z_x$ is Gaussian), and $\ms^2\, {\rho_n^*}^2(\xb)$ is the Mean-Squared Prediction Error (MSPE) at $\xb$. This construction corresponds to simple kriging; see, e.g., \cite{AuffrayBM2012, VazquezB2011}.
IMSPE-optimal designs minimize the integrated squared error $\IMSPE(\Xb_n)= \ms^2\, \int_\SX {\rho_n^*}^2(\xb) \dd\mu(\xb)$, with $\mu$ generally taken as the uniform probability measure on~$\SX$,
see, e.g., \cite{GP-SIAM_2014, GorodetskyK2016, SacksWMW89}.

IMSPE-optimal designs $\Xb_n^*$ depend on the chosen $K$. It is well known that the asymptotic rate of decrease of $\IMSPE(\Xb_n^*)$ with $n$ depends on the regularity properties of $K$ (the same is true for the integration problem); see for instance \cite{RitterWW95}. It is rather usual to take $K$ stationary (translation invariant), i.e., satisfying $K(\xb,\xb')=\Psi(\xb-\xb')$ for all $\xb$ and $\xb'$, with $\Psi$ in some parametric class selected according to prior knowledge on the smoothness properties of $f$. A typical example is the Mat\'ern class of covariances, see \cite[Chap.~2]{Stein99}. On the other hand, for reasons explained in Section~\ref{S:Intro}, computer experiments usually involve small values of $n$, and the asymptotic behavior of the approximation error is hardly observed. Its behavior on a short horizon is much more important and strongly depends on the correlation lengths in $K$, which are difficult to choose {\em a priori}. Robustness with respect to the choice of $K$ favours space-filling designs, where the $\xb_i$ are suitably spread over~$\SX$. Noticeably, it is shown in \cite{Schaback95} that for translation invariant and isotropic kernels (i.e., such that $K(\xb,\xb')=\Psi(\|\xb-\xb'\|)$, with $\|\cdot\|$ the Euclidean distance in $\mathds{R}^d$), one has $\rho_n^2(\xb) \leq S_K[h_r(\xb)]$ for some increasing function $S_K(\cdot)$. Here $h_r(\xb)=\max_{\|\xb-\xb'\|\leq r} \min_{1\leq i \leq n} \|\xb'-\xb_i\|$ measures the density of design points $\xb_i$ around $\xb$, with $r$ a fixed positive constant. It satisfies $\max_{\xb\in\SX} h_r(\xb) \geq \max_{\xb\in\SX} h_0(\xb)=\CR(\Xb_n)$, with
\be\label{CR}
\CR(\Xb_n) = \max_{\xb\in\SX} \min_{1\leq i \leq n} \|\xb-\xb_i\| \,,
\ee
the \emph{covering radius} of $\Xb_n$: $\CR(\Xb_n)$ defines the smallest $r$ such that the $n$ closed balls of radius $r$ centred at the $\xb_i$ cover~$\SX$. $\CR(\Xb_n)$ is also called the dispersion of $\Xb_n$ \cite[Chap.~6]{Niederreiter92} and corresponds to the minimax-distance criterion \cite{JohnsonMY90} used in space-filling design.
Loosely speaking, the property $\rho_n^2(\xb) \leq S_K[h(\xb)]$ quantifies the intuition that designs with a small value of $\CR$ provide precise predictions over~$\SX$ since for any $\xb$ in~$\SX$ there always exists a design point $\xb_i$ at proximity where $f(\xb_i)$ has been evaluated. Another standard geometrical criterion of spreadness is the \emph{packing radius}
\be\label{PR}
\PR(\Xb_n) = \frac12\, \min_{i\neq j} \|\xb_i-\xb_j\| \,.
\ee
It corresponds to the largest $r$ such that the $n$ open balls of radius $r$ centred at the $\xb_i$ do not intersect; $2 \PR(\cdot)$ corresponds to the maximin-distance criterion \cite{JohnsonMY90} often used in computer experiments.

In this paper, we shall adopt the following point of view. We do not intend to construct designs adapted to a particular $K$ chosen from {\em a priori} knowledge on $f$. Neither shall we estimate the parameters in $K$ (such as correlation lengths) when $K$ is taken from a parametric class. We shall rather consider the kernel $K$ as a tool for constructing a space-filling design, the quality of which will be measured in particular through the value of $\CR$. The motivation is twofold: ($i$) the construction will be much easier than the direct minimization of $\CR$, ($ii$) it will facilitate the construction of \emph{sequences of points} suitably spread over~$\SX$ (\emph{any-time} space-filling designs).

\subsection{Bayesian quadrature}
\label{S:BQ}

Denote by $\SM=\SM[\SX]$ the set of finite signed Borel measures on a nonempty set~$\SX$, and by $\SM(q)$, $q\in\mathds{R}$, the set of signed measures with total mass $q$: $\SM(q)=\{\mu\in\SM: \mu(\SX)=q\}$. The set of Borel probability measures on~$\SX$ is denoted by $\SM^+(1)$, $\SM^+$ is the set of finite positive measures on~$\SX$. Typical applications correspond to~$\SX$ being a compact subset of $\mathds{R}^d$ for some $d \in \mathds{N}$.

Suppose we wish to integrate a real function defined on~$\SX$ with respect to $\mu\in\SM^+(1)$.
Assume that $\Ex_\mu\{|f(X)|\}<+\infty$ and denote
$$
I_\mu(f)=\Ex_\mu\{f(\Xb)\}=\int_\SX f(\xb)\, \dd\mu(\xb) \,.
$$
We set a prior on $f$, and assume that $f$ is the realization of a Gaussian RF, with covariance $\ms^2\,K(\cdot,\cdot)$, $\ms^2>0$, and unknown mean $\beta_0$; that is, we consider the location model with correlated errors
\begin{equation}\label{model1}
  f(\xb)=\beta_0+Z_x \,,
\end{equation}
where $\Exx\{Z_x\}=0$ and $\Exx\{Z_x Z_{x'}\}=\ms^2\, K(\xb,\xb')$ for all $\xb,\xb'\in\SX$, $\ms>0$. Regression models more general than \eqref{model1} are considered in Appendix~\ref{S:BQ-several}. Here $K$ is a symmetric Positive Definite (PD) kernel; that is, $K(\xb,\xb')=K(\xb',\xb)$, and for all $n\in\mathds{N}$ and all pairwise different $\xb_1,\ldots,\xb_n\in\SX$, the matrix $\Kb_n$ is non-negative definite; if $\Kb_n$ is positive definite, then $K$ is called Strictly Positive Definite (SPD).
Note that $K^2(\xb,\xb') \leq K(\xb,\xb)K(\xb',\xb')<+\infty$ for all $\xb,\xb'\in\SX$ since $K$ corresponds to a covariance.
We will call a general kernel $K$ bounded when $K(\xb,\xb)<\infty$ for all $\xb\in\SX$, and uniformly bounded when there is a constant $C$ such that $K(\xb,\xb) \leq C$ for all $\xb\in\SX$. Any PD kernel is bounded.

Similarly to Section~\ref{S:GP regression}, we denote by $\SH_K$ the associated RKHS and by $\langle \cdot,\cdot \rangle_K$ the scalar product in $\SH_K$. The assumption that $K$ is bounded will be relaxed in Section~\ref{S:EP} where we shall also consider singular kernels, but throughout the paper we assume that $K$ is symmetric, $K(\xb,\xb')=K(\xb',\xb)$ for all $\xb,\xb'\in\SX$. Also, we always assume, as in \cite[Sect.~2.1]{Fuglede60}, that either $K$ is non-negative on $\SX \times \SX$, or~$\SX$ is compact.

We set a vague prior on $\beta_0$ and assume that $\beta_0 \sim \SN(\betah_0^0,\ms^2\,A)$ with $A\ra+\infty$. This amounts to setting $1/A=0$ in all Bayesian calculations; the choice of $\betah_0^0$ is then irrelevant. Suppose that $f$ has been evaluated at the $n$-point design $\Xb_n=\{\xb_1,\ldots,\xb_n\}\in\SX^n$. We assume that $\Kb_n$ has full rank. For any $\xb\in\SX$, the posterior distribution of $f(\xb)$ (conditional on $\ms^2$ and $K$) is normal, with mean
\begin{equation*}\label{predict1}
  \hat \eta_n(\xb) = \betah_0^n + \kb_n^T(\xb)\Kb_n^{-1}(\yb_n-\betah_0^n \1b_n)
\end{equation*}
and variance (mean-squared error)
\begin{equation}\label{MSE1}
\ms^2\rho_n^2(\xb)= \ms^2\,\left[K(\xb,\xb)-\kb_n^T(\xb)\Kb_n^{-1}\kb_n(\xb)+\frac{(1-\kb_n^T(\xb)\Kb_n^{-1}\1b_n)^2}{\1b_n^T\Kb_n^{-1}\1b_n} \right] \,,
\end{equation}
where
\begin{equation}\label{beta01}
\betah_0^n = \frac{\1b_n^T\Kb_n^{-1}\yb_n}{\1b_n^T\Kb_n^{-1}\1b_n}
\end{equation}
and $\1b_n$ is the $n$-dimensional vector $(1,\ldots,1)^T$,
see for instance \cite[Chap.~4]{SantnerWN2003}. The posterior mean of $I_\mu(f)$ is thus
\begin{equation}\label{In1}
\widehat I_n = \int_\SX \hat \eta_n(\xb)\, \, \dd\mu(\xb) = \Ex_\mu\{\hat \eta_n(\Xb)\} = \betah_0^n + \pb_n(\mu)^T\Kb_n^{-1}(\yb_n-\betah_0^n \1b_n)\,,
\end{equation}
with
\be\label{hb}
\pb_n(\mu)=(P_\mu(\xb_1),\ldots,P_\mu(\xb_n))^T \,,
\ee
where, for any $\nu\in\SM$ and $\xb\in\SX$, we denote
\be \label{kh}
P_\nu(\xb)= \int_\SX K(\xb,\xb')\, \dd\nu(\xb') \,.
\ee
$P_\nu(\cdot)$ is called the kernel imbedding of $\nu$ into $\SH_K$, see \cite[Def.~9]{SejdinovicSGF2013};
$P_\nu(\xb)$ is well defined and finite for any $\nu\in\SM$ and $\xb\in\SX$ when $K$ is uniformly bounded. On the other hand, there always exists $\nu\in\SM$ such that $P_\nu(\xb)$ is infinite for all $\xb\in\SX$ when $K$ is not uniformly bounded on~$\SX$.
The function $P_\nu(\cdot)$ is called potential in potential theory, see Section~\ref{S:EP}.

Similarly to \eqref{MSE1}, we obtain that the posterior variance of $I_\mu(f)$ becomes
\begin{eqnarray}
\ms^2 s_n^2  &=&  \ms^2\, \left[ \SE_K(\mu) - \pb_n^T(\mu) \Kb_n^{-1} \pb_n(\mu) + \frac{(1-\pb_n^T(\mu)\Kb_n^{-1}\1b_n)^2}{\1b_n^T\Kb_n^{-1}\1b_n} \right]  \,, \label{sn1}
\end{eqnarray}
where, for any $\nu\in\SM$, we denote
\be \label{overline K}
\SE_K(\nu) = \int_{\SX^2} K(\xb,\xb')\, \dd\nu(\xb)\dd\nu(\xb') \,.
\ee
This is one of the key notions in potential theory, called the energy of $\nu$; see Section~\ref{S:EP}.
For $\mu$ in $\SM^+(1)$, we have $\SE_K(\mu) = \Ex_\mu\{K(\Xb,\Xb')\}$
where $\Xb$ and $\Xb'$ are independently identically distributed (i.i.d.) with $\mu$. The quantity $-\SE_K(\mu)$ corresponds to the quadratic entropy introduced by C.R. Rao \cite{Rao1982a}; see also Remark~\ref{R:Bregman}.
Define
\be\label{M_K^a}
\SM_K^\ma = \left\{\nu\in\SM: \int_\SX K^\ma(\xb,\xb)\, \dd|\nu|(\xb) < +\infty \right\}\,, \ \ma>0 \,. 
\ee
When $\mu\in\SM_K^{1/2}$, the reproducing property and Cauchy-Schwarz inequality imply that
\be
\SE_K(\mu) &=& \int_{\SX^2} \langle K(\cdot,\xb),K(\cdot,\xb') \rangle_K\, \dd\mu(\xb)\dd\mu(\xb') \nonumber \\
&& \leq \left[\int_\SX K^{1/2}(\xb,\xb)\,\dd|\mu|(\xb)\right]^2 < +\infty\,. \label{CS:M1/2}
\ee
When $\beta_0$ is assumed to be known (equal to zero for instance), we simply substitute $\beta_0$ for $\betah_0^n$ in \eqref{In1} and the posterior variance is
\be\label{sn10}
\ms^2 s_{n,0}^2 = \ms^2\, \left[ \SE_K(\mu) - \pb_n^T(\mu) \Kb_n^{-1} \pb_n(\mu) \right] \,.
\ee

Bayesian quadrature relies on the estimation of $I_\mu(f)$ by $\widehat I_n$. An optimal design for estimating $I_\mu(f)$ should minimize $s_n^2$ given by \eqref{sn1}. One may refer to \cite{Diaconis88} for a historical perspective and to \cite{HennigOG2015} for a recent exposition on Bayesian numerical computation. The framework presented above corresponds to that considered in \cite{O'Hagan91}, restricted to the case (recommended in that paper) where the known trend function is simply the constant 1 (which corresponds to the presence of an unknown mean $\beta_0$ in the model \eqref{model1}). In Section~\ref{S:Empirical}, we shall see that $s_{n,0}^2$ is equal to the minimum value of a (squared) kernel discrepancy between the measure $\mu$ and a signed measure supported on $\Xb_n$, and that $s_n^2$ corresponds to the minimum of a squared discrepancy for signed measures that are constrained to have total mass one, and also corresponds to the minimum of an energy functional for a modified kernel $K_\mu$. Note that $\ms^2 s_n^2 \leq \IMSPE(\Xb_n)=\ms^2 \int_\SX \rho_n^2(\xb)\, \dd\mu(\xb)$ (which requires $\mu\in\SM_K^1\subset\SM_K^{1/2}$ to be well defined).
One of the key ideas of the paper is that space-filling design may be based on the minimization of $s_n^2$ rather than the minimization of $\IMSPE(\Xb_n)$.


\section{Kernel discrepancy, energy and potentials} \label{S:KDEP}
\subsection{Maximum mean discrepancy}\label{S:MMD}
Suppose now that $K$ is bounded and $f$ belongs to the RKHS $\SH_K$. Let $\mu$ and $\nu$ be two probability measures in $\SM^+(1)\cap\SM_K^{1/2}$. Since $f\in\SH_K$, using the reproducing property, we obtain $I_\mu(f)=\int_\SX \langle f, K_\xb \rangle_K \, \dd\mu(\xb)$, $I_\nu(f)=\int_\SX \langle f, K_\xb \rangle_K \, \dd\nu(\xb)$ and
$$
\left|I_\mu(f) - I_{\nu}(f) \right|  = \left| \int_\SX \langle f, K_\xb \rangle_K \, \dd(\mu-\nu)(\xb)\right| = \left| \langle f, P_\mu - P_\nu \rangle_K \right| \,,
$$
with $P_\mu(\cdot)$ and $P_\nu(\cdot)$ the kernel imbeddings \eqref{kh}. Define
\be\label{pseudoM-kernels}
\mg_K(\mu,\nu) = \| P_\mu - P_\nu\|_{\SH_K} \,.
\ee
Cauchy-Schwarz inequality yields the Koksma-Hlawka type inequality \cite[Chap.~2]{Niederreiter92}
$\left|I_\mu(f) - I_{\nu}(f) \right| \leq \| f\|_{\SH_K} \mg_K(\mu,\nu)$, and
\be\label{pseudoM-PD}
\mg_K(\mu,\nu) = \sup_{\|f\|_{\SH_K}=1} |I_\mu(f)-I_\nu(f)|  \,,
\ee
see, e.g., \cite[Th.~1]{SriperumbudurGFSL2010}. Also, the expansion of $\| P_\mu - P_\nu\|^2_{\SH_K}$ gives
\begin{eqnarray}
\mg_K(\mu,\nu) &=& \left(\|P_\mu\|_{\SH_K}^2 + \|P_\nu\|_{\SH_K}^2 - 2 \langle P_\mu, P_\nu \rangle_K \right)^{1/2} \nonumber \\
&=& \left( \int_{\SX^2} K(\xb,\xb')\, \dd(\nu-\mu)(\xb)\, \dd(\nu-\mu)(\xb') \right)^{1/2} \,. \label{kernel discrepancy}
\end{eqnarray}
Therefore, $\mg_K(\cdot,\cdot)$ is at the same time a pseudometric between kernel imbeddings \eqref{pseudoM-kernels} and an integral pseudometric on probability distributions \eqref{pseudoM-PD}. It defines a kernel discrepancy between distributions \eqref{kernel discrepancy}, $\mg_K(\cdot,\cdot)$ is also called the \emph{Maximum Mean Discrepancy} (MMD) between $\mu$ and $\nu$ in $\SM^+(1)\cap\SM_K^{1/2}$, see \cite[Def.~10]{SejdinovicSGF2013}.

To define a metric on the whole $\SM^+(1)$, we need $P_\mu$ to be well defined and so that $P_\mu=P_\nu$ for $\mu$ and $\nu$ in $\SM^+(1)$ implies $\mu=\nu$.
This corresponds to the notion of \emph{characteristic kernel}, see \cite[Def.~6]{SriperumbudurGFSL2010},
which is closely connected to the following definitions.

\begin{defi}\label{D:ISPD}
A kernel $K$ is \emph{Integrally Strictly Positive Definite} (ISPD) on $\SM$ when $\SE_K(\nu)>0$ for any nonzero measure $\nu\in\SM$.
\end{defi}

\begin{defi}\label{D:CISPD}
A kernel $K$ is \emph{Conditionally Integrally Strictly Positive Definite} (CISPD) on $\SM$ when it is ISPD on $\SM(0)$; that is, when $\SE_K(\nu)>0$ for all nonzero signed measures $\nu\in\SM$ such that $\nu(\SX)=0$.
\end{defi}

An ISPD kernel is CISPD.
A bounded ISPD kernel is SPD and defines an RKHS. In \cite[Lemma~8]{SriperumbudurGFSL2010}, the authors show that a uniformly bounded kernel is characteristic if and only if it is CISPD. The proof is a direct consequence of the expression \eqref{kernel discrepancy} for the MMD $\mg_K(\mu,\nu)$. They also give (Corollary~4) a spectral interpretation of $\mg_K(\mu,\nu)$ and show that a translation-invariant kernel such that $K(\xb,\xb')=\Psi(\xb-\xb')$, with $\Psi$ a uniformly bounded continuous real-valued positive-definite function, satisfies, for any $\mu$ and $\nu$ in $\SM^+(1)$,
$$
\mg_K(\mu,\nu) = \left[\int_{\mathds{R}^d} \left|\phi_\mu(\mob)-\phi_\nu(\mob)\right|^2\, \dd\Lambda(\mob) \right]^{1/2} \,.
$$
Here,  $\phi_\mu$ and $\phi_\nu$ denote the characteristic functions of $\mu$ and $\nu$ respectively and $\Lambda$ is the spectral Borel measure on $\mathds{R}^d$, defined by
\be\label{Lambda}
\Psi(\xb)= \int_{\mathds{R}^d} e^{-i\xb^T\mob}\, \dd\Lambda(\mob)\,.
\ee
Using this spectral representation, they prove (Th.~9) that $K$ is characteristic if and only if the support of the associated $\Lambda$ coincides with $\mathds{R}^d$. For example, the sinc squared kernel $K(x,x')=\sin^2[\mt(x-x')]/(x-x')^2$, $\mt>0$, is SPD but is not characteristic (and therefore not CISPD) since the support of $\Lambda$ equals $[-2\mt,2\mt]$.
When $\mg_K(\mu,\delta_\xb)$ is well defined for all $\xb\in\SX$, with $\delta_\xb$ the Dirac delta measure at $\xb$ (and thus in particular when $K$ is characteristic), we may consider the empirical measure $\xi_{n,e}=\xi_{n,e}(\Xb_n)=(1/n) \sum_{i=1}^n \delta_{\xb_i}$ associated with a given design $\Xb_n=\{\xb_1,\ldots,\xb_n\}$, and $\mg_K(\mu,\xi_{n,e})$ of \eqref{pseudoM-PD} gives the worst-case integration error for $\xi_{n,e}$ when $f$ has norm one in $\SH_K$; see Section~\ref{S:empirical}.

Typical examples of uniformly bounded ISPD, and therefore characteristic, kernels are the squared exponential kernel $K_t(\xb,\xb')=\exp(-t\,\|\xb-\xb'\|^2)$, $t>0$, and the isotropic Mat\'ern kernels, in particular
\be\label{K32}
K_{3/2,\mt}(\xb,\xb')=(1+\sqrt{3}\mt\,\|\xb-\xb'\|)\, \exp(-\sqrt{3}\mt\,\|\xb-\xb'\|) \ \mbox{ (Mat\'ern 3/2)}\,,
\ee
and $K_{5/2,\mt}(\xb,\xb')=[1+\sqrt{5}\mt\, \|\xb-\xb'\|+ 5\mt^2\, \|\xb-\xb'\|^2/3]\,\exp(-\sqrt{5}\mt\, \|\xb-\xb'\|)$ (Mat\'ern 5/2), see, e.g., \cite{Stein99}.
Two other important examples are given hereafter.

\begin{example}[Generalized multiquadric kernel]\label{Ex:integral-Gaussian-kernels}
The sum of ISPD kernels is ISPD. Since the squared exponential kernel $K_t(\xb,\xb')$ is ISPD for any $t>0$, the integrated kernel obtained by setting a probability distribution on $t$ is ISPD too. One may thus consider
$K(\xb,\xb')= \int K_t(\xb,\xb')\, \dd\pi(t)$ for $\pi$ bounded and non decreasing on $[0,+\infty)$, which generates the class of continuous isotropic autocovariance functions in arbitrary dimension, see \cite{Schoenberg38} and \cite[p.~44]{Stein99}. In particular, for any $\me>0$ and $s>0$, we obtain
$$
K(\xb,\xb')= \int_0^{+\infty} K_t(\xb,\xb')\, t^{s/2-1}\, \exp(-\me\,t) \, \dd t = \frac{\Gamma(s/2)}{(\|\xb-\xb'\|^2+\me)^{s/2}} \,,
$$
showing that the generalized multiquadric kernel
\be\label{Kes}
K_{s,\me}(\xb,\xb')=(\|\xb-\xb'\|^2+\me)^{-s/2}\,, \ \me > 0\,, \ s>0 \,,
\ee
is ISPD, see also \cite[Sect.~3.2]{SriperumbudurGFSL2010}.
\fin
\end{example}

\begin{example}[distance-induced kernels]\label{Ex:Riesz<0}
Consider the kernels defined by
\be\label{Riesz<0}
K^{(s)}(\xb,\xb') = - \|\xb-\xb'\|^{s} \,, \ s>0 \,,
\ee
which are CISPD for $s\in(0,2)$ \cite{SzekelyR2013}, and the related distance-induced kernels
$$
K'^{(s)}(\xb,\xb')=\|\xb\|^{s}+\|\xb'\|^{s}-\|\xb-\xb'\|^{s} \,, \ s>0 \,.
$$
Note that $\SE_{K'^{(s)}}(\mu)=\SE_{K^{(s)}}(\mu)$ when $\mu(\SX)=0$; in \cite{SzekelyR2013} $\SE_{K'^{(s)}}$ is called energy distance for $s=1$ and generalized energy distance for general $s\in(0,2]$.
For $s>0$, the set $\SM_{K'^{(s)}}$ contains all signed measures $\mu$ such that $\int_\SX \|\xb-\xb_0\|^{s}\, \dd|\mu|(\xb)<+\infty$ for some $\xb_0\in\SX$. This result is a direct consequence of the triangular inequality when $s\in(0,1]$; for $s>1$ it follows from considerations involving semimetrics generated by kernels, see \cite[Remark~21]{SejdinovicSGF2013}.
$K'^{(s)}$ is CISPD for $s\in(0,2)$ ($K'^{(s)}/2$ corresponds to the covariance function of the fractional Brownian motion), but is not SPD (one has in particular, $K'^{(s)}(\0b,\0b)=0$); $K'^{(2)}$ is not CISPD since $\SE_{K'^{(2)}}(\mu)=[\int_\SX \xb^T\,\dd\mu(\xb)][\int_\SX \xb\,\dd\mu(\xb)]$, $\mu\in\SM$.
$K(x,x')=1-K^{(1)}(x,x')=1-|x-x'|$ is ISPD for $\SX=[0,1]$.
\fin
\end{example}

\subsection{Energy and potentials}\label{S:EP}

In this section we extend the considerations of previous section to signed measures and kernels which may have singularity on the diagonal. Definitions~\ref{D:ISPD} and \ref{D:CISPD} extend to singular kernels, with Riesz kernels as typical examples.

\begin{example}[Riesz kernels]\label{Ex:Riesz}
These fundamental kernels of potential theory are defined by
\be\label{Riesz}
K_{(s)}(\xb,\xb')=\|\xb-\xb'\|^{-s} \,, \ s>0 \,, \mbox{ and } K_{(0)}(\xb,\xb')=-\log \|\xb-\xb'\| \,,
\ee
with $\xb,\xb'\in\SX\subset\mathds{R}^d$ and $\|\cdot\|$ the Euclidean norm. When $s\geq d$, $\SE_{K_{(s)}}(\mu)$ is infinite for any nonzero signed measure, but for $s\in(0,d)$ $K_{(s)}$ is ISPD.
Since the logarithmic kernel $K_{(0)}(\xb,\xb')$ has  singularity at zero and tends to $-\infty$ when $\|\xb-\xb'\|$ tends to $+\infty$, it will only be considered for~$\SX$ compact; $K_{(0)}$ is CISPD, see \cite[p.~80]{Landkof1972}.
\fin
\end{example}

\vsp
Consider again $\SE_K(\mu)$ given by \eqref{overline K}, with $\mu\in\SM$. In potential theory, this quantity is called the energy of the signed measure $\mu$ for the kernel $K$. Denote
$$
\SM_K = \left\{\nu\in\SM:  |\SE_K(\nu)| < +\infty \right\} \,.
$$
In the following, we shall only consider kernels that are at least CISPD.
When $K$ is ISPD, $\SE_K(\mu)$ is positive for any nonzero $\mu\in\SM$, but when $K$ is only CISPD, $\SE_K(\mu)$ can be negative; this is the reason for the presence of absolute value in the definition of $\SM_K$. Note that $\SM_K$ is the set of measures such that $\SE_K(\mu^+)$, $\SE_K(\mu^-)$ and $\SE_K(\mu^+,\mu^-)=\int_{\SX^2} K(\xb,\xb')\,\dd\mu^+(\xb)\dd\mu^-(\xb')$ are all finite, with $\mu^+$ and $\mu^-$ denoting the positive and negative parts of the Hahn-Jordan decomposition $\mu=\mu^+-\mu^-$ of $\mu$, see \cite[Sect.~2.1]{Fuglede60}. Also note that when $K$ is bounded and defines an RKHS, $\SM_K^\ma \subset \SM_K$ for any $\ma\geq 1/2$, see \eqref{M_K^a} and \eqref{CS:M1/2}; when $K$ is uniformly bounded, $\SM_K=\SM$.

For any $\mu\in\SM_K$, $P_\mu(\xb)$ given by \eqref{kh} is called the potential at $\xb$ associated with $\SE_K(\mu)$. It is well-defined, with values in $\mathds{R}\cup\{-\infty,+\infty\}$, when $P_{\mu^+}(\xb)$ and $P_{\mu^-}(\xb)$ are not both infinite. Also, $P_\mu(\xb)$ is finite for $\mu$-almost any $\xb$, even if $K$ is singular, when $\mu\in\SM^+(1)\cap\SM_K^{1/2}$.

When $K$ is ISPD, we can still define MMD through \eqref{kernel discrepancy},
\be\label{kernel discrepancyB}
\mg_K(\mu,\nu) = \SE_K^{1/2}(\nu-\mu) \,,
\ee
since $\SE_K(\nu-\mu)$ is nonnegative whenever defined. The set $\SM_K$ forms a pre-Hilbert space, with scalar product the mutual energy $\SE_K(\mu,\nu)=\int_{\SX^2} K(\xb,\xb')\, \dd\mu(\xb)\dd\nu(\xb')$ and norm $\SE_K^{1/2}(\mu)$.
Denote by $\SP_K$ the linear space of potential fields  $P_\mu(\cdot)$, $\mu\in\SM_K$; when $K$ defines an RKHS $\SH_K$, $\|P_\mu\|_{\SH_K}=\SE_K^{1/2}(\mu)$ so that $\SP_K\subset\SH_K$, and $\SP_K$ is dense in $\SH_K$. For $\SP_K$ to contain all functions $K_\xb(\cdot)=K(\cdot,\xb)$, $\xb\in\SX$, we need $\delta_\xb\in\SM_K$ for all $\xb$, which requires $K(\xb,\xb)<\infty$ for all $\xb\in\SX$.

For $\mu,\nu\in\SM_K$, $\SE_K(\mu,\nu)$ defines a scalar product $\langle P_\mu,P_\nu \rangle_{\SP_K}$ on $\SP_K$, with $\mg_K(\mu,\nu)=\|P_\mu-P_\nu\|_{\SP_K}$. Similarly to Section~\ref{S:MMD}, we obtain
\be
\mg_K(\mu,\nu) &=& \sup_{\xi\in\SM_K,\, \SE_K(\xi)=1} \left| \int_{\SX^2} K(\xb,\xb')\, \dd\xi(\xb)\dd(\mu-\nu)(\xb') \right| \label{MMD-bis} \\
&=& \sup_{\|h\|_{\SP_K}\leq 1} \left|I_\mu(h) - I_{\nu}(h) \right| \,; \nonumber
\ee
that is, a result that extends \eqref{pseudoM-PD} to general ISPD kernels.
If $K$ is only CISPD, we can also define $\mg_K(\mu,\nu)$ in the same way when considering measures $\mu,\nu\in\SM(1)$; we then define $\SP_K$ as the linear space of potentials fields $P_\mu(\cdot)$, $\mu\in\SM_K\cap\SM(0)$, and in \eqref{MMD-bis} we restrict $\xi$  to be in $\SM(0)$.

When $K$ is singular, there always exists $\nu$ in $\SM_K$ such that $P_\nu(\xb_0)=+ \infty$ for some $\xb_0$. Consider for example Riesz kernel $K_{(s)}(\xb,\xb')$ with $s\in(0,d)$; $\SM_K$ contains in particular all signed measures with compact support $\SSS(\mu)$ whose potential $P_\mu(\xb)$ is bounded on $\SSS(\mu)$, see \cite[p.~81]{Landkof1972}. Take $\nu$ as the measure with density $c/\|\xb-\xb_0\|^{s-d}$ on~$\SX$, with $\xb_0\in\SX$; we have $\SE_{K_{(s)}}(\nu)<\infty$ for~$\SX$ compact, but $P_\nu(\xb_0)=+\infty$. As a consequence, as noted in \cite{DamelinHRZ2010}, singular kernels have little interest for integration. Indeed, take $\mu,\nu\in\SM_K$ and $h=P_\nu\in\SP_K$, then $|I_\mu(h)| \leq \|h\|_{\SP_K}\,\SE_K^{1/2}(\mu)=\SE_K^{1/2}(\nu)\SE_K^{1/2}(\mu)<\infty$ whereas $|I_{\xi_n}(h)|$ may be infinite for some discrete approximation $\xi_n$ of $\mu$ as $h$ can be infinite at some points. Singular kernels may nevertheless be used for the construction of space-filling designs, see for instance the example in Section~\ref{Ex:4.1}, and this is our motivation for considering them in the following.

The key difficulty with singular kernels is the fact that delta measures do not belong to $\SM_K$.
An expedient solution to circumvent the problem is replace a singular kernel with a bounded surrogate. For instance, in space-filling design we may replace Riesz kernel $K_{(s)}$, $s>0$, by a generalized inverse multiquadric kernel $K_{s,\me}$ given by \eqref{Kes}, and consider the limiting behaviour of the designs obtained when $\me \ra 0$, see Section~\ref{S:empirical}.

\subsection{Minimum energy and equilibrium measures}\label{S:minE-measures}

In this section, we show that there exist strong connections between results in potential theory and optimal design theory, where one minimizes a convex functional of $\mu\in\SM^+(1)$, with the particularity that here the functional is quadratic.

\subsubsection{ISPD kernels and convexity of $\SE_K(\cdot)$}

\begin{lem}\label{P:convexity1} \mbox{}
$K$ is ISPD if and only if $\SM_K$ is convex and $\SE_K(\cdot)$ is strictly convex on $\SM_K$.
\end{lem}

\noindent{\em Proof.} \mbox{}
For any $K$, any $\mu$ and $\nu$ in $\SM_K$ and any $\ma\in[0,1]$, direct calculation gives
\be \label{convexity1}
(1-\ma)\,\SE_K(\mu)+\ma\,\SE_K(\nu) - \SE_K[(1-\ma)\mu+\ma\nu] = \ma(1-\ma)\, \SE_K(\nu-\mu) \,.
\ee

Assume that $K$ is ISPD. For any $\mu$ and $\nu$ in $\SM_K$, the mutual energy $\SE_K(\mu,\nu)$ satisfies $|\SE_K(\mu,\nu)|\leq \sqrt{\SE_K(\mu)\SE_K(\nu)} < +\infty$. Therefore, $\SE_K(\mu-\nu)=\SE_K(\mu)+\SE_K(\nu)-2\,\SE_K(\mu,\nu)$ is finite and \eqref{convexity1} implies that $\SE_K[(1-\ma)\mu+\ma\nu]$ is finite, showing that $\SM_K$ is convex. Since $K$ is ISPD, $\SE_K(\nu-\mu)>0$ for $\mu,\nu\in\SM$, $\nu\neq \mu$, and \eqref{convexity1} implies that $\SE_K(\cdot)$ is strictly convex on $\SM_K$.

Conversely, assume that $\SM_K$ is convex and $\SE_K(\cdot)$ is strictly convex on $\SM_K$. Any $\xi\in\SM_K$ can be written as $\xi=\nu-\mu$ with, for instance, $\nu=2\xi$ and $\mu=\xi$, both in $\SM_K$. If $\SE_K(\cdot)$ is strictly convex on $\SM_K$, \eqref{convexity1} with $\ma\in(0,1)$ implies that $\SE_K(\xi)>0$ when $\nu\neq \mu$, that is, when $\xi\neq 0$. Therefore, $K$ is ISPD.
\carre

\vsp
Lemma~\ref{P:convexity1} also applies to singular kernels. The lemma below concerns CISPD kernels, which are assumed to be uniformly bounded.

\begin{lem}\label{P:convexity2} \mbox{}
Assume that $K$ is uniformly bounded. Then, $K$ is CISPD if and only if $\SE_K(\cdot)$ is strictly convex on $\SM(1)$.
\end{lem}

\noindent{\em Proof.} \mbox{}
Since $K$ is uniformly bounded, $\SM_K=\SM$. Assume that $K$ is CISPD. Then, $\SE_K(\nu-\mu)>0$ for any $\mu\neq\nu\in\SM(1)$, and \eqref{convexity1} implies that $\SE_K(\cdot)$ is strictly convex on $\SM(1)$.

Assume now that $\SE_K(\cdot)$ is strictly convex on $\SM(1)$. Take any non-zero signed measure $\xi$ in $\SM(0)$ and consider the Hahn-Jordan decomposition $\xi=\xi^+-\xi^-$, with $\xi^+(\SX)=\xi^-(\SX)=c>0$. Denote $\nu=\xi^+/c$, $\mu=\xi^-/c$, with $\nu$ and $\mu$ in $\SM^+(1)$ ($\nu$ and $\mu$ are in $\SM_K$ since $K$ is uniformly bounded). Then, for any $\ma\in(0,1)$, \eqref{convexity1} and the strict convexity of $\SE_K(\cdot)$ on $\SM(1)$ gives
$\SE_K(\xi) = c^2\, \SE_K(\nu-\mu) >0$.
\carre

\vsp
Note that one may replace $\SM(1)$ by $\SM^+(1)$, or by any $\SM(\mg)$ with $\mg\neq 0$, in Lemma~\ref{P:convexity2}.

\subsubsection{Minimum-energy probability measures}

In the remaining part of Section~\ref{S:minE-measures}, we assume that $K$ is such that $\SE_K(\cdot)$ is strictly convex on $\SM^+(1)\cap\SM_K$ and $\SM(1)\cap\SM_K$, which is true under the conditions of Lemma~\ref{P:convexity1} or Lemma~\ref{P:convexity2}.

For $\mu,\nu\in\SM_K$, denote by $F_K(\mu;\nu)$ the directional derivative of $\SE_K(\cdot)$ at $\mu$ in the direction $\nu$,
$$
F_K(\mu;\nu) = \lim_{\ma\ra 0^+} \frac{\SE_K[(1-\ma)\mu+\ma\nu]-\SE_K(\mu)}{\ma} \,.
$$
Straightforward calculation gives
\be \label{Dir-der}
F_K(\mu;\nu) = 2 \left[ \int_{\SX^2} K(\xb,\xb')\, \dd\nu(\xb)\dd\mu(\xb') - \SE_K(\mu) \right] \,.
\ee
In particular, for any $\xb\in\SX$, the potential $P_\mu(\xb)$ associated with $\mu$ at $\xb$ satisfies
$$
P_\mu(\xb) = \frac12 \, F_K(\mu;\delta_\xb) + \SE_K(\mu) \,.
$$

\begin{remark}[Bregman divergence and Jensen difference]\label{R:Bregman}
The strict convexity of $\SE_K(\cdot)$ implies that $\SE_K(\nu) \geq \SE_K(\mu)+F_K(\mu,\nu)$ for any $\mu,\nu\in\SM_K$, with equality if and only if $\nu=\mu$. This can be used to define a Bregman divergence between measures in $\SM_K$ (and thus between probability measures in $\SM^+(1)\cap\SM_K$), as
$$
B_K(\mu,\nu) = \SE_K(\nu)-[\SE_K(\mu)+F_K(\mu,\nu)] \,;
$$
see \cite{RaoN85}.
Direct calculation gives $B_K(\mu,\nu)=\SE_K(\nu-\mu)$ (with therefore $B_K(\mu,\nu)=B_K(\nu,\mu)$), providing another interpretation for the MMD $\mg_K(\mu,\nu)$, see \eqref{kernel discrepancyB}.

The squared MMD is also proportional to dissimilarity coefficient, or Jensen difference, $\Delta_J(\mu,\nu)=(1/2)[\SE_K(\mu)+\SE_K(\nu)]-\SE_K[\mu+\nu)/2]$ of \cite{Rao1982a}; indeed, direct calculation gives $\mg_K^2(\mu,\nu)=\SE_K(\nu-\mu)=4\,\Delta_J(\mu,\nu)$.
\fin
\end{remark}

\vsp
Since $\SE_K(\cdot)$ is strictly convex on $\SM^+(1)$, there exists a unique minimum-energy probability measure. The measure $\mu_K^+\in\SM^+(1)$ is the minimum-energy measure if and only if $F_K(\mu_K^+;\nu)\geq 0$ for all $\nu\in\SM^+(1)$, or equivalently, since $\nu$ is a probability measure, if and only if $F_K(\mu_K^+;\delta_\xb)\geq 0$ for all $\xb\in\SX$. We thus obtain the following property, called equivalence theorem in the optimal-design literature.

\begin{theo}\label{P:ET-P}
When $\SE_K(\cdot)$ is strictly convex on $\SM^+(1)\cap\SM_K$, $\mu_K^+\in\SM^+(1)$ is the minimum-energy probability measure on~$\SX$ if and only if
$$
\forall\,  \xb\in\SX\,, \ P_{\mu_K^+}(\xb) \geq \SE_K(\mu_K^+) \,.
$$
\end{theo}

Note that, by construction, $\int_\SX P_{\mu_K^+}(\xb) \,\dd\mu_K^+(\xb) = \SE_K(\mu_K^+)$, implying $P_{\mu_K^+}(\xb)=\SE_K(\mu_K^+)$ on the support of $\mu_K^+$.
The quantity $C_K^+=[\inf_{\mu\in\SM^+(1)} \SE_K(\mu)]^{-1}$, with $K$ an ISPD kernel, is called the \emph{capacity} of~$\SX$ in potential theory; note that $C_K^+ \geq 0$.
The minimizing measure $\mu_K^+\in\SM^+(1)$ is called the \emph{equilibrium measure} of~$\SX$ ($\mu_K^+$ is sometimes renormalized into $C_K^+\,\mu_K^+$, see \cite[p.~138]{Landkof1972}). Theorem~\ref{P:ET-P} thus gives a necessary and sufficient condition for a probability measure $\mu$ to be the equilibrium measure of~$\SX$.

\begin{example}[Continuation of Example~\ref{Ex:Riesz<0}]\label{Ex:Riesz<0,d>=2} 
Properties of minimum-energy probability measures $\mu^+=\mu_{K^{(s)}}^+$ for $K^{(s)}$ given by \eqref{Riesz<0} with~$\SX$ a compact subset of $\mathds{R}^d$, $d\geq 2$,  are investigated in \cite{Bjorck56} and \cite{PWZ2016-SP}. The mass of $\mu^+$ is concentrated on the boundary of~$\SX$, and its support only comprises extreme points of the convex hull of~$\SX$ when $s>1$; for $0<s<2$, $\mu^+$ is unique; it is supported on no more than $d+1$ points when $s>2$.

Take $\SX=\SB_d(0,1)$. For symmetry reasons, $\mu^+$ for $0<s<2$ is uniform on the unit sphere $\SSp_d(\0b,1)$ and
$$
\SE_{K^{(s)}}(\mu^+) = -\int_{\SX^2} \|\xb-\xb'\|^{s} \, \dd\mu^+(\xb)\dd\mu^+(\xb') = -\int_{\SX} \|\xb_0-\xb'\|^{s} \, \dd\mu^+(\xb')\,,
$$
where $\xb_0=(1,0,\ldots,0)^T$. Denote by $\psi_d(\cdot)$ the density of the first component $t=x'_1$ of $\xb'=(x'_1,\ldots,x'_d)^T$. We obtain
$\psi_d(t) = (d-1)\, V_{d-1}\,(1-t^2)^{(d-3)/2}/(d\, V_d)$ and
$$
\SE_{K^{(s)}}(\mu^+) = -\int_{-1}^1 [(1-t)^2+1-t^2]^{s/2}\, \psi_d(t)\dd t = -\frac {{2}^{d-q-2}\Gamma  (d/2) \Gamma[(d+s-1)/2]}{\sqrt {\pi }\Gamma(d+s/2-1)} \,.
$$
In particular, $\SE_{K^{(1)}}(\mu^+)=-4/\pi$ when $d=2$ and is a decreasing function of $d$. When $s=2$, the uniform distribution on the unit sphere is also optimal, and the minimum energy equals $-2$ for all $d\geq 1$, but $\mu^+$ is not unique and the measure allocating equal weight $1/(d+1)$ at each of the $d+1$ vertices of a $d$ regular simplex with vertices on the unit sphere is optimal too.
\fin
\end{example}

\begin{example}[Continuation of Example~\ref{Ex:Riesz}]\label{Ex:Riesz-B_2}
Consider Riesz kernels $K_{(s)}$, see \eqref{Riesz}, for $\SX=\SB_d(\0b,1)$, the closed unit ball in $\mathds{R}^d$.
When $s\geq d$, $\SE_{K_{(s)}}(\nu)$ is infinite for any non-zero $\nu\in\SM$, but for $0<s<d$ there exists a minimum-energy probability measure $\mu^+=\mu_{K_{(s)}}^+$. When $d>2$ and $s\in(0,d-2]$, $\mu^+$ is uniform on the unit sphere $\SSp_d(\0b,1)$ (the boundary of~$\SX$); the potential at all interior points satisfies $P_{\mu^+}(\xb) \geq \SE_{K_{(s)}}(\mu^+)$ with strict inequality when $s\in(0,d-2)$. When $s\in(d-2,d)$, $\mu^+$ has a density $\varphi_s(\cdot)$ in $\SB_d(\0b,1)$,
$$
\varphi_s(\xb)=  \frac{\pi^{-d/2}\,\Gamma(1+s/2)}{\Gamma[1-(d-s)/2]}\; \frac{1}{(1-\|\xb\|^2)^{(d-s)/2}} \,,
$$
and the potential $P_{\mu^+}(\cdot)$ is constant in $\SB_d(\0b,1)$, see, e.g., \cite[p.~163]{Landkof1972}.

When $d\leq 2$ and $s=0$, $\mu^+$ has a density in $\SB_2(\0b,1)$ and $P_{\mu^+}(\cdot)=\SE_{K_{(0)}}(\mu^+)$ in $\SB_2(\0b,1)$. In particular, for $d=1$, $\mu^+$ has the arsine density $1/(\pi\sqrt{1-x^2})$ in $[-1,1]$ with potential $P_{\mu^+}(x) = \log(2)-\log(||x|+\sqrt{x^2-1}|)$, $x\in\mathds{R}$ (and $P_{\mu^+}(x)=\log(2)$ for $x\in[-1,1]$).

The energy $\SE_{K_{(s)}}$ is infinite for empirical measures associated with $n$-point designs $\Xb_n$. One may nevertheless consider the ``physical'' energy
\be\label{physical energy}
\widetilde\SE_{K_{(s)}}(\Xb_n)=[2/n(n-1)] \sum_{1\leq i<j\leq n} \|\xb_i-\xb_j\|^{-s}
\ee
($\widetilde\SE_{K_{(s)}}(\Xb_n)=-[2/n(n-1)] \sum_{1\leq i<j\leq n} \log\|\xb_i-\xb_j\|$ when $s=0$), which is finite provided that all $\xb_i$ are distinct, see \cite{DamelinHRZ2010}. An $n$-point set $\Xb_n^*$ minimizing $\widetilde\SE_{K_{(s)}}(\Xb_n)$ is called a set of Fekete points, and the limit $\lim_{n\ra\infty} \widetilde\SE_{K_{(s)}}^{-1}(\Xb_n^*)$ exists and is called the transfinite diameter of~$\SX$. A major result in potential theory, see, e.g., \cite{HardinS2004}, is that the transfinite diameter coincides with the capacity $C_{K_{(s)}}^+$ of~$\SX$. If $C_{K_{(s)}}^+>0$, then $\mu_{K_{(s)}}^+$ is the weak limit of a sequence of empirical probability measures associated with Fekete points in $\Xb_n^*$. In the example considered, $\widetilde\SE_{K_{(s)}}(\Xb_n^*)$ tends to infinity when $s\geq d$, but any sequence of Fekete points is asymptotically uniformly distributed in~$\SX$; $\widetilde\SE_{K_{(s)}}(\Xb_n^*)$ grows like $n^{s/d-1}$ for $s>d$ (and like $\log n$ for $s=d$).
\fin
\end{example}

\begin{remark}[Stein variational gradient descent and energy minimization]
Variational inference using smooth transform based on kernelized Stein discrepancy provides a gradient descent method for the approximation of a target distribution; see \cite{LiuW2016} and the references therein. The fact that the construction does not require knowledge of the normalizing constant of the target distribution makes the method particularly attractive for approximating a posterior distribution in Bayesian inference. Direct calculation shows that when the kernel is translation invariant and the target distribution is uniform, then the method corresponds to steepest descent for the minimization of $\SE_K(\xi_{n,e})$; that is, at iteration $k$ each design point $\xb_i^{(k)}$ is updated into
$$
\xb_i^{(k+1)}=\xb_i^{(k)}+ \mg\, \sum_{i<j} \frac{\mp K(\xb,\xb_j^{(k)})}{\mp \xb}\bigg|_{\xb=\xb_i^{(k)}}
$$
for some $\mg>0$.
The construction of space-filling design through energy minimization has already been considered in the literature. For instance, it is suggested in \cite{AudzeE77} to construct designs in a compact subset~$\SX$ of $\mathds{R}^d$ by minimizing $\widetilde\SE_{K_{(2)}}(\Xb_n)$ given by \eqref{physical energy} (note that for $d\geq 3$ design points constructed in this way are not asymptotically uniformly distributed in~$\SX$). This approach tends to push points to the border of~$\SX$, similarly to the maximization of the packing radius $\PR(\Xb_n)$ defined by \eqref{PR}. This is generally not desirable, especially when $d$ is large.
\fin
\end{remark}

\subsubsection{Minimum-energy signed measures}

The situation is slightly different from that in previous section when we consider measures in $\SM(1)$.
In that case, $\mu_K^*$ is the minimum-energy measure in $\SM(1)$ if and only if $F_K(\mu_K^*;\nu)= 0$ for all $\nu\in\SM(1)$, this condition being equivalent
to $F_K(\mu_K^*;\delta_\xb) = 0$ for all $\xb\in\SX$. We thus obtain the following property.

\begin{theo}\label{P:ET-M(1)}
When $\SE_K(\cdot)$ is strictly convex on $\SM(1)\cap\SM_K$, $\mu_K^*\in\SM(1)$ is the minimum-energy signed measure with total mass one on~$\SX$ if and only if
\be\label{ET-M(1)}
\forall\, \xb\in\SX\,, \ P_{\mu_K^*}(\xb) = \SE_K(\mu_K^*) \,.
\ee
\end{theo}

If we define now a \emph{signed equilibrium measure} on~$\SX$ as a measure $\mu\in\SM(1)$ such that $P_{\mu}(\xb)$ is constant on~$\SX$, from the definition of $P_\mu(\cdot)$, when such a measure exists it necessarily satisfies the condition of Theorem~\ref{P:ET-M(1)} and therefore coincides with $\mu_K^*$. Similarly to the case where one considers probability measures in $\SM^+(1)$, we can define the (generalized) capacity of~$\SX$ for measures in $\SM(1)$ as $C_K^*=[\inf_{\mu\in\SM^(1)} \SE_K(\mu)]^{-1}$, with $C_K^*=1/\SE_K(\mu_K^*)$ when $\mu_K^*$ exists, see \cite[p.~824]{DamelinHRZ2010} (note that $C_K^*$ may be negative).
However, $\mu_K^*$ may not exist. Notice in particular that $\SM(1)$ is not vaguely compact, contrarily to $\SM^+(1)$ (and for Riesz kernels \eqref{Riesz} with $s<d-1$, $\SM_{K_{(s)}}$ is not complete contrarily to $\SM_{K_{(s)}}\cap\SM^+$  \cite[Th.~1.19]{Landkof1972}).

\begin{example}[Continuation of Examples~\ref{Ex:Riesz<0} and \ref{Ex:Riesz<0,d>=2}]\label{Ex:Riesz<0,d=1} Take $K(x,x')=K^{(s)}(x,x')=-|x-x'|^{s}$ on $\SX=[0,1]$, $s\in(0,2)$, see \eqref{Riesz<0}. $K$ is CISPD, and there exists a unique minimum-energy probability measure $\mu^+=\mu_{K^{(s)}}^+$ in $\SM^+(1)$. On the other hand, below we show that minimum-energy signed measures in $\SM(1)$ do not belong to $\SM^+(1)$ when $s\in(1,2)$ and that there is no minimum-energy signed measure in $\SM(1)$ when $s \geq 2$.

When $s\in(0,1)$, $\mu^+$ has a density $\varphi^{(s)}(\cdot)$ with respect to the Lebesgue measure on $[0,1]$,
$$
\varphi^{(s)}(x)=  \frac{\Gamma[1-s/2]}{2^s\, \sqrt{\pi}\,\Gamma[(1-s)/2]}\; \frac{1}{[x(1-x)]^{(1+s)/2}} \,,
$$
and $P_{\mu^+}(x)=\SE(\mu^+)=-\sqrt{\pi}\,\Gamma(1-s/2)/\{2^s\,\Gamma[(1-s)/2]\,\cos(\pi s/2)\}$ for all $x\in\SX$ (and $\SE(\mu^+) \ra -1/2$ as $s\ra 1^-$). The fact that $P_{\mu^+}(x)=\SE(\mu^+)$ for all $x\in\SX$ indicates that $\mu^+$ is the minimum-energy signed measure with total mass one when $s\in(0,1]$.

When $s\in[1,2)$, $\mu^+=(\delta_0 + \delta_1)/2$; the associated potential is $P_{\mu^+}(x)=-(|x|^{s}+|1-x|^{s})/2 \geq \SE(\mu^+)=-1/2$, $x\in\SX$ (note that $P_{\mu^+}(x)=-1/2$ for all $x\in\SX$ when $s=1$).

Consider now the signed measure $\mu_w=[(1+w)/2](\delta_0+\delta_1)-w \delta_{1/2}$, $w>0$, so that $\mu_w(\SX)=1$ (i.e., $\mu_w\in\SM(1)$). Direct calculation gives
$\SE_{K^{(s)}}(\mu_w)=-(1+w)(1+w-2^{2-s}w)$, which is minimum for $w=w_*(s)=(1-2^{1-s})/(2^{2-s}-1)$ when $s<2$, with $\SE_{K^{(s)}}(\mu_{w_*(s)})=2(1-2^{2-s})/(4-2^{s})^2$. For $s\in(1,2)$ we get $\SE_{K^{(s)}}(\mu_{w_*(s)})< \SE(\mu^+)=-1/2$, and there exist signed measures in $\SM(1)$ such that $\SE_{K^{(s)}}(\mu_w)<\SE(\mu^+)$. Therefore, minimum-energy signed measures with total mass one are not probability measures. For $s\geq 2$, $\lim_{w\ra+\infty}\SE_{K^{(s)}}(\mu_w)=-\infty$, and there is no minimum-energy signed measure; in particular, $\SE_{K^{(s)}}(\mu_w)=-(w+1)/2$ for $s=2$.
\fin
\end{example}

\begin{example}[Continuation of Examples~\ref{Ex:Riesz} and \ref{Ex:Riesz-B_2}]\label{Ex:Riesz3}
Consider Riesz kernels $K_{(s)}$, see \eqref{Riesz}, for $\SX=\SB_d(\0b,1)$, $d>2$ and $s\in(0,d-2)$; the minimum-energy probability measure $\mu^+$ is then uniform on the unit sphere $\SSp_d(\0b,1)$ and the potential at all interior points satisfies $P_{\mu^+}(\xb) > \SE_{K_{(s)}}(\mu^+)$. Consider the signed measure $\mu_w=(1+w) \mu^+ -w \mu^{(r)}$, with $\mu^{(r)}$ uniform on the sphere $\SSp_d(\0b,r)$ with radius $r\in(0,1)$. Calculations similar to those in the proof of \cite[Th.~1.32]{Landkof1972} show that
$\SE_{K_{(s)}}(\mu_w)<\SE_{K_{(s)}}(\mu^+)$ for $w$ small enough, indicating that $\mu^+$ is not the minimum-energy signed measure with total mass one.
\fin
\end{example}

\subsubsection{When minimum-energy signed measures are probability measures}\label{S:subharmonic}

\begin{theo}\label{Th:MAIN}
Assume that $K$ is ISPD and translation invariant, with $K(\xb,\xb')=\Psi(\xb-\xb')$ and $\Psi$ continuous,
twice differentiable except at the origin, with $\Psi(\0b)<\infty$ and Laplacian $\Delta_\Psi(\xb)=\sum_{i=1}^d \mp^2 \Psi(\xb)/\mp x_i^2 \geq 0$, $\xb\neq \0b$.
Then there exists a unique minimum-energy signed measure $\mu_K^*$ in $\SM(1)$, and $\mu_K^*$ is a probability measure.
\end{theo}

\noindent{\em Proof.} \mbox{}
The conditions of Theorem~\ref{P:ET-P} are satisfied, and there exists a unique minimum-energy probability measure $\mu^+$ such that
$P_{\mu^+}(\xb) \geq \SE_K(\mu^+)$ for all $\xb\in\SX$. It also satisfies $P_{\mu^+}(\xb) = \SE_K(\mu^+)$ on the support of $\mu^+$.
On the other hand, the conditions on $K$ imply that for any $\mu$ in $\SM^+(1)$, $P_\mu(\cdot)$ is subharmonic outside the support of $\mu$, see, e.g., \cite[Sect.~I.2]{Landkof1972}. The first maximum principle of potential theory thus holds \cite[Th.~1.10]{Landkof1972}: $P_\mu(\xb)\leq c$ on the support of $\mu$ implies $P_\mu(\xb)\leq c$ everywhere. Applying this to $\mu^+$, we obtain that $P_{\mu^+}(\xb) \leq  \SE_K(\mu^+)$ everywhere; therefore, $P_{\mu^+}(\xb) = \SE_K(\mu^+)$ for all $\xb\in\SX$. Theorem~\ref{P:ET-M(1)} implies that $\mu^+$ is the minimum-energy signed measure with total mass one.
\carre

\vsp
The central argument for the proof of the property above is that $P_\mu(\cdot)$ is subharmonic outside the support of $\mu$ for any probability measure $\mu$ with finite energy. Weaker conditions than those in the theorem may be sufficient in particular situations, such as $\Psi(\xb-\xb')=\psi(\|\xb-\xb'\|)$ with $\psi(\cdot)$ convex on $(0,\infty)$, which generalizes a result of H\'ajek (1956), see also Section~\ref{S:continuous-BLUE}.

Another generalization is to consider CISPD kernels. For example, for the kernels $K^{(s)}$ of \eqref{Riesz<0}, we have $\Delta(-\|\xb\|^{s})=s(2-s-d)/\|\xb\|^{2-s}$, $\xb\neq\0b$. Potentials are superharmonic for $d\geq 2$. When $d=1$, they are superharmonic for $s\in[1,2)$;  they are subharmonic and satisfy the maximum principle for $s\in(0,1)$, see Example~\ref{Ex:Riesz<0,d=1}.

Extension of Theorem~\ref{Th:MAIN} to singular kernels requires advanced results from potential theory; see especially \cite{FugledeZ2018, Landkof1972}. In particular, for the Riesz kernels $K_{(s)}$ of \eqref{Riesz}, we have $\Delta(\|\xb\|^{-s})=s(s+2-d)/\|\xb\|^{s+2}$, $\xb\neq\0b$.
When $d>2$ and $s\in(0,d-2]$, $P_\mu$ can be proved to be superharmonic in $\mathds{R}^d$, and when
$s\in[d-2,d)$, $P_\mu$ can be proved to be subharmonic outside the support of $\mu$, $\mu^+$ being then the minimum-energy signed measure.
This is also true for the logarithmic kernel for $d\leq 2$, with $\Delta(-\log\|\xb\|)= (2-d)/\|\xb\|^2$, $\xb\neq\0b$.
Examples~\ref{Ex:Riesz-B_2} and \ref{Ex:Riesz3} give an illustration.

\subsection{Best Linear Unbiased Estimator (BLUE) of $\beta_0$}\label{S:BLUE}

\subsubsection{Continuous BLUE}\label{S:continuous-BLUE}

Consider again the situation of Section~\ref{S:BQ} where $\ms^2\,K$ corresponds to the covariance of a random field $Z_x$.
Suppose that we may observe $f(\cdot)$ over~$\SX$ in order to estimate $\beta_0$ in the regression (location) model with correlated errors \eqref{model1}. Any linear estimator of $\beta_0$ takes the general form
$$
\hat\beta_0 = \hat\beta_0(\xi) = \int_\SX f(\xb) \,\dd\xi(\xb) = I_\xi(f)
$$
for some $\xi\in\SM$, and $\hat\beta_0(\xi)$ is unbiased when $\xi\in\SM(1)$. Its variance is
$$
V_\xi=\Exx\{(\hat\beta_0(\xi)-\beta_0)^2\} = \ms^2\, \SE_K(\xi) \,;
$$
see \cite[Sect.~4.2]{Nather85}. The existence of a minimum-energy signed measure $\mu_K^*$ is then equivalent to the existence of the continuous BLUE $\hat\beta_0^*$ for $\beta_0$, with $\hat\beta_0^*=\hat\beta_0(\mu_K^*)$; the variance of $\hat\beta_0^*$ is proportional to the minimum energy $\SE_K(\mu_K^*)$, and Theorem~\ref{P:ET-M(1)} corresponds to Grenander's theorem \cite{Grenander50}.
Also, from that theorem, the existence of $\mu_K^*$ is equivalent to the existence of an equilibrium measure that yields a constant potential on~$\SX$. It can be related to a property of the generalized capacity $C_K^*$, as shown in the following theorem.

\begin{theo}\label{P:capacity-1}
When $K$ is ISPD, the constant function $1_\SX$ equal to 1 on~$\SX$ belongs to the space $\SP_K$ of potential fields if and only if there exists a minimum-energy signed measure $\mu_K^*\in\SM(1)$, with $\SE_K(\mu_K^*) \neq 0$. Moreover, the generalized capacity $C_K^*$ is finite and nonzero, and satisfies $\|1_\SX\|^2_{\SP_K}=C_K^*$.
\end{theo}

\noindent {\em Proof.}
Suppose that $1_\SX\in\SP_K$. There exists $\mu\in\SM_K$ such that $P_\mu=1_\SX$; that is, $P_\mu(\xb)=1$ for all $\xb\in\SX$. The definition of $P_\mu$ yields $\SE_K(\mu)=\mu(\SX)$, which is finite and strictly positive since $K$ is ISPD and $\mu\neq 0$. Denote $\mu'=\mu/\mu(\SX) \in\SM(1)$. We obtain $P_{\mu'}(\xb)=1/\mu(\SX)=\SE_K(\mu')>0$ for all $\xb\in\SX$. Theorem~\ref{P:ET-M(1)} implies that $\mu'$ is the minimum-energy measure $\mu_K^*$. Also, $C_K^* = 1/\SE_K(\mu')=\mu(\SX) \neq 0$, with $\|1_\SX\|^2_{\SP_K}=\SE_K(\mu)$, see Section~\ref{S:EP}.

Suppose now that there exists a minimum-energy signed measure $\mu_K^*\in\SM(1)$ with $\SE_K(\mu_K^*) \neq 0$. Theorem~\ref{P:ET-M(1)} implies that
$P_{\mu_K^*}(\xb) = \SE_K(\mu_K^*)$ for all $\xb\in\SX$. For $\mu=\mu_K^*/\SE_K(\mu_K^*)$, we get $P_\mu(\xb)=1$ for all $\xb\in\SX$, and $\|1_\SX\|^2_{\SP_K}=\SE_K(\mu)=1/\SE_K(\mu_K^*)$.
\carre

\vsp
Under the conditions of Theorem~\ref{Th:MAIN}, the BLUE exists, $\hat\beta_0^*=\hat\beta_0(\mu_K^+)$, with $\mu_K^+$ the minimum-energy probability measure, and its variance equals $\ms^2\SE_K(\mu_K^+)$. This is also true when $K(\xb,\xb')=\psi(\|\xb-\xb'\|)$ with $\psi(\cdot)$ convex on $(0,\infty)$. For $d=1$, this property was known to H\'ajek (1956), see \cite[p.~56]{Nather85}. The existence of a minimum-energy signed measure is not guaranteed in other circumstances, in particular when $K(\xb,\xb')=\Psi(\xb-\xb')$ and $\Psi$ is differentiable at 0; see Example~\ref{Ex:Antoniadis} below.

\subsubsection{Discrete BLUE}\label{S:discreteBLUE}
Consider the framework of Section~\ref{S:BQ}, with the same notation, and suppose that the $n$ design points $\xb_i$ in $\Xb_n$ are fixed.
Any linear estimator of $\beta_0$ in \eqref{model1} has then the form $\tilde \beta_0^n=\wb_n^T \yb_n$, with $\wb_n=(w_1,\ldots,w_n)^T\in\mathds{R}^n$. The unbiasedness constraint imposes $\wb_n^T\1b_n=1$. The variance of $\tilde \beta_0^n$ equals $\ms^2\wb_n^T \Kb_n \wb_n$, and the BLUE corresponds to the estimator $\betah_0^n$ given by \eqref{beta01} (we assume that $\Kb_n$ is nonsingular). The minimum-energy signed measure in $\SM(1)$ (here discrete) $ \mu_K^*$ is defined by the weights $\wb_n^*=\Kb_n^{-1}\1b_n/(\1b_n^T\Kb_n^{-1}\1b_n)$ set on the points in $\Xb_n$; its energy is $\SE_K( \mu_K^*)={\wb_n^*}^T \Kb_n \wb_n^* = 1/(\1b_n^T\Kb_n^{-1}\1b_n)$
and the variance of the BLUE equals $\ms^2 \SE_K( \mu_K^*)$. Note that some components of $\wb_n^*$ may be negative and that the potential associated with the measure $ \mu_K^*/\SE_K( \mu_K^*)$ on $\SX=\Xb_n$ gives the constant function $1_\SX=\1b_n$, see Theorem~\ref{P:capacity-1}. The optimal design problem for the discrete BLUE thus corresponds to the determination of the $n$-point set maximizing $\1b_n^T\Kb_n^{-1}\1b_n$.

\begin{example}\label{Ex:Antoniadis} Consider $K(x,x')=\exp(-\mt |x-x'|)$, $\mt>0$, for $x,x'\in\SX=[0,1]$. $K$ is ISPD and satisfies
$$
1_\SX = \frac{K(\cdot,0)+K(\cdot,1)}{2} + \frac{\mt}{2}\, \int_0^1 K(\cdot,x) \dd x \,,
$$
so that $1_\SX\in\SP_K$, see \cite{Antoniadis84}.
The minimum-energy measure in $\SM(1)$ is $\mu_K^*=(\delta_0+\delta_1+\mt\mu_L)/(\mt+2)$, with $\mu_L$ the Lebesgue measure on~$\SX$, and $\mu_K^*\in\SM^+(1)$.
The BLUE of $\beta_0$ in \eqref{model1} is $\hat \beta_0^*=\int_\SX f(x)\, \dd\mu_K^*(x)$, its variance equals $\ms^2\SE_K(\mu_K^*)=2\ms^2/(\mt+2)$, see \cite[p.~56]{Nather85}.
Note that $K'=K-2/(\mt+2)$ is still positive definite, but $1_\SX \not\in \SH_{K'}$ since $c^2 K'-1$ is not positive definite for any $c \neq 0$, see, e.g., \cite[p.~30]{BerlinetT-A2004}, \cite[p.~20]{Paulsen2009}.

Consider now the squared exponential kernel $K(x,x')=\exp(-\mt|x-x'|^2)$, $\mt>0$. The constant $1_\SX$ does not belong to $\SH_K$
\cite{SteinwartHS2006} and the BLUE of $\beta_0$ in \eqref{model1} is not defined for that kernel. On the other hand, the discrete BLUE \eqref{beta01} is well defined for any set of $n$ distinct points $x_i$,
$\hat \beta_0^n={\wb_n^*}^T\yb_n=\1b_n^T\Kb_n^{-1}\yb_n/(\1b_n^T\Kb_n^{-1}\1b_n)$. Suppose that the $n$ points $x_i$ are equally spaced in $\SX=[0,1]$. The process $Z_x$ in \eqref{model1} has mean square derivatives of all orders, and, roughly speaking, for large $n$ the construction of the BLUE mimics the estimation of derivatives of $f$ and the weights $w_i^*$ strongly oscillate between large positive and negative values. Figure~\ref{F:wstar}-Left shows the optimal weights $(w_i^*/|w_i^*|)(\log_{10}(\max\{|w_i^*|,1\})$, truncated to absolute values larger than 1 and in log scale, when $x_i=(i-1)/(n-1)$, $i=1,\ldots,n=101$. In Figure~\ref{F:wstar}-Right, the kernel is $K(x,x')=(1+\sqrt{5}|x-x'|+5|x-x'|^2/3)\exp(-\sqrt{5}|x-x'|)$ (Mat\'ern 5/2), so that $Z_x$ is twice mean-square differentiable; the construction of the BLUE mimics the estimation of the first and second order derivatives of $f$ at $0$ and $1$, suggesting that $1_\SX\not\in\SH_K$ in that case too; see \cite{DettePZ2016} for more details.
\fin

\begin{figure}[ht!]
\begin{center}
\includegraphics[width=.49\linewidth]{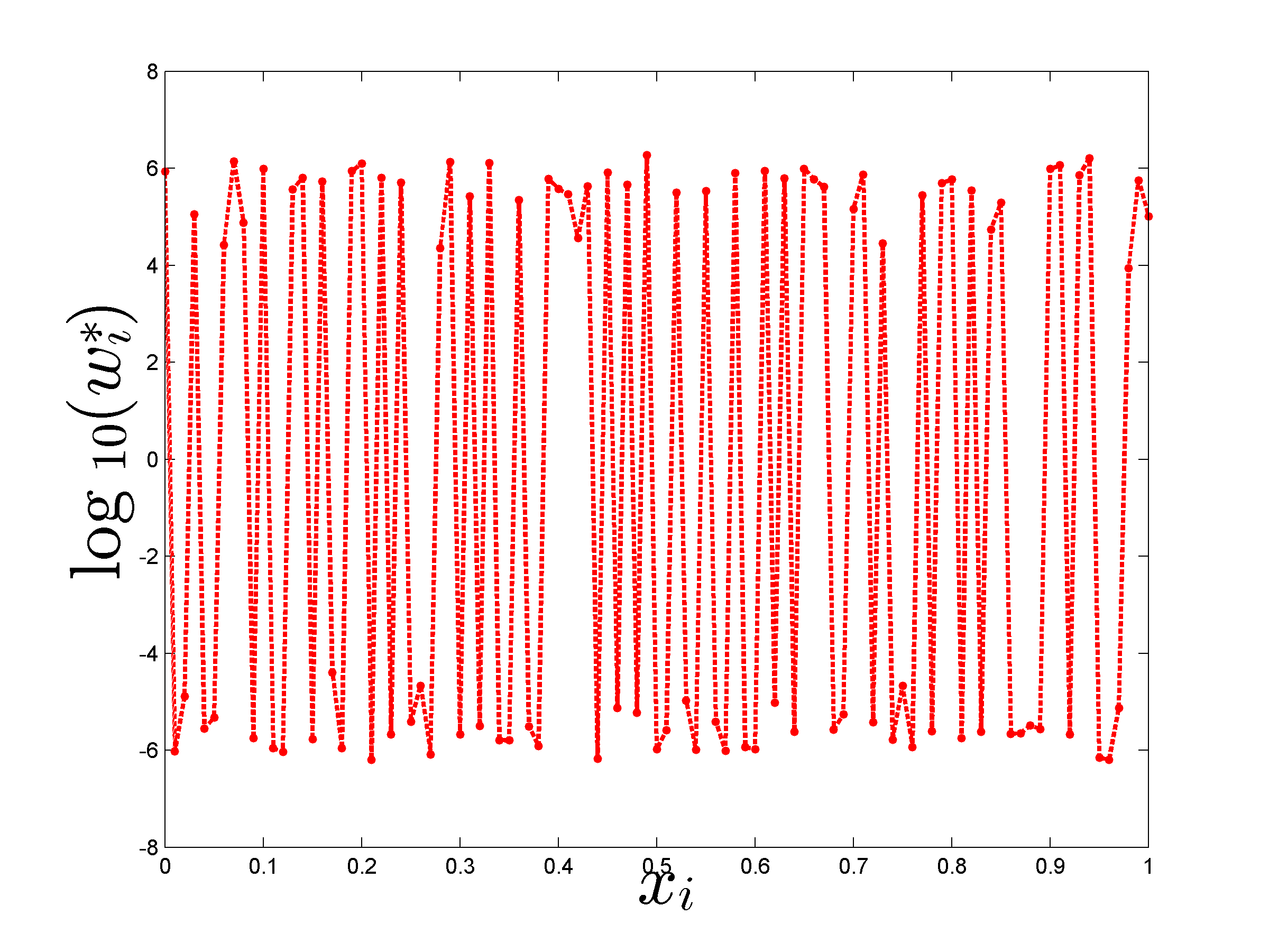} \includegraphics[width=.49\linewidth]{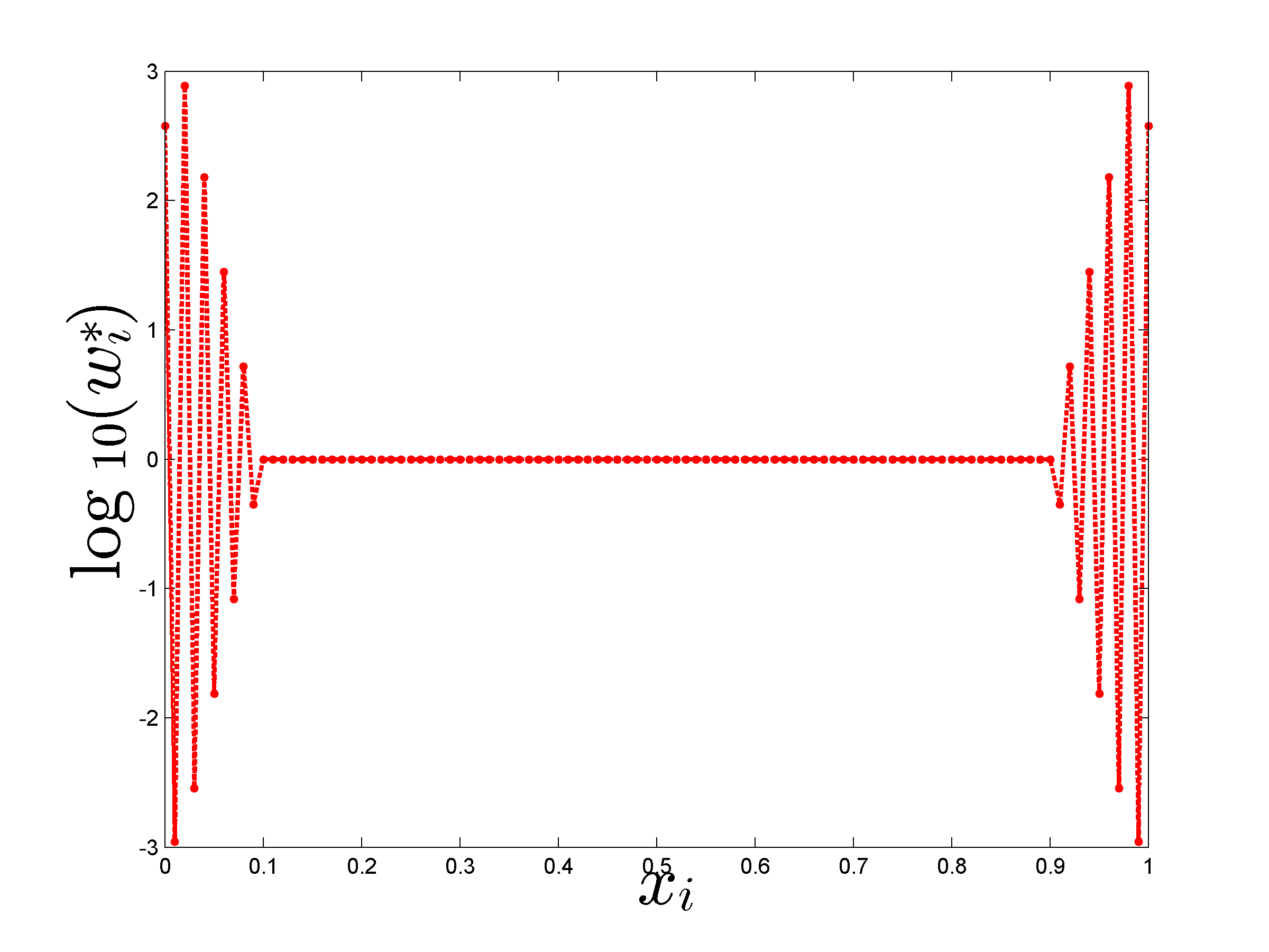}
\end{center}
\caption{\small BLUE weights $(w_i^*/|w_i^*|)(\log_{10}(\max\{|w_i^*|,1\})$ for $x_i=(i-1)/(n-1)$, $i=1,\ldots,n=101$. Left: $K(x,x')=\exp(-|x-x'|^2)$, Right: $K(x,x')=(1+\sqrt{5}|x-x'|+5|x-x'|^2/3)\exp(-\sqrt{5}|x-x'|)$ (Mat\'ern 5/2).}
\label{F:wstar}
\end{figure}

\end{example}

\vsp
Although a minimum-energy signed measure may not exist, in the next section we shall see how, for any measure $\mu\in\SM(1)$ and any CISPD kernel $K$, we can modify $K$ in such a way that the minimum-energy signed measure for the modified kernel exists (and coincides with $\mu$).

\subsection{Equilibrium measure and kernel reduction}\label{S:D-EM-RK}

Minimum-energy signed measures, when they exist, satisfy the following property.

\begin{lem}\label{P:ET-2}
If $K$ is CISPD and if a minimum-energy signed measure $\mu_K^*$ exists in $\SM(1)$, we have
$\SE_K(\xi) = \SE_K[\xi-\xi(\SX)\mu_K^*] + [\xi(\SX)]^2 \, \SE_K(\mu_K^*) \,, \ \forall\, \xi\in\SM_K$.
\end{lem}

\noindent {\em Proof.}
For any $\xi\in\SM_K$, direct calculation gives
\bea
\SE_K[\xi-\xi(\SX)\mu_K^*] &=& \SE_K(\xi)+ [\xi(\SX)]^2\, \SE_K(\mu_K^*)  - 2 \xi(\SX)\, \int_{\SX^2} K(\xb,\xb')\,\dd\mu_K^*(\xb)\dd\xi(\xb') \\
&=& \SE_K(\xi) - [\xi(\SX)]^2 \, \SE_K(\mu_K^*) \,,
\eea
where the second equality follows from \eqref{ET-M(1)}.
\carre

\vsp
Under the conditions of Lemma~\ref{P:ET-2}, any $\xi\in\SM(1)$ satisfies
$$
\SE_K(\xi)=\SE_K(\xi-\mu_K^*)+ \SE_K(\mu_K^*) \,,
$$
where the first term on the right-hand side equals the squared MMD $\mg_K^2(\xi,\mu_K^*)$, see \eqref{kernel discrepancyB}, and the second term does not depend on $\xi$. Minimizing the energy $\SE_K(\xi)$ is thus equivalent to minimizing the MMD $\mg_K(\xi,\mu_K^*)$. However, ($i$) $\mu_K^*$ may not exist, ($ii$) in many situations we wish to select a measure $\xi$ having small MMD $\mg_K(\xi,\mu)$ for \emph{a given} measure $\mu\in\SM_K$. This is the case in particular when one aims at evaluating the integral of a function with respect to some $\mu\in\SM^+(1)$ (Section~\ref{S:BQ}), or when we want construct a space-filling design in~$\SX$, $\mu$ being then uniform.

\subsubsection{Kernel reduction}

Take any $\mu\in\SM_K$ such that $\mu(\SX)\neq 0$. Without any loss of generality, we assume $\mu\in\SM(1)$. Following \cite{DamelinHRZ2010}, we show how to modify the kernel $K$ in such a way that minimizing the energy $\SE_{K_\mu}(\xi)$, $\xi\in\SM(1)$, for the new (reduced) kernel $K_\mu$ is equivalent to minimizing $\mg_{K_\mu}(\xi,\mu)$.

Define
\be\label{K_mu}
K_\mu(\xb,\xb')=K(\xb,\xb') - P_\mu(\xb) - P_\mu(\xb') + \SE_K(\mu)\,,
\ee
see \cite{Schaback199}.
One can readily check that the energy for this new reduced kernel $K_\mu$ satisfies $\SE_{K_\mu}(\beta\mu)=0$ for any real $\beta$ and that the potential for $\mu$ associated with $K_\mu$ satisfies $\widetilde P_\mu(\xb)=\int_\SX K_\mu(\xb,\xb')\, \dd\mu(\xb') = 0$ for all $\xb$. 

Next theorem indicates that, for any given $\mu$ in $\SM(1)\cap\SM_K$, when considering signed measures $\xi$ with total mass one, minimizing the energy $\SE_{K_\mu}(\xi)$ is equivalent to minimizing the MMD $\mg_K(\xi,\mu)$, provided that $K$ is CISPD.

\begin{theo}\label{P:Energy-RK}
If $K$ is CISPD, then for any $\mu\in\SM(1)\cap\SM_K$, we have\\
($i$) the reduced kernel $K_\mu$ defined by \eqref{K_mu} is CISPD; \\
($ii$) $\mu$ is the minimum-energy measure in $\SM(1)$ for $K_\mu$, and
\bea\label{Energy-RK}
\forall\, \xi\in\SM_K\,, \ \SE_{K_\mu}(\xi)=\SE_K[\xi-\xi(\SX)\mu] = \SE_{K_\mu}[\xi-\xi(\SX)\mu]  \,.
\eea
\end{theo}

\noindent{\em Proof.} For any nonzero $\xi\in\SM_K$, direct calculation using \eqref{K_mu} gives
\be
\SE_{K_\mu}(\xi) &=& \SE_K(\xi) - 2\xi(\SX)\, \int_{\SX^2} K(\xb,\xb')\,\dd\mu(\xb)\dd\xi(\xb')+ [\xi(\SX)]^2 \, \SE_K(\mu) \nonumber \\
&=& \SE_K[\xi-\xi(\SX)\mu] \,. \label{Energy-RKB}
\ee
($i$) When $\xi(\SX)=0$ we get $\SE_{K_\mu}(\xi)=\SE_K(\xi)$ which is strictly positive when $\xi\neq 0$, showing that $K_\mu$ is CISPD.
($ii$) Since $[\xi-\xi(\SX)\mu](\SX)=0$ and $K$ is CISPD, $\SE_{K_\mu}(\xi)\!>\!0$ for $\xi\neq \xi(\SX)\mu$, showing that $\mu$ is the (unique) minimum-energy signed measure in $\SM(1)$ for $K_\mu$.
%
Since $\SE_{K_\mu}(\mu)=0$, Lemma~\ref{P:ET-2} with $K_\mu$ substituted for $K$ implies that $\SE_{K_\mu}(\xi)=\SE_{K_\mu}[\xi-\xi(\SX)\mu]$ for any $\xi\in\SM_K$, which, together with \eqref{Energy-RKB}, concludes the proof.
\carre

\subsubsection{Kernel reduction, BLUE and Bayesian integration}\label{S:BLUE-Kreduction}
Consider again the situation of Section~\ref{S:BLUE}, and define $\mathcal{P}_1$ as the orthogonal projection of $L^2(\SX,\mu)$ onto the linear space spanned by the constant $1$; see \cite{GP-CSDA2016}. The model \eqref{model1} can then we written as
\be\label{model1p}
f(\xb)= \beta_0 + {\mathcal{P}_1}Z_x + (\Id_{L^2}-{\mathcal{P}_1})Z_x = \beta_0'+ \widetilde Z_x\,,
\ee
where $\beta_0'=\beta_0 + {\mathcal{P}_1}Z_x$ and $\widetilde Z_x=(\Id_{L^2}-{\mathcal{P}_1})Z_x$, with $\widetilde Z_x$ having zero mean and covariance $\Exx\{\widetilde Z_x \widetilde Z_{x'}\} = \ms^2 K_\mu(\xb,\xb')$.
We have seen in Section~\ref{S:BLUE} that the variance of the continuous BLUE of $\beta_0$ equals $\ms^2 \SE_K( \mu_K^*)$ provided that the minimum-energy signed measure $ \mu_K^*$ exists. (Note that the prior on $\beta_0'$ remains non-informative when the prior on $\beta_0$ is non-informative.) On the other hand, we obtain now that the continuous BLUE of $\beta_0'$ always exists: it coincides with $I_\mu(f)$ and its variance is $\ms^2 \SE_{K_\mu}(\mu)=0$. Therefore, as mentioned in introduction, Bayesian integration for the model \eqref{model1} with correlated errors is equivalent to parameter estimation in a location model with different correlation structure.

\subsection{Tensor product kernels}\label{S:TensorP}

From $d$ kernels $K_i$ respectively defined on $\SX_i\times\SX_i$, $i=1,\ldots,d$, we can construct a tensor product kernel as
\be\label{TP}
K^\otimes(\xb,\xb')=\prod_{i=1}^d K_i(x_i,x'_i)\,,
\ee
where $\xb=(x_1,\ldots,x_d)^T$ and $\xb'=(x'_1,\ldots,x'_d)^T$ belong to the product space $\SX=\SX_1 \times \cdots\times \SX_d$. The construction is particularly useful when considering product measures on~$\SX$, since, in some sense, it allows us to decompose an integration or space-filling design problem in a high dimensional space into its one-dimensional counterparts.
Suppose that each $K_i$ is uniformly bounded and CISPD on $\SM^{(i)}=\SM[\SX_i]$; that is, $K_i$ is ISPD on $\SM^{(i)}(0)$, see Definitions~\ref{D:ISPD} and \ref{D:CISPD}. One can show that this is equivalent to $K^\otimes$ being ISPD on $\otimes_{i=1}^d \SM^{(i)}(0)$, see \cite[Th.~2]{SzaboS2017}. In the same paper, the authors prove (Th.~4) that if each $K_i$ is moreover continuous and translation invariant, then $K^\otimes$ is ISPD on $\SM(0)$; that is, $K^\otimes$ is CISPD on $\SM$. Their proof relies on the equivalence between the CISPD and characteristic properties for uniformly bounded kernels, and on the characterization of characteristic continuous, uniformly bounded and translation invariant kernels through a property of the support of the measure $\Lambda$ defined in \eqref{Lambda}; see Section~\ref{S:MMD}.

An important property of tensor product kernels $K^\otimes$ is that kernel reductions $K^\otimes_\mu$, see \eqref{K_mu}, are easily obtained explicitly. Indeed, when $\mu=\otimes_{i=1}^d \mu^{(i)}$ is a product measure on~$\SX$, then, for all $\xb\in\SX$,
\be
\SE_{K^\otimes}(\mu) &=& \prod_{i=1}^d \SE_{K_i}(\mu^{(i)}) \,, \label{E-T} \\
P_\mu(\xb) &=& \prod_{i=1}^d \int_{\SX_i} K_i(x_i,x'_i)\, \dd\mu^{(i)}(x'_i) = \prod_{i=1}^d P_{\mu^{(i)}}(x_i) \,,  \label{Pot-T}
\ee
which facilitates the calculation of $\SE_{K^\otimes_\mu}(\xi)$, in particular when $\xi$ is a discrete measure as considered in Section~\ref{S:Empirical}. Table~\ref{Tb:EandP} gives the expressions of $\SE_K(\mu)$ and $P_{\mu}(x)$ obtained for a few kernels, with $\mu$ uniform on $\SX=[0,1]$; the expressions for the squared exponential and Mat\'ern kernels can be found in \cite{GinsbourgerRSDL2014}. Note that in each case $\SE_{K_\mu}(\xi)>0$ for any $\xi\in\SM(1)$, $\xi\neq \mu$.

\begin{table}[ht!]
\begin{center}
{\scriptsize
\caption{\small Energy $\SE_K(\mu)$ and potential $P_{\mu}(x)$ for different kernels $K$ with $\mu$ uniform on $\SX=[0,1]$; $P_{\mu}(x)=S_\mu(x)+S_\mu(1-x)+T_\mu(x)$; $S_\mu(\cdot)$ is continuously differentiable in $(0,1]$, $T_\mu=0$ when $K$ is translation invariant.}
\label{Tb:EandP}
\begin{tabular}{lll}
\toprule
$K(x,x')$ & $\SE_K(\mu)$ & $S_{\mu}(x)$ [and $T_\mu(x)$] \\ 
\midrule
$\exp(-\mt |x-x'|)$ & $2(\mt+\e1^{-\mt}-1)/\mt^2$ & $x(1-\e1^{-\mt|x|})/(\mt|x|)$ \\
$K_{3/2,\mt/\sqrt{3}}(x,x')$ in \eqref{K32} & $ 2[\mt(2+\e1^{-\mt})+3(\e1^{-\mt}-1)]/\mt^2$ & $x[2-(2+\mt|x|)\e1^{-\mt|x|}]/(\mt|x|)$ \\
%
$[(x-x')^2+\me]^{-1}$ ($\me\geq 0$) & $(2/\sqrt{\me})\,\arctan(1/\sqrt{\me})-\log(1+1/\me)$ & $(1/\sqrt{\me})\arctan(x/\sqrt{\me})$ \\
%
$(|x-x'|+\me)^{-1}$ ($\me> 0$)& $2\, [(1+\me)\,\log(1+1/\me)-1]$ & $\mathrm{sign}(x) \log(1+|x|/\me)$\\
%
$(|x-x'|+\me)^{-1/2}$ ($\me> 0$)& $4\me^{3/2}\, [2(1+1/\me)^{3/2}-2-3/\me]/3$ & $2\sqrt{\me}\,\mathrm{sign}(x)(\sqrt{1+|x|/\me}-1)$ \\
%
$1-\mt\, |x-x'|$ ($0<\mt\leq 1$) & $1-\mt/3$ & $1/2-\mt x|x|/2$ \\
%
$|x-x'|^{-s}$ ($0<s<1$) & $2/(s^2-3s+2)$ & $x/[(1-s)|x|^s]$\\
%
$-\log |x-x'|$ & 3/2 & $1/2-x\log|x|$\\
$|x|+|x'|-|x-x'|$ & $2/3$ & $1/4-x|x|/2$ \ [$T_\mu(x)=|x|$]\\
$\sqrt{|x|}+\sqrt{|x'|}-\sqrt{|x-x'|}$ & $4/5$ & $1/3-2x\sqrt{|x|}/3$ \ [$T_\mu(x)=\sqrt{|x|}$]\\
\bottomrule
\end{tabular}
}
\end{center}
\end{table}

\begin{remark}[super- and subharmonicity for tensor product kernels]\label{R:subharmonicTPK} In complement of Theorem~\ref{Th:MAIN}, we may notice that when each $K_i$ in \eqref{TP} is stationary and satisfies $K_i(x,x')=\Psi_i(x-x')$, then $K^\otimes(\xb,\xb')=\Psi(\xb-\xb')$ with $\Psi(\zb)=\prod_{i=1}^d \Psi_i(z_i)$. Assume that all $\Psi_i$ are twice continuously differentiable. Then, the Laplacian of $\Psi$ is
$$
\Delta_\Psi(\xb)=\sum_{i=1}^d \left[\frac{\mp^2 \Psi_i(x_i)}{\mp x_i^2} \prod_{j\neq i} \Psi_j(x_j)\right] \,. \quad\quad \fin
$$
\end{remark}

\section{Experimental design}\label{S:Empirical}

Consider an $n$-point design $\Xb_n=\{\xb_1,\ldots,\xb_n\}$, with $\xb_i\in\SX$ for all $i$. In this section, we shall restrict our attention to finite signed measures $\xi_n=\sum_{i=1}^n w_i\, \delta_{\xb_i}$ supported on $\Xb_n$, and denote $\wb_n=(w_1,\ldots,w_n)^T$. As in Section~\ref{S:BQ}, we consider a measure $\mu\in\SM^+(1)$, with special attention to space-filling design for which $\mu$ is uniform on a compact subset~$\SX$ of $\mathds{R}^d$. We assume that $K$ is SPD and $\mu$ has finite energy $\SE_K(\mu)$, see \eqref{overline K}.
Direct calculation gives
\be
\mg_K^2(\xi_n,\mu) = \SE_K(\xi_n-\mu) &=&  \wb_n^T \Kb_n \wb_n - 2 \wb_n^T \pb_n(\mu)+ \SE_K(\mu) \,, \nonumber \\
&=& \sum_{i,j} w_i w_j\, K(\xb_i,\xb_j) - 2\, \sum_{i=1}^n w_i\,P_\mu(\xb_i) + \SE_K(\mu) \,, \label{JK_BQ2}
\ee
where $\{\Kb_n\}_{i,j}=K(\xb_i,\xb_j)$, $i,j=1,\ldots,n$,
and $\pb_n(\mu)$ is given by \eqref{hb}. Note that $\SE_K(\mu)$ and the $P_\mu(\xb_i)$ have simple expressions when $K$ is a tensor product kernel and $\mu=\otimes_{i=1}^d \mu^{(i)}$ is a product measure on $\SX=\SX_1\times\cdots\times\SX_d$, see (\ref{E-T}, \ref{Pot-T}). Monte-Carlo approximation, based on a large i.i.d.\ sample from $\mu$, or a low-discrepancy sequence, can always be used instead.

\subsection{One-shot designs}\label{S:n-fixed}

\subsubsection{Support of empirical measures}\label{S:empirical}

Denote by $\xi_{n,e}=\xi_{n,e}(\Xb_n)$ the empirical measure associated with a given design $\Xb_n=\{\xb_1,\ldots,\xb_n\}\subset\mathds{R}^{nd}$, $\xi_{n,e}=(1/n) \sum_{i=1}^n \delta_{\xb_i}$. As indicated hereafter, the literature on space-filling design provides several examples of construction of $n$-point designs through the minimization of the squared MMD $\SE_K(\xi_{n,e}-\mu)$ with respect to $\Xb_n$.

For $\SX=[0,1]^d$, tensorised kernels based on variants of Brownian motion covariance yield $L_2$ discrepancies (symmetric, centred, wrap-around and so on); see, e.g., \cite{Hickernell1998}, \cite[Chap.~3]{FangLS2006}. For instance, for $\SX=[0,1]$ and $K(x,x')=1-|x-x'|$ (for which the expressions of $\SE_K(\mu)$ and $P_{\mu}(x)$ are given in Table~\ref{Tb:EandP}), $\SE_{K_\mu}(\xi_{n,e})$ is twice the squared $L_2$ star discrepancy for $d=1$.

The ISPD kernel $K^\otimes_{s,\me}(\xb,\xb') = \prod_{i=1}^d K_{s,\me}(x_i,x'_i)$, with $K_{s,\me}$ given by \eqref{Kes} with $s>0$ and $\me>0$, is called \emph{projection kernel} in \cite{MakJ2017b}. For very small $\me$, the minimization of  $\SE_{K^\otimes_{1,\me}}(\xi_{n,e})$ corresponds to the construction of a maximum-projection design, as defined in \cite{JosephGB2015}. Note that minimizing $\SE_{K^\otimes_{s,\me}}(\xi_{n,e})$ is not equivalent to minimizing $\SE_{K^\otimes_{s,\me}}(\xi_{n,e}-\mu)$: in particular, when $\mu$ is uniform on~$\SX$, which is assumed to be compact and convex, the former tends to push design points to the boundary of~$\SX$ whereas the latter keeps all points in the interior of~$\SX$; see \cite{MakJ2017b}.

In \cite{MakJ2017}, space-filling designs in a compact set $\SX\subset \mathds{R}^d$ are constructed by minimizing $\SE_{K^{(1)}}(\xi_{n,e}-\mu)$ for $\mu$ uniform on~$\SX$, see \eqref{Riesz<0}. They call \emph{support points} the optimal support $\Xb_n^*$, which they determine via a majorization-minimization algorithm using the property that the problem can be formulated as a difference-of-convex optimization problem. Values of  $\SE_{K^{(1)}}(\mu)$ and $P_{\mu}(\xb)$ are not available even for $\SX=[0,1]^d$ and Monte-Carlo approximation is used.

\subsubsection{Space-filling design through Bayesian quadrature}\label{S:BQ2}

Since $s_n^2$ given by \eqref{sn1} does not depend on the function $f$ considered, a design $\Xb_n$ for Bayesian integration can in principle be chosen beforehand, by direct minimization of $s_n^2$. This corresponds to the approach followed in \cite{O'Hagan91} where several quadrature rules are tabulated (for several values of $n$). Next theorem shows the connection between the minimum of $\SE_K(\xi_n-\mu)$ with respect to weights $\wb_n$ and the posterior variances $s_n^2$ and $s_{n,0}^2$, see \eqref{sn1} and \eqref{sn10}. We assume that all points in $\Xb_n$ are pairwise different and $\mu$ is not fully supported on $\Xb_n$.

\begin{theo}\label{P:BQ-KD} Let $K$ be an SPD kernel and let $\mu\in\SM^+(1)\cap\SM_K$.

($i$) The optimal unconstrained weights $\wb_n^*$ that minimize $\SE_K(\xi_n-\mu)$ are $\wb_n^*= \Kb_n^{-1}\pb_n(\mu)$ and the corresponding measure $\xi_n^*$, with weights $\wb_n^*$, satisfies
\be\label{JKopt}
\SE_K(\xi_n^*-\mu) = s_{n,0}^2\,,
\ee
with $s_{n,0}^2$ given by \eqref{sn10}.

($ii$) The optimal weights $\hat\wb_n$ that minimize $\SE_K(\xi_n-\mu)$ under the constraint $\wb_n^T\1b_n=\sum_{i=1}^n w_i=1$ are
\be\label{wopt-sum=1}
\hat \wb_n = \left(\Kb_n^{-1} - \frac{\Kb_n^{-1} \1b_n\1b_n^T \Kb_n^{-1}}{\1b_n^T\Kb_n^{-1}\1b_n} \right) \pb_n(\mu) + \frac{\Kb_n^{-1} \1b_n}{\1b_n^T\Kb_n^{-1}\1b_n} \,,
\ee
and the corresponding measure $\hat\xi_n$, with weights $\hat\wb_n$, satisfies
\be\label{JKopt-sum=1}
\SE_K(\hat\xi_n-\mu) = s_{n}^2\,,
\ee
with $s_{n}^2$ given by \eqref{sn1}; the estimator \eqref{In1} of the integral $I_\mu(f)$ is $\widehat I_n = \hat \wb_n^T \yb_n$.

($iii$) For any bounded signed measure $\xi_n=\sum_{i=1}^n w_i\, \delta_{\xb_i}$ we can write
\be\label{JK-sum=any}
\SE_K(\xi_n-\mu)=(\wb_n-\wb_n^*)^T\Kb_n(\wb_n-\wb_n^*)+\SE_K(\xi_n^*-\mu) \,,
\ee
and when the weights $w_i$ sum to one, we have
\be\label{JK-sum=1}
\SE_K(\xi_n-\mu)=(\wb_n-\hat\wb_n)^T\Kb_n(\wb_n-\hat\wb_n)+\SE_K(\hat\xi_n-\mu).
\ee
\end{theo}

\noindent{\em Proof.}
The expression for $\wb_n^*$, \eqref{JKopt} and \eqref{JK-sum=any} directly follow from the fact that $\SE_K(\xi_n-\mu)$ is quadratic in $\wb_n$, see \eqref{JK_BQ2}.
Since $K$ is SPD, straightforward calculation using Lagrangian theory indicates that the minimization of $\SE_K(\xi_n-\mu)$ under the constraint $\wb_n^T\1b_n=1$ gives \eqref{wopt-sum=1}
and \eqref{JKopt-sum=1}. Suppose that $\wb_n^T\1b_n=1$, then $\SE_K(\xi_n-\mu)=(\wb_n-\hat\wb_n+\hat\wb_n-\wb_n^*)^T\Kb_n(\wb_n-\hat\wb_n+\hat\wb_n-\wb_n^*)+\SE_K(\xi_n^*-\mu)$ gives \eqref{JK-sum=1} since $\Kb_n(\hat\wb_n-\wb_n^*)$ is proportional to $\1b_n$ and $(\wb_n-\hat\wb_n)^T\1b_n=0$.
\carre

\vsp
In the discrete case considered here, the minimum-energy signed measure $\hat\xi_n$ with total mass one always exists, but note that it is not necessarily a probability measure; that is, some weights $\hat w_i$ may be negative. Theorem~\ref{P:BQ-KD} can be extended to the case where $K$ is only conditionally SPD, but the computation of optimal weights $\hat\wb_n$ is more involved when $\Kb_n$ is singular; see Remark~\ref{R:optimal weights}.

Denote by $\widetilde\Kb_n$ the $n\times n$ matrix with elements $\{\widetilde\Kb_n\}_{i,j}=K_\mu(\xb_i,\xb_j)$, where $K_\mu$ is the reduced kernel \eqref{K_mu}; the corresponding vector of potential values at the $\xb_i$ is then $\widetilde\pb_n=(\widetilde P_\mu(\xb_1),\ldots,\widetilde P_\mu(\xb_n))^T=\0b$.
For measures $\xi_n$ in $\SM(1)$, in complement of ($ii$) of Theorem~\ref{P:BQ-KD}, we also have the following property.

\begin{theo}\label{P:BQ-RK}
For $K$ an SPD kernel, $\mu\in\SM^+(1)\cap\SM_K$ and $\xi_n\in\SM(1)$, we have
\be\label{JK-RK}
\SE_K(\xi_n-\mu) = \SE_{K_\mu}(\xi_n) = \wb_n^T \widetilde\Kb_n \wb_n \,.
\ee
The posterior mean \eqref{In1} and variance \eqref{sn1} of $I_\mu(f)$ are respectively given by
\be
\widehat I_n &=& \frac{\1b_n^T\widetilde\Kb_n^{-1}\yb_n}{\1b_n^T\widetilde\Kb_n^{-1}\1b_n}\,, \label{In2} \\
\ms^2 s_n^2 &=& \ms^2 (\1b_n^T \widetilde\Kb_n^{-1} \1b_n)^{-1} \,. \label{sn2}
\ee
\end{theo}

\noindent{\em Proof.}
Equation \eqref{JK-RK} follows from Proposition~\ref{P:Energy-RK}. Since we assumed that $\mu$ is not fully supported on $\Xb_n$ and $K$ is SPD, \eqref{JK-RK} gives $\inf_{\|\wb_n\|=1} \wb_n^T \widetilde\Kb_n \wb_n >0$ , which implies that $\widetilde\Kb_n$ has full rank. Direct calculation using \eqref{K_mu} gives $\widetilde\Kb_n=\Kb_n - \pb_n(\mu)\1b_n^T - \1b_n\pb_n^T(\mu) + \SE_K(\mu)\, \1b_n\1b_n^T$. The expression for $\widetilde\Kb_n^{-1}$ then yields $\1b_n^T \widetilde\Kb_n^{-1} \1b_n=1/s_n^2$, with $s_n^2$ given by \eqref{sn1}, proving \eqref{sn2}. The expansion of $(\1b_n^T\widetilde\Kb_n^{-1}\yb_n)/(\1b_n^T\widetilde\Kb_n^{-1}\1b_n)$ gives \eqref{In1}, which proves \eqref{In2}.
\carre

\begin{remark}[Optimal weights for CISPD kernels]\label{R:optimal weights}
Lagrangian theory indicates that the solution $\hat\wb_n$ is obtained by solving the linear equation $\Rb_n(\hat\wb_n^T \ \lambda)^T=(\0b^T \ 1)^T$, where
$$
\Rb_n = \left(
  \begin{array}{cc}
    \widetilde\Kb_n & \1b_n \\
    \1b_n^T & 0 \\
  \end{array}
\right) \,.
$$
When $K$ is conditionally SPD, $K_\mu$ is conditionally SPD too, and the matrix $\Rb_n$ has full rank $n+1$.
Indeed, $\Rb_n (\zb_n^T \ z)^T=0$ implies $\1b_n^T\zb_n=0$ and $\widetilde\Kb_n\zb_n+z\1b_n=0$. Multiplying the second equation by $\zb_n^T$, we get $\zb_n^T\widetilde\Kb_n\zb_n=0$. Since $K_\mu$ is conditionally SPD, this is incompatible with $\1b_n^T\zb_n=0$ unless $\zb_n=\0b$ and $z=0$. We obtain
$$
\hat\wb_n = \frac{(\widetilde\Kb_n+\1b_n\1b_n^T)^{-1}\1b_n}{\1b_n^T(\widetilde\Kb_n+\1b_n\1b_n^T)^{-1}\1b_n}\,,
$$
and $s_n^2=\hat\wb_n^T\widetilde\Kb_n\hat\wb_n= (\1b_n^T(\widetilde\Kb_n+\1b_n\1b_n^T)^{-1}\1b_n)^{-1}-1$. When $K$ is SPD and $\widetilde\Kb_n$ has full rank (Theorem~\ref{P:BQ-RK}), we recover $\hat\wb_n = \widetilde\Kb_n^{-1}\1b_n/(\1b_n^T\widetilde\Kb_n^{-1}\1b_n)$ and $\widehat I_n=\hat\wb_n^T \yb_n$ given by \eqref{In2}.
\fin
\end{remark}

\begin{remark}[BLUE and kernel reduction]\label{R:reducedK}
Equations \eqref{In2} and \eqref{sn2} indicate that $\widehat I_n$ is the BLUE of $\beta_0'$ and $\ms^2 s_n^2$ is its variance in the model
\eqref{model1p}, $f(\xb)= \beta_0'+ \widetilde Z_x$, see Sections~\ref{S:discreteBLUE} and \ref{S:BLUE-Kreduction}. A possible interpretation is as follows.
Predictions are not modified when using the reduced kernel $K_\mu$ instead of $K$,
that is, when considering model $f(\xb)= \beta_0'+ \widetilde Z_x$ instead of \eqref{model1},
see \cite[Sect.~5.4]{GP-CSDA2016}. It implies that the expressions \eqref{In1} and \eqref{sn1} of $\widehat I_n$ and $s_n^2$ are unchanged when replacing $K$ by $K_\mu$. Since, by construction, $\widetilde \pb_n(\mu)=\0b$ and $\SE_{K_\mu}(\mu)=0$ ($(\Id_{L^2}-{\mathcal{P}_1})Z_x$ has no contribution to the integral of $f$), we directly obtain \eqref{In2} and \eqref{sn2}.
\fin
\end{remark}

\begin{remark}[IMSPE for tensor product kernels]
The use of a tensor product kernel \eqref{TP} and a product measure $\mu=\otimes_{i=1}^d \mu^{(i)}$ on $\SX=\SX_1\times\cdots\times\SX_d$ facilitates the calculations of $\widetilde\Kb_n$ and $\SE_K(\xi_n-\mu)$, see \eqref{JK_BQ2}, since $\SE_K(\mu)$ and $P_\mu(\xb_i)$ have the simple expressions (\ref{E-T}, \ref{Pot-T}). The calculation of the IMSPE is facilitated too, but to a lesser extend. Indeed, we have
\bea
\int_\SX \rho_n^2(\xb)\, \dd\mu(\xb) = \SE_K(\mu) + \frac{1}{\1b_n^T\Kb_n^{-1}\1b_n}-2\,\frac{\pb_n^T(\mu)\Kb_n^{-1}\1b_n}{\1b_n^T\Kb_n^{-1}\1b_n} - \tr\left[\Kb_n^{-1}\Qb_n^\bot\Hb_n(\mu)\right],
\eea
see \eqref{MSE1}, where $\Qb_n^\bot = \Ib_n-\1b_n\1b_n^T\Kb_n^{-1}/(\1b_n^T\Kb_n^{-1}\1b_n)$, with $\Ib_n$ the $n$-dimensional identity matrix, is a projector onto the linear space orthogonal to $\1b_n$, and where $\Hb_n(\mu)$ is the symmetric non-negative definite $n\times n$ matrix with elements
\bea
\{\Hb_n(\mu)\}_{j,k}=\int_\SX K(\xb,\xb_j)K(\xb,\xb_k)\,\dd\mu(\xb) = \prod_{i=1}^d \int_{\SX_i} K_i(x_i,{x_j}_i)K_i(x_i,{x_k}_i)\,\dd\mu^{(i)}(x_i) \,. \fin
\eea
\end{remark}

\vsp
Theorems~\ref{P:BQ-KD} and \ref{P:BQ-RK} indicate that, if $K$ is SPD, $(\1b_n^T \widetilde\Kb_n^{-1} \1b_n)^{-1}$ is the minimum value of $\SE_K(\xi_n-\mu)$ for measures $\xi_n\in\SM(1)$. Hence, we can construct space-filling designs on a compact and convex subset~$\SX$ of $\mathds{R}^d$ by maximizing $\1b_n^T \widetilde\Kb_n^{-1} \1b_n$ with respect to $\Xb_n\in\mathds{R}^{nd}$, taking $\mu$ uniform on~$\SX$. This can be performed using any unconstrained nonlinear programming algorithm, as Example~\ref{Ex:4.1} will illustrate. Note that, from \eqref{JK-RK} and Cauchy-Schwarz inequality, $(\1b_n^T \widetilde\Kb_n^{-1} \1b_n)^{-1}\leq \SE_K(\xi_{n,e}-\mu) = (\1b_n^T\widetilde\Kb_n\1b_n)/n^2$, the minimization of which was considered in Section~\ref{S:empirical}.

\subsection{Any-time designs}\label{S:any-time}

There exist situations where the number $n$ of design points ultimately used (for integration, or function approximation) differs from that initially planned, say $N$. It is the case in particular when function evaluations are computationally more expensive than expected, and numerical experimentation is stopped after $n<N$ simulations, or when simulations fail at some design points and testing at more than $N$ points is required to obtain $N$ valid evaluations in total. In such circumstances, it is convenient to have sequences of nested designs at one's disposal. The objective is then to construct any-time designs; that is, ordered sequences $\xb_1,\xb_2,\ldots$ of designs points such that any design $\Xb_n=\{\xb_1,\ldots,\xb_n\}$ made of the first $n$ points of the sequence has good space-filling properties. A typical example is given by Low Discrepancy Sequences (LDS) in $[0,1]^d$, see \cite{Niederreiter92}.

\vsp
When $K$ is SPD, we may exploit expression \eqref{sn2} of the conditional variance of $I_\mu(f)$ in a greedy sequential construction: at step $n$ we choose $\xb_{n+1}$ that minimizes $s_{n+1}^2$. This sequential construction, called Sequential Bayesian Quadrature in \cite{BriolOGO2015}, is straightforward to implement compared with global minimization of $s_n^2$, see \eqref{sn1}. Direct calculation, using formulae for the inversion of the block matrix
$$
\widetilde\Kb_{n+1}=\left(
                      \begin{array}{cc}
                        \widetilde\Kb_n & \widetilde\kb_n(\xb_{n+1}) \\
                        \widetilde\kb_n^T(\xb_{n+1}) &  K_\mu(\xb_{n+1},\xb_{n+1})\\
                      \end{array}
                    \right)\,,
$$
where $\{\widetilde\Kb_n\}_{i,j}=K_\mu(\xb_i,\xb_j)$ and $\{\widetilde\kb_n(\xb)\}_i=K_\mu(\xb,\xb_i)$, $i,j=1,\ldots,n$, $\xb\in\SX$, gives
\be\label{sn+1^2}
s_{n+1}^2 = \left[\1b_n^T\widetilde\Kb_n^{-1}\1b_n+ \frac{(1-\widetilde\kb_n(\xb_{n+1})^T\widetilde\Kb_n^{-1}\1b_n)^2} {K_\mu(\xb_{n+1},\xb_{n+1})-\widetilde\kb_n^T(\xb_{n+1})\widetilde\Kb_n^{-1}\widetilde\kb_n(\xb_{n+1})} \right]^{-1} \,.
\ee
The sequential construction is thus
$$
\xb_{n+1}\in\Arg\max_{\xb\in\SX} \frac{(1-\widetilde\kb_n(\xb)^T\widetilde\Kb_n^{-1}\1b_n)^2} {K_\mu(\xb,\xb)-\widetilde\kb_n^T(\xb)\widetilde\Kb_n^{-1}\widetilde\kb_n(\xb)} \,.
$$

The conditional gradient algorithm of \cite{FrankW56} yields a simpler construction, particularly well adapted to the situation and also applicable when $K$ is unbounded. It relies on the sequential selection of points that minimize the current directional derivative of $\SE_K(\xi-\mu)=\SE_{K_\mu}(\xi)$, with $\xi$ supported on design points previously selected. The algorithm is initialized at a measure $\xi^{(n_0)}$ supported on $\Xb_{n_0}\in\SX^{n_0}$ (with for instance $n_0=1$ and $\xi^{(1)}=\delta_{\xb_1}$ for some $\xb_1\in\SX$). Let $\xi^{(n)}$ denote the measure associated with the current design $\Xb_n$ of iteration $n$, with weights $w_i^{(n)}$, i.e., $\xi^{(n)}=\sum_{i=1}^n w_i^{(n)} \delta_{\xb_i}$. Next design point is chosen in $\Arg\min_{\xb\in\SX} F_{K_\mu}(\xi^{(n)},\delta_\xb)$ (any minimizer can be selected in case there are several). Straightforward calculation using \eqref{Dir-der} gives $\xb_{n+1} \in\Arg\min_{\sb\in\SX} [ P_{\xi^{(n)}}(\xb) - P_\mu(\xb)]$, that is,
\be\label{x-n+1}
\xb_{n+1} \in\Arg\min_{\xb\in\SX} \left[ \sum_{i=1}^n w_i^{(n)} K(\xb,\xb_i) - P_\mu(\xb) \right] \,.
\ee
Note that this construction is well defined even if $K$ is singular: in that case, it ensures that all design points are different ($\xb_i\neq \xb_j$ for all $i,j$); the same is true for all one-dimensional canonical projections when $K$ is the tensorised product of singular kernels.

After choosing $\xb_{n+1}$, the measure $\xi^{(n)}$ is updated into
\be\label{xi-n+1}
\xi^{(n+1)} = (1-\ma_n) \xi^{(n)} + \ma_n \delta_{\xb_{n+1}}
\ee
for some $\ma_n\in[0,1]$, so that $\xi^{(n+1)}\in\SM^+(1)$ when $\xi^{(n)}\in\SM^+(1)$.
When $\xi^{(n_0)}$ is the empirical (uniform) measure on $\Xb_{n_0}$, the choice $\ma_n=1/(n+1)$ implies that $\xi^{(n)}$ remains uniform on its support $\Xb_n$ for all $n$, see \cite{Wynn70} for an early contribution in the design context. The method is called \emph{kernel herding} in the machine-learning literature, see \cite{BachLJO2012, ChenWS2012, HuszarD2012}. It is shown in \cite{ChenWS2012} that $\SE_K(\xi^{(n)}-\mu)=\SO(1/n^2)$ when $\SH_K$ is finite dimensional, but we only have the weaker result $\SE_K(\xi^{(n)}-\mu)=\SO(1/n)$ when $\SH_K$ is infinite dimensional, see \cite{BachLJO2012}.

\begin{remark}
Denote $\xi^{(n+)}(\xb)=[n/(n+1)]\xi^{(n)}+[1/(n+1)]\delta_\xb$. The direct minimization of $\SE_K(\xi^{(n+)}(\xb)-\mu)$ with respect to $\xb$ yields
$$
\xb^{(n+1)} \in \Arg\min_{\xb\in\SX} \left[ \frac{1}{n+1}\, \sum_{i=1}^k K(\xb,\xb_i) - P_\mu(\xb) + \frac{1}{2(k+1)}\, K(\xb,\xb) \right] \,,
$$
that is, a selection very close to \eqref{x-n+1} when $K(\xb,\xb)$ is constant (Mat\'ern kernel for instance). Note that this construction requires $K(\xb,\xb)<\infty$ for all $\xb\in\SX$, contrary to \eqref{x-n+1}.
\fin
\end{remark}

In practice $n$ is always smaller than some given $n_{\max}$, and to facilitate the construction we can restrict the choice of the $\xb_i$ to a finite subset $\SX_\mO=\{\sb_1,\ldots,\sb_\mO\}$ of~$\SX$, with $\mO\gg n_{\max}$ (when $\SX=[0,1]^d$, $\SX_\mO$ can be given by the first $\mO$ points of a LDS).
For any $n\leq n_{\max}$, we can write $\Xb_n=\{\xb_1,\ldots,\xb_n\}=\{\sb_{i_1},\ldots,\sb_{i_n}\}$,
the construction being initialized at some $n_0$-point design $\Xb_{n_0}\subset\SX_\mO$.
A measure $\xi$ supported on $\Xb_n$ can thus be written as
$\xi=\sum_{i=1}^\mO \mo_i \delta_{\sb_i}$, with $\mo_{i}=0$ when $\sb_i\not\in\Xb_n$. Therefore, for all $n$, $\xi^{(n)}$ is fully characterized by a $\mO$-dimensional vector $\mob^{(n)}=(\mo_1^{(n)},\ldots,\mo_\mO^{(n)})^T$, with $\mob^{(n)}$ in the probability simplex $\mathds{P}_\mO$ when $\xi^{(n)}\in\SM^+(1)$. The updating equations (\ref{x-n+1}, \ref{xi-n+1}) then imply that $\mob^{(n+1)}$ is obtained by moving $\mob^{(n)}$ in the direction of a vertex of $\mathds{P}_\mO$, hence the name \emph{vertex-direction} given to methods based on \eqref{xi-n+1} in the literature on optimal design, see, e.g., \cite[Chap.~9]{PP2013} and the references therein.
A summary of results on the rate of decrease of $\SE_K(\xi^{(n)}-\mu)$ in this situation is given in Appendix~A. The cost of the determination of $\xb_{n+1}$ in \eqref{x-n+1} is $\SO(\mO)$ (we need to compute $K(\xb,\xb_n)$ for all $\xb\in\SX_\mO$), and the cost for $n$ iterations scales as $\SO(n\mO)$ (including the initial cost for the computation of $P_\mu(\xb)$ for all $\xb\in\SX_\mO$). An $n$-point design constructed in this way can be used as initialization for the (unconstrained) minimization of $\SE_K(\xi_{n,e}-\mu)=\1b_n^T\widetilde\Kb_n\1b_n/n^2$ (Section~\ref{S:empirical}), or the maximization of $\1b_n^T \widetilde\Kb_n^{-1} \1b_n$ (Section~\ref{S:BQ2}), with respect to $\Xb_n$. The resulting design $\Xb_n^*$ can in turn be used as candidate set for the greedy construction of \cite{Gonzalez85} (also called coffee-house design in  \cite{Muller2000}), yielding a sequence of nested designs $\Xb_1,\ldots,\Xb_n=\Xb_n^*$.

\subsection{Illustrative examples}\label{Ex:4.1}
We take $\SX=[0,1]^2$, $n_{\max}=100$, $n_0=1$ with $\Xb_1=\{(0.5,0.5)\}$; $\mu$ is uniform on~$\SX$ and $\SX_\mO$ is given by the first $2^{12}$ points of Sobol' LDS. The kernel $K$ is the tensor product of uni-dimensional Mat\'ern 3/2 covariance functions $K_{3/2,\mt}$, see \eqref{K32}.

\begin{figure}[ht!]
\begin{center}
\includegraphics[width=.49\linewidth]{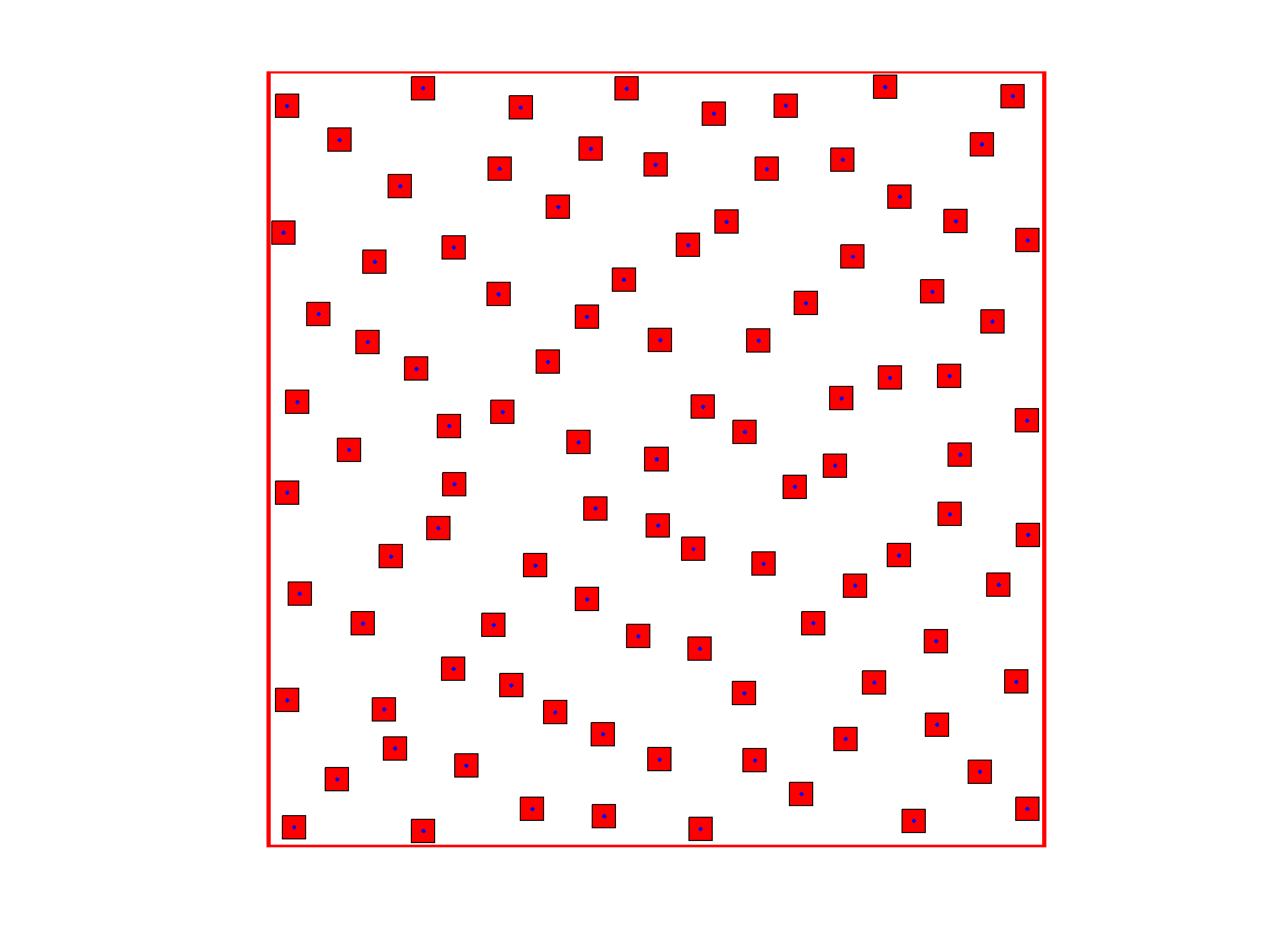} \includegraphics[width=.49\linewidth]{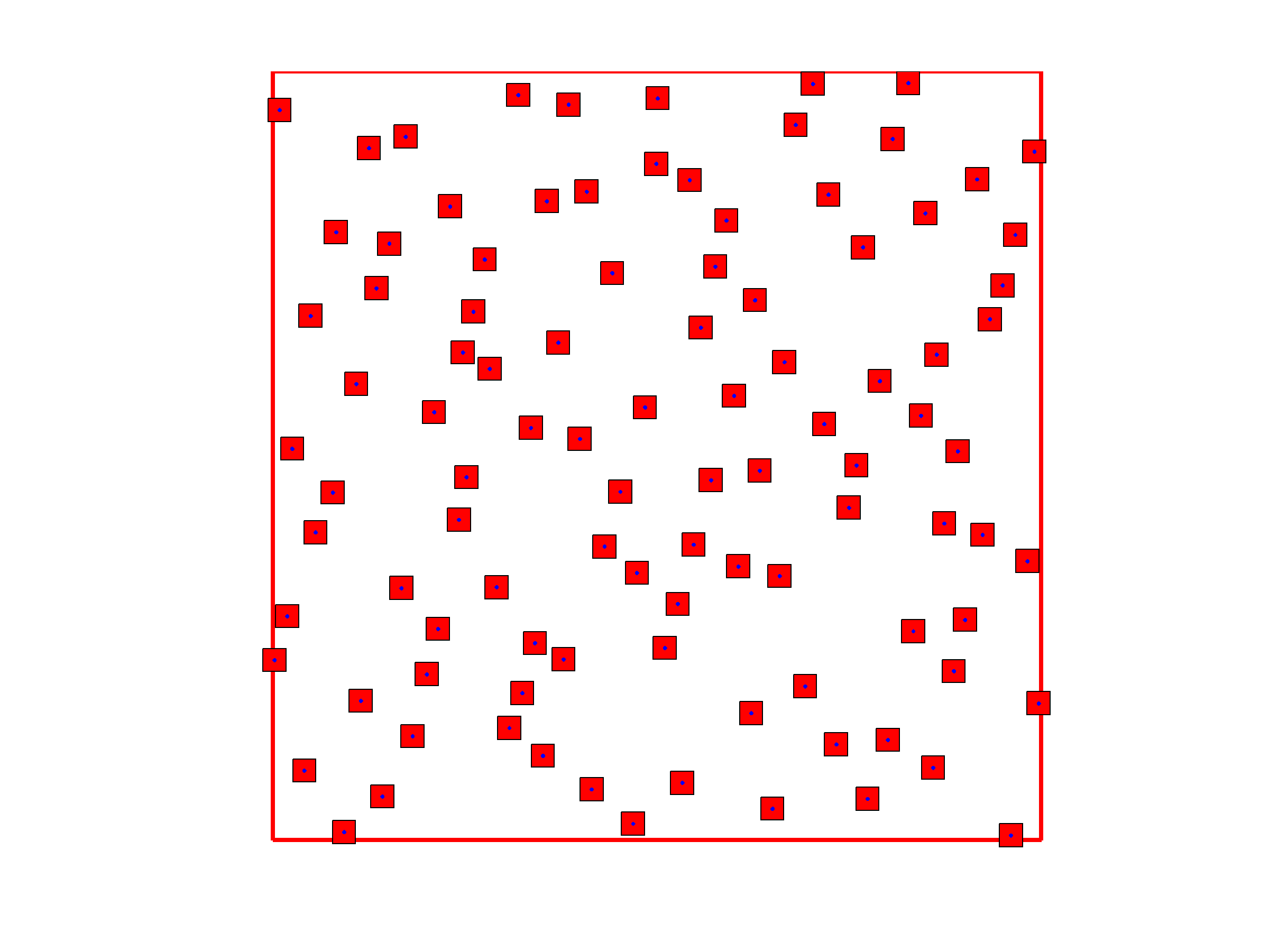}
\end{center}
\caption{\small Left: $\Xb_{100}=\Supp(\xi^{(100)})$ generated by (\ref{x-n+1}, \ref{xi-n+1}) with $\ma_n=1/(n+1)$ and $\mt=10$ in $K_{3/2,\mt}$. Right: first 100 points $\Sb_{100}$ of Sobol' LDS.}
\label{F:designs-theta10-d2-n100}
\end{figure}

Figure~\ref{F:designs-theta10-d2-n100}-Left shows the design $\Xb_{100}$ obtained after 99 iterations of (\ref{x-n+1}, \ref{xi-n+1}) with $\ma_n=1/(n+1)$ ($\Xb_{100}$ is the support of $\xi^{(100)}$), for $\mt=10$ in $K_{3/2,\mt}$. This design has visually better space-filling properties than the first 100 points $\Sb_{100}$ of Sobol' sequence presented on the right part of the figure. This is confirmed  by the numerical values of the covering and packing radii, respectively given by \eqref{CR} and \eqref{PR}: $\CR(\Xb_{100}) \simeq 0.0925 < \CR(\Sb_{100}) \simeq 0.1377$, and $\PR(\Xb_{100}) \simeq 0.0262 > \PR(\Sb_{100}) \simeq 0.0204$. Local maximization of $\1b_n^T \widetilde\Kb_n^{-1} \1b_n$ (Section~\ref{S:BQ2}) with respect to $\Xb_n$, initialized at $\Xb_{100}$, yields a design $\Xb_{100}^*$ with better space-filling properties: $\CR(\Xb_{100}^*) \simeq 0.0869$ and $\PR(\Xb_{100}^*) \simeq 0.0364$. When minimizing $\SE_K(\xi_{n,e}-\mu)$ with respect to $\Xb_n$ (Section~\ref{S:empirical}) we obtain $\CR(\Xb_{100}^*) \simeq 0.0955$ and $\PR(\Xb_{100}^*) \simeq 0.0359$. The performance is significantly worse, both in terms of $\CR$ and $\PR$, when the optimization is initialized at $\Sb_{100}$.

\begin{figure}[ht!]
\begin{center}
\includegraphics[width=.49\linewidth]{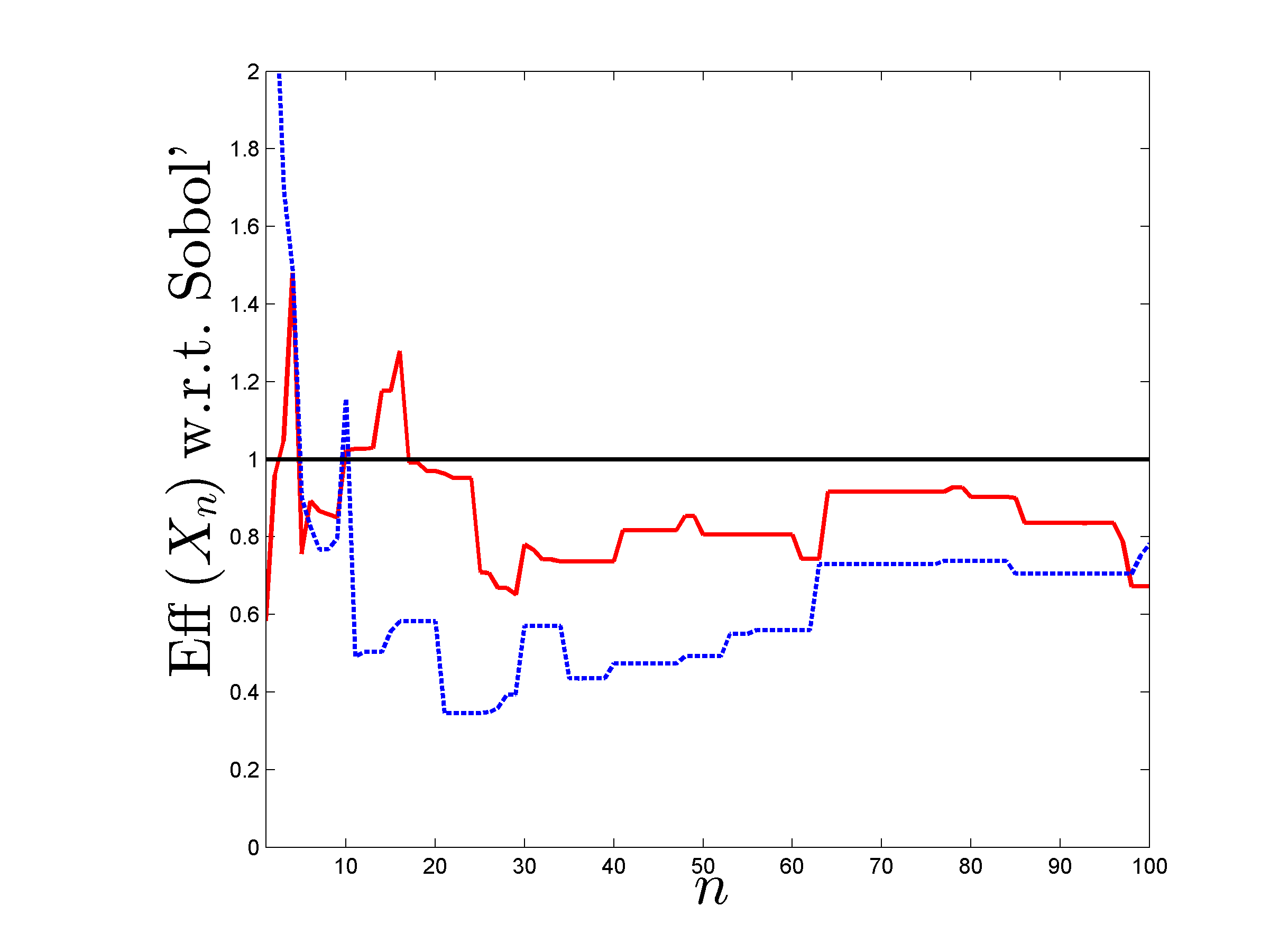} \includegraphics[width=.49\linewidth]{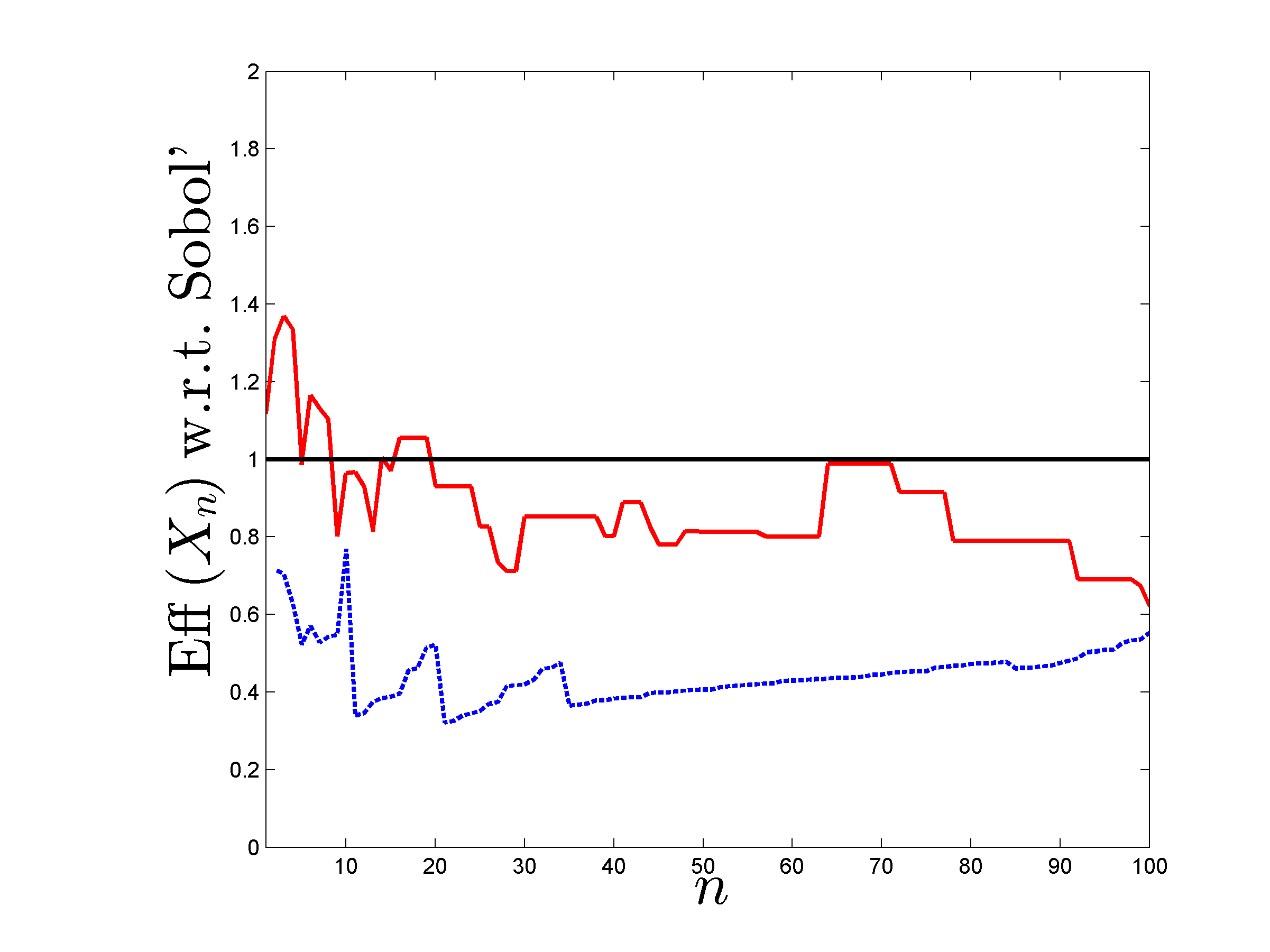}
\end{center}
\caption{\small Efficiencies $\CR(\Xb_{n})/\CR(\Sb_{n})$ (solid line) and $\PR(\Sb_{n})/\PR(\Xb_{n})$ (dashed line), with $\Sb_n$ corresponding to the first $n$ points of Sobol' sequence. Left: $\Xb_n=\Supp(\xi^{(n)})$ is generated by (\ref{x-n+1}, \ref{xi-n+1}) with $\ma_n=1/(n+1)$ and $\mt=10$ in $K_{3/2,\mt}$. Right: $\Xb_n$ is generated by coffee-house design with candidate set $\Xb_{100}^*$ obtained by local maximization of $\1b_n^T \widetilde\Kb_n^{-1} \1b_n$.}
\label{F:filling-theta10-d2-n100}
\end{figure}

Figure~\ref{F:filling-theta10-d2-n100}-Left presents the efficiencies, in terms of covering radius \eqref{CR} (solid line) and packing radius \eqref{PR} (dashed line) of $\Sb_n$ corresponding to the first $n$ points of Sobol' sequence relatively to $\Xb_n=\Supp(\xi^{(n)})$, when $\xi^{(n)}$ is generated by (\ref{x-n+1}, \ref{xi-n+1}) with $\ma_n=1/(n+1)$ and $\mt=10$ in $K_{3/2,\mt}$. Values smaller than one (shown by an horizontal line) indicate that $\Xb_n$ has better space-filling properties than $\Sb_n$.
A greedy coffee-house construction \cite{Gonzalez85, Muller2000} ensures $\CR(\Xb_n) \leq 2 \min_{\Xb'_n} \CR(\Xb'_n)$ and $\PR(\Xb_n) \geq (1/2)\, \max_{\Xb'_n} \PR(\Xb'_n)$, but the method is computationally expensive and tends to choose points on the border of~$\SX$. On the other hand, the construction is straightforward (and efficient) when we use the design $\Xb_{n_{\max}}^*$ obtained by maximizing $\1b_n^T \widetilde\Kb_n^{-1}\1b_n$ as candidate set. The first point $\xb^{(1)} $ is chosen as the point in $\Xb_{n_{\max}}^*$ closest to the center of~$\SX$ (here $(1/2,1/2)$); then $\xb^{(n+1)}$ is the point in $\Xb_{n_{\max}}^*$ furthest away from $\Xb_n=\{\xb^{(1)},\ldots,\xb^{(n)}\}$, for $n=1,\ldots,n_{\max}-1$. The efficiencies of $\Sb_n$, in terms of covering radius (solid line) and packing radius (dashed line), relatively to this new sequence, are plotted in  Figure~\ref{F:filling-theta10-d2-n100}-Right.

\begin{figure}[ht!]
\begin{center}
\includegraphics[width=.49\linewidth]{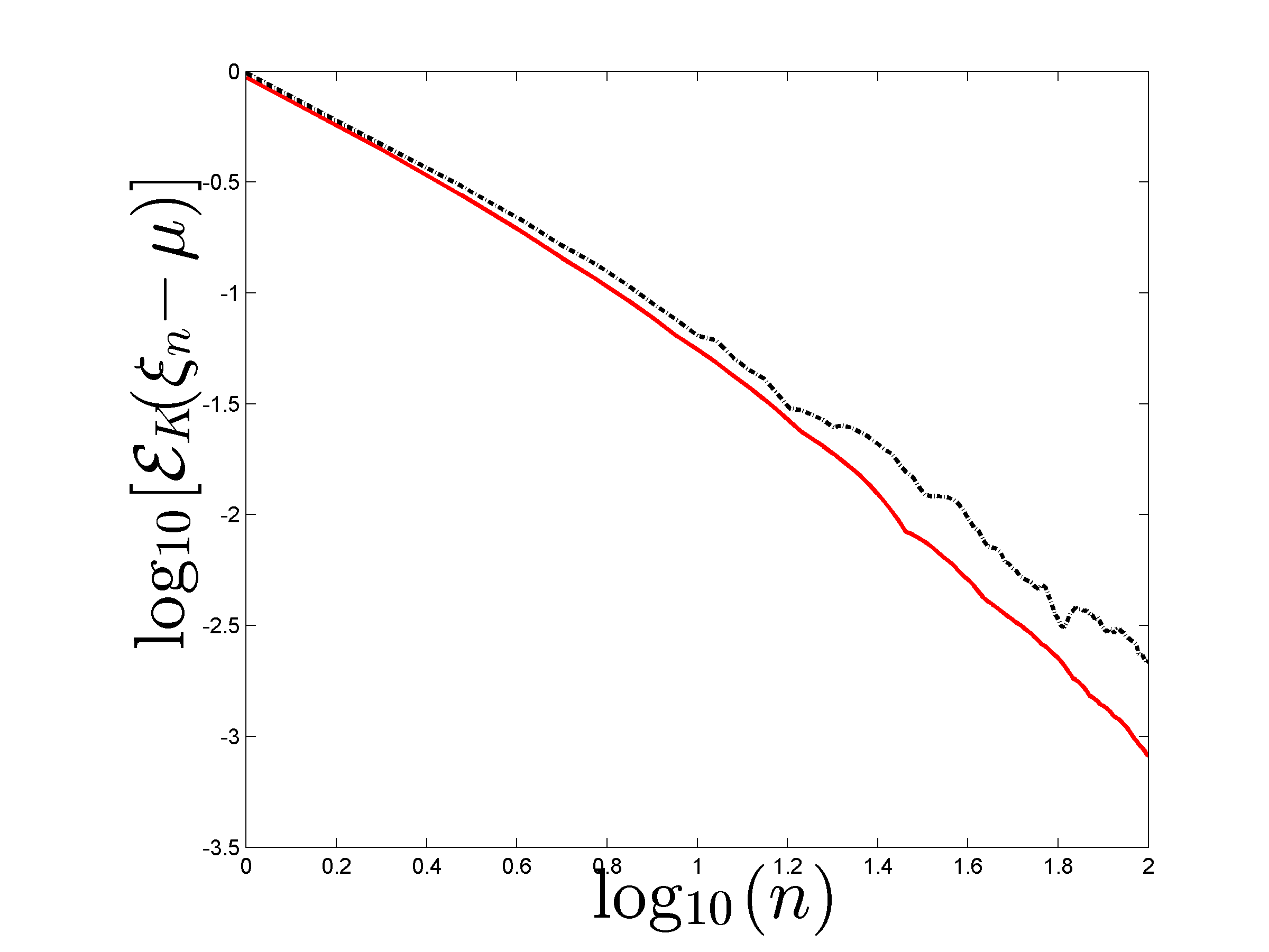} \includegraphics[width=.49\linewidth]{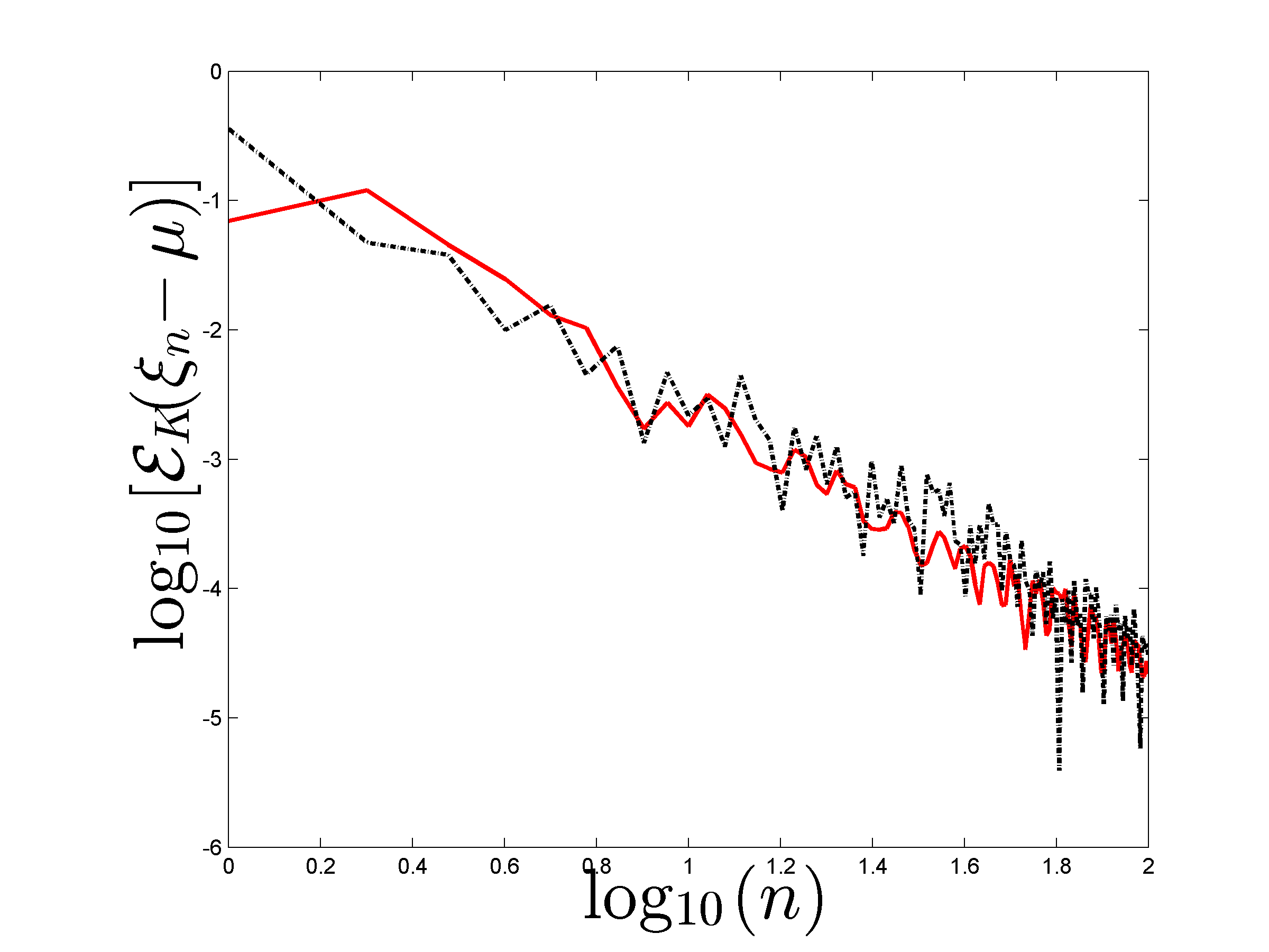}
\end{center}
\caption{\small $\SE_K(\xi^{(n)}-\mu)$ (solid line) and $\SE_K[\xi_{n,e}(\Sb_n)-\mu]$ (dashed line) as functions of $n$ (log scale). $\Xb_n=\Supp(\xi^{(n)})$ is generated by (\ref{x-n+1}, \ref{xi-n+1}) with $\ma_n=1/(n+1)$ and the kernel $K_{3/2,\mt}$; $\Sb_n$ corresponds to the first $n$ points of Sobol' sequence. Left: $\mt=10$ in $K_{3/2,\mt}$; Right: $\mt=1$.}
\label{F:performance-theta10-d2-n100}
\end{figure}

\begin{figure}[ht!]
\begin{center}
\includegraphics[width=.49\linewidth]{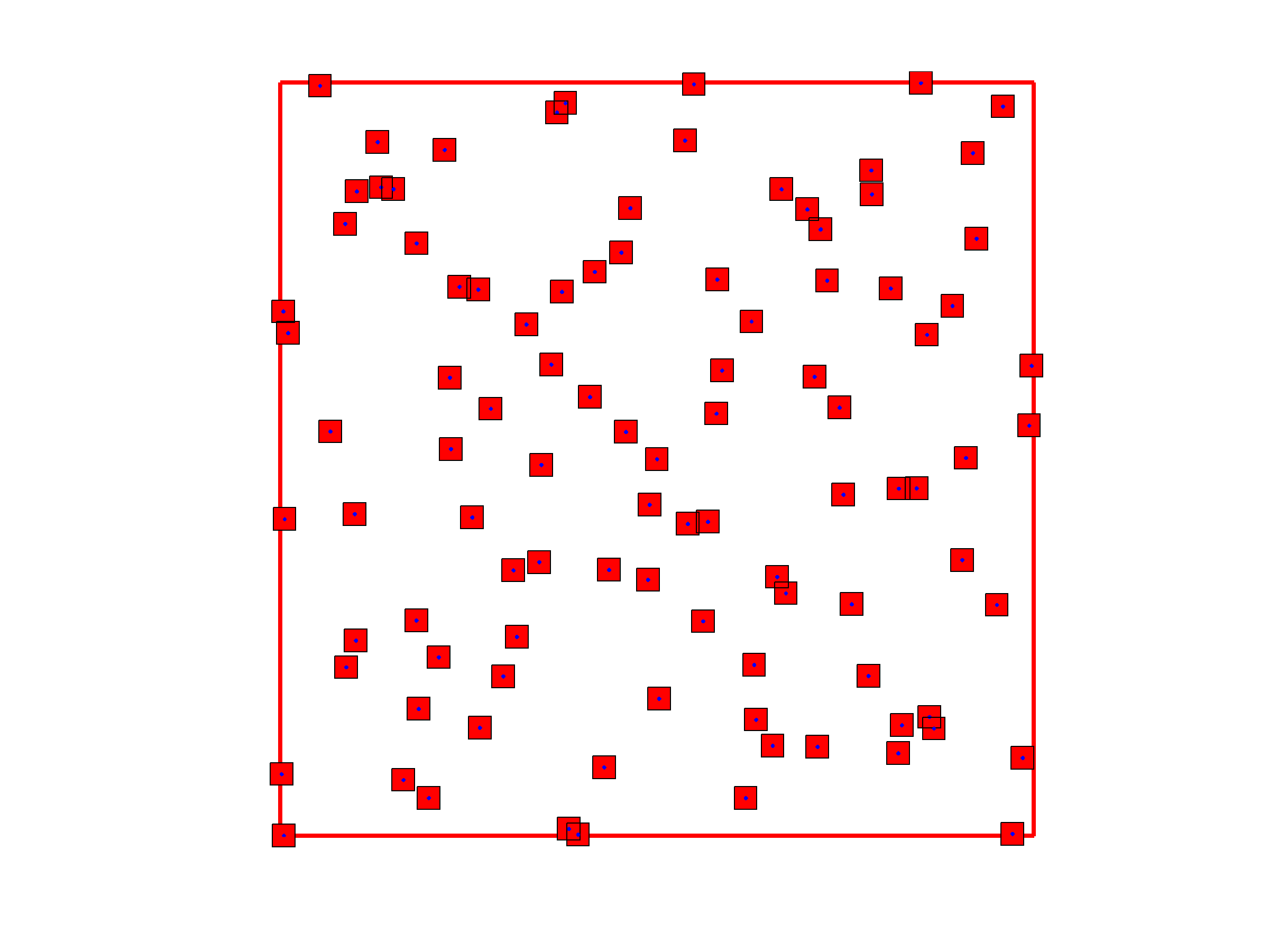} \includegraphics[width=.49\linewidth]{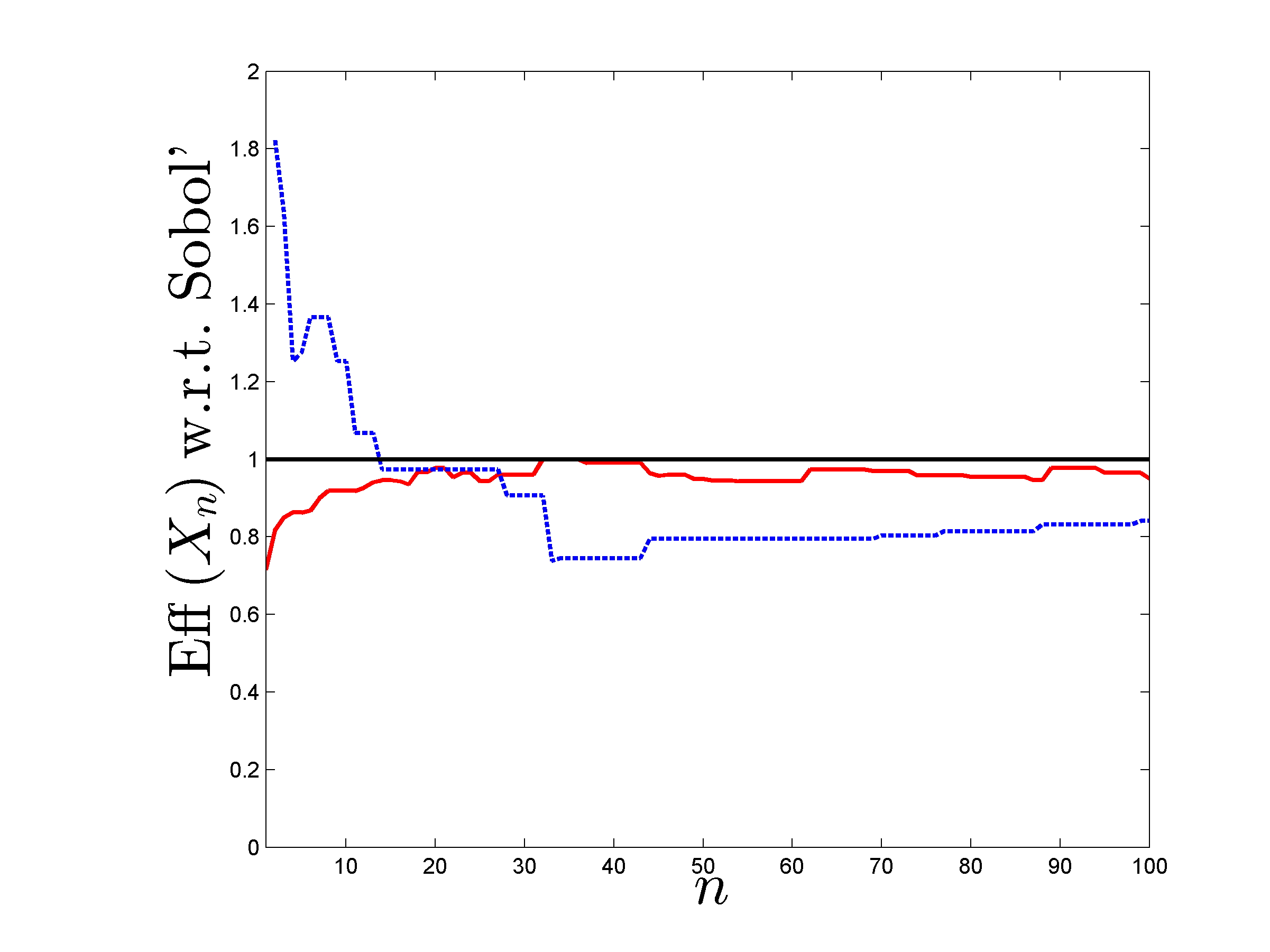}
\end{center}
\caption{\small Left: $\Xb_{100}=\Supp(\xi^{(100)})$; $\SX=[0,1]^2$, $\Xb_n$ is generated by (\ref{x-n+1}, \ref{xi-n+1}) with $\ma_n=1/(n+1)$ and the kernel $K_{3/2,\mt}$ with $\mt=1$. Right: Efficiencies $\CR(\Xb_{n})/\CR(\Sb_{n})$ (solid line) and $\PR(\Sb_{n})/\PR(\Xb_{n})$ (dashed line); $\SX=[0,1]^{10}$, $\Xb_n=\Supp(\xi^{(n)})$ generated by (\ref{x-n+1}, \ref{xi-n+1}) with $\ma_n=1/(n+1)$, and $K(\xb,\xb')=\prod_{i=1}^{10} \log(1/|x_i-x'_i|)$, $\Sb_n$ corresponds to the first $n$ points of Sobol' sequence.}
\label{F:designs-theta1-d2-n100}
\end{figure}

Not surprisingly, the performance of $\xi^{(n)}$ in terms of MMD are better than those of the empirical measure $\xi_{n,e}(\Sb_n)$ associated with $\Sb_n$, see Figure~\ref{F:performance-theta10-d2-n100}-Left.
Figure~\ref{F:performance-theta10-d2-n100}-Right illustrates the fact that a larger correlation length in $K$ yields a faster decrease of $\SE_K[\xi_{n,e}(\Sb_n)-\mu]$ (compare with the figure on the left). On the other hand, this faster decrease does not mean that design points are better distributed, compare Figure~\ref{F:designs-theta1-d2-n100}-Left with Figures~\ref{F:designs-theta10-d2-n100}-Left.

Finally, Figure~\ref{F:designs-theta1-d2-n100}-Right illustrates an application of Algorithm (\ref{x-n+1}, \ref{xi-n+1}) in higher dimension with a singular kernel: $\SX=[0,1]^{10}$ and $K$ is the tensor product of the one-dimensional logarithmic kernel $K_{(0)}$ in \eqref{Riesz}. Similarly to Figure~\ref{F:filling-theta10-d2-n100}, the figure presents the efficiencies $\CR(\Xb_{n})/\CR(\Sb_{n})$ and $\PR(\Sb_{n})/\PR(\Xb_{n})$ for $n=1,\ldots,100$. Due to the large value of $d$, for a given design $\Zb_n$, $\CR(\Zb_n)$ defined by \eqref{CR} is under-approximated by $\CR'(\Zb_n)=\max_{\xb\in\SX'_{\mO'}}\min_{\zb\in\Zb_n} \|\xb-\zb\|$, with $\SX'_{\mO'}$ formed by the first $2^{19}$ points of Sobol' sequence complemented with a $3^d$ full factorial design.

\appendix
\section{Some convergence properties of conditional gradient algorithms}\label{S:CV-properties}

We consider a conditional gradient algorithm with iterations given by \eqref{xi-n+1} when~$\SX$ is replaced by a finite set $\SX_\mO=\{\sb_1,\ldots,\sb_\mO\}$. We do not assume that $\SH_K$ is finite dimensional. A measure $\xi$ on $\SX_\mO$ is characterized by a vector of weights $\mob$ in the probability simplex $\mathds{P}_\mO$, and we shall write $J_K(\mob)=\SE_K(\xi-\mu)$. When $\mob^T\1b_\mO=1$, equation \eqref{JK-sum=1} gives
$$
J_K(\mob)=\|\mob-\hat\mob\|_\Kb^2+J_K(\hat\mob) \,,
$$
where we write $\Kb=\Kb_\mO$, $\|\mob-\hat\mob\|_\Kb^2=(\mob-\hat\mob)^T\Kb(\mob-\hat\mob)$ and where $\hat\mob=\hat\mob_\mO$ is given by \eqref{wopt-sum=1}. We denote by $\ml_{\max}(\Kb)$ the largest eigenvalue of $\Kb$.

For $i=1,\ldots,\mO$, we denote by $\eb_i$ the $i$-th basis vector, with component number $i$ equal to one. Iteration \eqref{xi-n+1} has the form
$$
\mob^{(n+1)}=\mob^{(n)} + \alpha_n \Delta_n
$$
for some step-size $\alpha_n$ and direction $\Delta_n= \eb_{i_n^+}-\mob^{(n)}$, with the index $i_n^+$ taken in $\Arg\min_{i=1,\ldots,\mO} \eb_i^T \nabla J_K(\mob^{(n)})$, where the gradient $\nabla J_K(\mob)$ is given by
$$
\nabla J_K(\mob) = 2 \Kb (\mob-\hat\mob)\,.
$$
This is equivalent to $\sb_{i_n^+} \in\Arg\min_{\sb\in\SX_\mO} \left[ P_{\xi^{(n)}}(\sb) - P_\mu(\sb) \right]$, see \eqref{x-n+1}.

\subsection{Vertex-direction, predefined step-size}\label{A:VD-predefined}

Take $\ma_n=1/(n+1)$ in \eqref{xi-n+1}. We first mention a simple result indicating that $\|\mob^{(n)}-\1b_\mO/\mO\|_\Kb^2=\SO(1/n)$ during the initial $n_1\leq \mO$ iterations when all $i_n^+$ are distinct for $n\leq n_1$.

\begin{lem}\label{L:alpha_k-predefined_b}
Algorithm \eqref{xi-n+1} with $\ma_n=1/(n+1)$, initialized at $\mob^{(1)}=\eb_{i_0}$ for some $i_0\in\{1,\ldots,\mO\}$, satisfies
$$
\|\mob^{(n)}-\1b_\mO/\mO\|_\Kb^2 \leq  \frac{\ml_{\max}(\Kb)}{n}\,, \quad 1\leq n \leq n_1 \leq \mO\,,
$$
where $n_1$ is such that all $i_n^+$ are distinct for $n\leq n_1$.
\end{lem}

\noindent{\em Proof.} For $n\leq n_1$, after a suitable reordering of indices we have $\mob^{(n)}=(1/n,\ldots,1/n,$ $0\ldots,0)^T$. Therefore,
$\|\mob^{(n)}-\1b_\mO/\mO\|_\Kb^2 \leq \ml_{\max}(\Kb) \|\mob^{(n)}-\1b_\mO/\mO\|^2 = \ml_{\max}(\Kb)\, (\mO-n)/(n\mO) \leq \ml_{\max}(\Kb)/n$.
\carre

\vsp
Note that this property is independent of the order in which the vertices of $\mathds{P}_\mO$ (the $\eb_{i_n}$) are selected. It is therefore also valid for MC sampling without replacement within $\SX_\mO$. Also note that the optimal step-size $\hat\ma_n$ at iteration $n$ for the minimization of $\|\mob^{(n)}-\1b_\mO/\mO\|^2$ equals $\ma_n=1/(n+1)$.

The following lemma shows that \eqref{xi-n+1} with $\ma_n=1/(n+1)$ ensures that $\|\mob^{(n)}-\hat\mob\|_\Kb^2=\SO(\log n/n)$, independently of $\mO$ and of the positions of $\mob^{(1)}$ and $\hat\mob$ in $\mathds{P}_\mO$.

\begin{lem}\label{L:alpha_k-predefined}
Algorithm \eqref{xi-n+1} with $\ma_n=1/(n+1)$, initialized at any $\mob^{(1)}$ in $\mathds{P}_\mO$, satisfies
\be\label{rate-alpha_k-predefined}
\|\mob^{(n)}-\hat\mob\|_\Kb^2 \leq 2\, \ml_{\max}(\Kb)\, \frac{1+2\,\log(n+1)}{n}\,, \quad n\geq 1\,.
\ee
\end{lem}

\noindent{\em Proof.} The proof follows the same lines as in \cite[Sect.~3]{Clarkson2010}.
Denote $g(\mob)=\|\mob-\hat\mob\|_\Kb^2$ and $\mob^{(n+)}(\ma)=\mob^{(n)}+\ma\Delta_n$. Notice that $\mob^{(n)}\in\mathds{P}_\mO$ for all $n\geq 1$. We have
\bea
g[\mob^{(n+)}(\ma)] &=& g(\mob^{(n)})+2\ma \Delta_n^T\Kb(\mob^{(n)}-\hat\mob)+\ma^2 \|\Delta_n\|_\Kb^2 \\
&&\leq g(\mob^{(n)})+2\ma \Delta_n^T\Kb(\mob^{(n)}-\hat\mob)+\ma^2 \ml_{\max}(\Kb) \|\Delta_n\|^2 \\
&&\leq g(\mob^{(n)})+2\ma \Delta_n^T\Kb(\mob^{(n)}-\hat\mob)+ 2\,\ma^2 \ml_{\max}(\Kb) \,.
\eea
The convexity of $g(\cdot)$ and the definition of $\Delta_n$ imply that
$$
g(\mob^{(n)}) \geq g(\hat\mob) =0 \geq g(\mob^{(n)})+(\hat\mob-\mob^{(n)})^T\nabla J_K(\mob^{(n)})\geq g(\mob^{(n)})+\Delta_n^T\nabla J_K(\mob^{(n)})\,.
$$
Therefore, $\Delta_n^T\nabla J_K(\mob^{(n)}) = 2 \Delta_n^T\Kb(\mob^{(n)}-\hat\mob)\leq -g(\mob^{(n)})$ and
\be\label{main-alpha_k-predefined}
g[\mob^{(n+)}(\ma)] \leq (1-\ma) g(\mob^{(n)})+ 2\,\ma^2 \ml_{\max}(\Kb) \,.
\ee

The rest of the proof is by induction on $n$.
The bound \eqref{rate-alpha_k-predefined} is valid for $n=1$ since $\|\mob^{(1)}-\hat\mob\|_\Kb^2 \leq 2\,\ml_{\max}(\Kb)$. Suppose that it is satisfied by $\mob^{(n)}$; \eqref{main-alpha_k-predefined} gives
\bea
\|\mob^{(n+1)}-\hat\mob\|_\Kb^2 &\leq& 2\, \ml_{\max}(\Kb)\, \left\{ \frac{1+2\,\log(n+2)}{n+1}- \frac{2(n+1)\log[1+1/(n+1)]-1}{(n+1)^2} \right\} \\
&& \leq 2\, \ml_{\max}(\Kb)\, \frac{1+2\,\log(n+2)}{n+1}
\eea
since $\log(1+t)\geq t/2$ for $t\in[0,1]$.
\carre

\vsp
Using $\ma=2/(n+3)$ in \eqref{main-alpha_k-predefined}, one can easily prove by induction that $g(\mob^{(n)}) \leq 8\,\ml_{\max}(\Kb)/(n+3)$ for all $n$, see \cite{Clarkson2010}, which means that \eqref{xi-n+1} with $\ma_n=2/(n+3)$ instead of $1/(n+1)$ satisfies $\|\mob^{(n)}-\hat\mob\|_\Kb^2 \leq 8\,\ml_{\max}(\Kb)/(n+3)$, $n\geq 1$, with thus a much faster decrease than \eqref{rate-alpha_k-predefined}. Using a different approach, it is shown in \cite{DunnH78} that a rate of decrease of $\SO(1/n)$ is also obtained when $\ma_n$ corresponds to the sequence $\ma_{n+1}=\ma_n-\ma_n^2/2$ with $\ma_0=1$.

Next lemma, based on \cite{ChenWS2012}, shows that $\|\mob^{(n)}-\hat\mob\|_\Kb^2$ decreases in $C/n^2$ when $\hat\mob$ lies in the interior of $\mathds{P}_\mO$. Here, contrary to \cite{ChenWS2012}, we do not assume that $\SH$ is finite dimensional and use instead the finite dimensionality of $\mob$.

\begin{lem}\label{L:alpha_k-predefined-B}
When $\hat\mob$ is in the interior of $\mathds{P}_\mO$, \eqref{xi-n+1} with with $\ma_n=1/(n+1)$, initialized at any $\mob^{(1)}$ in $\mathds{P}_\mO$,
satisfies
\bea\label{rate-alpha_k-predefined_B}
\|\mob^{(n)}-\hat\mob\|_\Kb^2 \leq 4R_*^2 \left(1+ \frac{R_*^2}{\ma_*^2} \right)  \, \frac{1}{n^2}\,, \quad n\geq 1\,,
\eea
where $R_* = [\ml_{\max}(\Kb)\, (1-1/\mO)]^{1/2}$ and $\ma_*=w_*/L$, with $w_*=\min_{i=1,\ldots,\mO} \{\hat\wb\}_i$ (so that $w_* \leq 1/\mO$) and $L=(\max_{i=1,\ldots,\mO}\{\Kb^{-1}\}_{ii})^{1/2}$.
\end{lem}

\noindent{\em Proof.} Denote $\vb(\ma)= \hat\mob - \ma (\mob^{(n)}-\hat\mob)/\|\mob^{(n)}-\hat\mob\|_\Kb$, $\ma>0$. Then, for any $i=1,\ldots,\mO$,
\bea
\frac{|\eb_i^T(\mob^{(n)}-\hat\mob)|}{\|\mob^{(n)}-\hat\mob\|_\Kb} &\leq& \max_{\ub:\, \ub^T\1b_\mO=0,\, \ub^T\Kb\ub=1} |\eb_i^T\ub| \\
&\leq& \max_{\ub:\, \ub^T\Kb\ub=1} |\eb_i^T\ub| = \sqrt{\eb_i^T\Kb^{-1}\eb_i} \leq L = w_*/\ma_* \,,
\eea
so that $\{\vb(\ma)\}_i \geq \{\hat\mob\}_i-w_* \geq 0$, and $\vb(\ma)\in\mathds{P}_\mO$, for any $\ma \leq \ma_*$. The definition of $\eb_{i_n^+}$ then implies that
\be\label{main-alpha_k-predefined-B}
\hspace{-0.5cm} (\eb_{i_n^+}-\hat\mob)^T \Kb (\mob^{(n)}-\hat\mob) \leq [\vb(\ma_*) -\hat\mob)]^T\Kb (\mob^{(n)}-\hat\mob) = - \ma_* \|\mob^{(n)}-\hat\mob\|_\Kb \,.
\ee

The rest of the proof is based on \cite{ChenWS2012}. Denote $\eb_{i_0}^+=\mob^{(1)}$
and $\zb_n=\sum_{i=1}^n (\eb_{i_{n-1}^+}-\hat\mob)$. We can write $\mob^{(n)}=(1/n) \sum_{i=1}^n \eb_{i_{n-1}^+}$, so that
$\zb_n=n(\mob^{(n)}-\hat\mob)$, $n^2 \|\mob^{(n)}-\hat\mob\|_\Kb^2 = \|\zb_n\|_\Kb^2$, and we only need to bound $\|\zb_n\|_\Kb^2$. We have
$$
\|\zb_n\|_\Kb^2-\|\zb_{n+1}\|_\Kb^2 = -2 (\eb_{i_n^+}-\hat\mob)^T\Kb\zb_n - \|\eb_{i_n^+}-\hat\mob\|_\Kb^2 \,,
$$
where $\|\eb_{i_n^+}-\hat\mob\|_\Kb \leq 2\, R_*$ and
$(\eb_{i_n^+}-\hat\mob)^T\Kb\zb_n \leq - \ma_* \|\zb_n\|_\Kb$ from \eqref{main-alpha_k-predefined-B}. Therefore,
$$
\|\zb_{n+1}\|_\Kb^2 \leq \|\zb_n\|_\Kb^2- 2\ma_* \left(\|\zb_n\|_\Kb - 2 R_*^2/\ma_* \right) \,.
$$
Suppose that $\|\zb_n\|_\Kb>2 R_*^2/\ma_*$. Then, $\|\zb_{n+1}\|_\Kb^2 \leq \|\zb_n\|_\Kb^2$, and $\|\zb_n\|_\Kb^2$ decreases until some $n_0$ when $\|\zb_{n_0}\|_\Kb \leq 2 R_*^2/\ma_*$. But then,
\bea
\|\zb_{n_0+1}\|_\Kb^2 \leq \|\zb_{n_0}\|_\Kb^2- 2\ma_* \left(\|\zb_{n_0}\|_\Kb - 2 R_*^2/\ma_* \right) \leq 4 R_*^2(1+R_*^2/\ma_*^2) \,,
\eea
so that $\|\zb_n\|_\Kb^2 \leq 4 R_*^2(1+R_*^2/\ma_*^2)$ for all $n>n_0$.
\carre

\vsp
Lemma~\ref{L:alpha_k-predefined-B} indicates that $\|\mob^{(n)}-\hat\mob\|_\Kb^2 \leq C/n^2$. However, for large $\mO$ the constant $C$ grows like $\SO(\mO^2)$
(since $\ma_*\leq 1/(L\mO)$) which makes this result of theoretical interest only. Note that kernel herding typically concerns situations where $\mO$ is very large (or even infinite when~$\SX$ is not discretized).

\subsection{Vertex-direction, optimal step-size}
The choice of a predefined step-size $\ma_n=1/(n+1)$ in \eqref{xi-n+1} does not ensure a monotonic decrease of $\SE_K(\xi^{(n)}-\mu)$. An alternative option is to choose $\ma_n$ that minimizes $\SE_K[\xi^{(n+)}(\ma)-\mu]$ with respect to $\ma\in[0,1]$, with $\xi^{(n+)}(\ma)=(1-\ma)\, \xi^{(n)} + \ma\, \delta_{\xb_{n+1}}$ and $\xb_{n+1}$ given by \eqref{x-n+1}. Straightforward calculation gives $\ma_n = \min\{1,\hat\ma_n\}$, with
\be\label{ma_n}
\hat\ma_n &=& \frac{\langle P_{\xi^{(n)}}-P_\mu, P_{\xi^{(n)}}-P_{\delta_{\xb_{n+1}}}\rangle_K}{\|P_{\xi^{(n)}}-P_{\delta_{\xb_{n+1}}}\|_{\SH_K}^2} \\
&=& \frac{ \SE_K(\xi^{(n)}) - P_{\xi^{(n)}}(\xb_{n+1}) - \sum_{i=1}^n w_i^{(n)} P_\mu(\xb_i) + P_\mu(\xb_{n+1})}{\SE_K(\xi^{(n)})-2P_{\xi^{(n)}}(\xb_{n+1}) + K(\xb_{n+1},\xb_{n+1})  } \nonumber 
\ee
(which requires that $\SE_K(\xi^{(n)})<\infty$).
Direct calculation shows that $\hat\ma_n$ given by \eqref{ma_n} satisfies
\be\label{ma_k-B}
\hat\ma_n = \frac{(\eb_{i_n^+}-\mob^{(n)})^T\Kb(\hat\mob-\mob^{(n)})}{\|\eb_{i_n^+}-\mob^{(n)}\|_\Kb^2} \,.
\ee
Next Lemma indicates that $\hat\ma_n\leq 1$ when $\hat\wb\in\mathds{P}_\mO$, so that setting $\ma_n=\hat\ma_n$ in \eqref{xi-n+1} ensures that $\mob^{(n)}$ remains in $\mathds{P}_\mO$ for all $n$. It should be noticed that the global decrease of $\|\mob^{(n)}-\hat\mob\|_\Kb^2$ over \emph{many} iterations with this optimal $\ma_n$ is not necessarily better that with the predefined step-size $\ma_n=1/(n+1)$ of Section~\ref{A:VD-predefined}; one may refer in particular to \cite{DunnH78} for such considerations; see also \cite{BachLJO2012}.

\begin{lem}\label{L:alpha_k<1}
When $\hat\mob\in\mathds{P}_\mO$, $\hat\ma_n$ given by \eqref{ma_k-B} is less than one.
\end{lem}

\noindent{\em Proof.}
We can write $(\eb_{i_n^+}-\mob^{(n)})^T\Kb(\hat\mob-\mob^{(n)})=\|\eb_{i_n^+}-\mob^{(n)}\|_\Kb^2+
(\eb_{i_n^+}-\hat\mob)^T\Kb(\mob^{(n)}-\hat\mob)-\|\eb_{i_n^+}-\hat\mob\|_\Kb^2$. When $\hat\mob\in\mathds{P}_\mO$, $\hat \mo_i \geq 0$ for all $i$, and $\sum_{i=1}^\mO \hat \mo_i (\eb_i-\hat\mob)=\0b$ implies that $\sum_{i=1}^\mO \hat \mo_i (\eb_i-\hat\mob)^T\Kb(\mob^{(n)}-\hat\mob)=0$. Therefore
$\min_{i=1,\ldots,\mO}(\eb_i-\hat\mob)^T\Kb(\mob^{(n)}-\hat\mob)=(\eb_{i_n^+}-\hat\mob)^T\Kb(\mob^{(n)}-\hat\mob)\leq 0$, which gives $\hat\ma_n \leq 1$.
\carre

\begin{lem}\label{L:alpha_k-opt}
Algorithm \eqref{xi-n+1} with $\ma_n=\hat\ma_n$ given by \eqref{ma_k-B}, initialized at any $\mob^{(1)}\in\mathds{P}_\mO$, satisfies
\be\label{rate-alpha_k-opt}
\|\mob^{(n)}-\hat\mob\|_\Kb^2 \leq 8\, \ml_{\max}(\Kb)\, \frac{1}{n+3}\,, \quad n\geq 1\,.
\ee
\end{lem}

\noindent{\em Proof.} The proof follows \cite[Sect.~2]{Clarkson2010} and uses the same notation as in the proof of Lemma~\ref{L:alpha_k-predefined}. The right-hand side of \eqref{main-alpha_k-predefined} is minimum for $\hat\ma=g(\mob^{(n)})/[4\,\ml_{\max}(\Kb)] \leq 1/2$. Therefore,
\bea
g(\mob^{(n+1)}) = \min_{\ma\in[0,1]} g[\mob^{(n+)}(\ma)] &\leq& (1-\hat\ma) g(\mob^{(n)})+ 2\,\hat\ma^2 \ml_{\max}(\Kb) \\
&&= g(\mob^{(n)})\left[1-\frac{g(\mob^{(n)})}{8\,\ml_{\max}(\Kb)}\right] \,.
\eea
Since $1-t\leq 1/(1+t)$ for all $t>-1$, we obtain
$$
g(\mob^{(n+1)})/[8\,\ml_{\max}(\Kb)] \leq \frac{1}{1+\{g(\mob^{(n)})/[8\,\ml_{\max}(\Kb)]\}^{-1}}
$$
which, by induction, implies that $g(\mob^{(n)}) \leq 8\, \ml_{\max}(\Kb)/(n+3)$; that is, \eqref{rate-alpha_k-opt}.
\carre

\begin{lem}\label{L:alpha_k-opt-B}
When $\hat\mob$ is in the interior of $\mathds{P}_\mO$, \eqref{xi-n+1} with $\ma_n=\hat\ma_n$ given by \eqref{ma_k-B}, initialized at any $\mob^{(1)}\in\mathds{P}_\mO$, satisfies
\be\label{rate-alpha_k-opt_B}
\|\mob^{(n+1)}-\hat\mob\|_\Kb^2 \leq \|\mob^{(1)}-\hat\mob\|_\Kb^2 \, \exp\left(-\frac{\ma_*^2\,k}{4R_*^2}\right)\,, \quad n\geq 1\,,
\ee
where $R_* = [\ml_{\max}(\Kb)\, (1-1/\mO)]^{1/2}$ and $\ma_*=w_*/L$, with $w_*=\min_{i=1,\ldots,\mO} \{\hat\wb\}_i$ (so that $w_* \leq 1/\mO$) and $L=(\max_{i=1,\ldots,\mO}\{\Kb^{-1}\}_{ii})^{1/2}$.
\end{lem}

\noindent{\em Proof.} We use the same approach as in \cite{BeckT2004} and use the same notation as in the proof of Lemma~\ref{L:alpha_k-predefined}. We can write
$g(\mob^{(n+1)})=g[\mob^{(n+)}(\hat\ma_n)]$, with $\hat\ma_n$ given by \eqref{ma_k-B}. Therefore,
\bea
g(\mob^{(n+1)}) &=& g(\mob^{(n)}) + 2 \hat\ma_n \Delta_n^T \Kb (\mob^{(n)}-\hat\mob)+ \ma_n^2 \|\Delta_n\|_\Kb^2 \\
&=& g(\mob^{(n)}) - \frac{[\Delta_n^T \Kb (\mob^{(n)}-\hat\mob)]^2}{\|\Delta_n\|_\Kb^2} \,.
\eea
Equation \eqref{main-alpha_k-predefined-B} implies that $[\Delta_n^T \Kb (\mob^{(n)}-\hat\mob)]^2 \geq \ma_*^2 g(\mob^{(n)})$, and thus
$$
g\mob^{(n+1)}) \leq g(\mob^{(n)}) \left[1 - \frac{\ma_*^2}{\|\Delta_n\|_\Kb^2} \right] \leq g(\mob^{(n)}) \left[1 - \frac{\ma_*^2}{4R_*^2} \right]\,.
$$
This implies $g(\mob^{(n+1)}) \leq g(\mob^{(1)})\, \exp[-\ma_*^2\,k /(4R_*^2)]$, that is, \eqref{rate-alpha_k-opt_B}.
\carre

\vsp
Similarly to Lemma~\ref{L:alpha_k-predefined-B}, the small value of the constant $\ma_*$ makes the linear convergence rate in \eqref{rate-alpha_k-opt_B} of theoretical interest only.

\subsection{Vertex-exchange}

Following \cite{MolchanovZ01, MolchanovZ02}, one may also use a vertex-exchange method based on the true steepest-descent direction, see also \cite{Bohning85, Bohning86}. The iterations are then
\be\label{iter3}
\xi^{(n+1)}=\xi^{(n)} + \ma_n\, (\delta_{\xb_{n+1}} - \delta_{\xb_n^-})\,,
\ee
where $\xb_{n+1}$ is given by \eqref{x-n+1} and
\be\label{xn-}
\xb_n^- \in \Arg\max_{\xb\in\Supp(\xi^{(n)})} \left[ P_{\xi^{(n)}}(\xb) - P_\mu(\xb) \right]\,,
\ee
with $\Supp(\xi^{(n)})=\Xb_n$ the support of $\xi^{(n)}$.
The step-size $\ma_n$ is then given by $\min\{\hat\ma_n,\xi^{(n)}(\xb_n^-)\}$, where
$\hat\ma_n$ minimizes $\SE_K[\{\xi^{(n)} + \ma\, (\delta_{\xb_{n+1}} - \delta_{\xb_n^-})\}-\mu]$ with respect to $\ma$ (the constraint $\ma_n\leq \xi^{(n)}(\xb_n^-)$ ensures that $\xi^{(n+1)}\in\SM^+(1)$ when $\xi^{(n)}\in\SM^+(1)$). Direct calculation gives
\be
\hat\ma_n &=& \frac{\langle P_{\xi^{(n)}}-P_\mu, P_{\delta_{\xb_n^-}}-P_{\delta_{\xb_{n+1}}}\rangle_K}{\|P_{\delta_{\xb_n^-}}-P_{\delta_{\xb_{n+1}}}\|_{\SH_K}^2} \nonumber\\
&=& \frac{ [P_{\xi^{(n)}}(\xb_n^-)-P_{\mu}(\xb_n^-)] - [P_{\xi^{(n)}}(\xb_{n+1})-P_{\mu}(\xb_{n+1})]}{K(\xb_n^-,\xb_n^-)+K(\xb_{n+1},\xb_{n+1})-2K(\xb_n^-,\xb_{n+1})}\,. \label{alpha_n3}
\ee

For the algorithm defined by (\ref{iter3}, \ref{xn-}), we have $\mob^{(n+1)}=\mob^{(n)}+\ma_n \Delta_n$ with now $\Delta_n=\eb_{i_n^+}-\eb_{i_n^-}$, where we take $i_n^+\in\Arg\min_{i=1,\ldots,\mO} \eb_i^T \Kb(\mob^{(n)}-\hat\mob)$ and $i_n^-\in\Arg\max_{i:\, \eb_i^T\mob^{(n)}>0} \eb_i^T \Kb(\mob^{(n)}-\hat\mob)$. The step size \eqref{alpha_n3} equals
\be\label{ma_k-VE}
\hat\ma_n = \frac{\Delta_n^T\Kb(\hat\mob-\mob^{(n)})}{\|\Delta_n\|_\Kb^2} \,.
\ee
Take $\ma_n=\min\{\mo^{(n)}_{i_n^-},\hat\ma_n\}$ in \eqref{iter3}, so that $\mob^{(n)}$ remains in $\mathds{P}_\mO$ for all $n$. Using the same notation as in the proof of Lemma~\ref{L:alpha_k-predefined}, we have
$$
g[\mob^{(n+)}(\ma)] \leq g(\mob^{(n)})+2\ma \Delta_n^T\Kb(\mob^{(n)}-\hat\mob)+ 2\,\ma^2 \ml_{\max}(\Kb) \,,
$$
and, since $\mob^{(n)}\in\mathds{P}_\mO$, the convexity of $g(\cdot)$ and the definition of $\Delta_n$ imply that
$$
g(\hat\mob) =0 \geq g(\mob^{(n)})+ 2(\hat\mob-\mob^{(n)})^T \Kb(\mob^{(n)}-\hat\mob)\geq g(\mob^{(n)})+2\Delta_n^T \Kb(\mob^{(n)}-\hat\mob)\,.
$$
We obtain the following property; the proof is identical to that of Lemma~\ref{L:alpha_k-opt}.

\begin{lem}\label{L:VE-alpha_k-predefined}
Suppose that $\hat\mob$ and $\Kb$ are such that $\hat\ma_n \leq \mo^{(n)}_{i_n^-}$ for any $\mob^{(n)}\in\mathds{P}_\mO$. Then, algorithm (\ref{iter3}, \ref{xn-}) with $\ma_n=\hat\ma_n$ given by \eqref{ma_k-VE}, initialized at any $\mob^{(1)}\in\mathds{P}_\mO$, satisfies
\bea\label{VE-rate-alpha_k-predefined}
\|\mob^{(n)}-\hat\mob\|_\Kb^2 \leq \frac{8\, \ml_{\max}(\Kb)}{n+3}\,, \quad n\geq 1\,.
\eea
\end{lem}

There exist situations where the condition $\hat\ma_n \leq \mo^{(n)}_{i_n^-}$ is not satisfied. Take for instance $\mO=3$, $\Kb$ the identity matrix and $\hat\mob=(0,0,1)^T$, $\mob^{(n)}=(1/3,1/3,1/3)^T$; then $\hat\ma_n=1/2 > \mo^{(n)}_{i_n^-}=1/3$). On the other hand, the condition is satisfied for instance for $\hat\mob=\1b_\mO/\mO$ and $\Kb$ the identity matrix (we have $i_{n^-}=\Arg\max_{i:\,\mo^{(n)}_i>0} (\mo^{(n)}_i-\hat \mo_i)$, and $\sum_{i=1}^\mO (\mo^{(n)}_i-\hat \mo_i)=0$ implies that $\mo^{(n)}_{i_n^-}>\hat \mo_{i_n^-}$ and similarly $\mo^{(n)}_{i_n^+}<\hat \mo_{i_n^+}$; we get $\hat \ma_n=(\mo^{(n)}_{i_n^-}-\mo^{(n)}_{i_n^+})/2 \leq \mo^{(n)}_{i_n^-}/2 < \mo^{(n)}_{i_n^-}$), and numerical experiments indicate that it holds true in most situations.

\section{Bayesian quadrature: several integrals}\label{S:BQ-several}

Following \cite{O'Hagan91}, consider a generalization of the situation considered in Section~\ref{S:BQ} where one wishes to estimate
$$
\Ib_\mu(f)=\Ex_\mu\{f(\Xb)\rb(\Xb)\}=\int_\SX f(\xb)\rb(\xb)\, \dd\mu(\xb) \,,
$$
with $\rb(\xb)=(r_0(\xb),\ldots,r_p(\xb))^T$ a vector of $p+1$ known functions of $\xb$, such that the $(p+1)\times (p+1)$ matrix
$$
\Mb_r = \Ex_\mu\{\rb(\Xb)\rb^T(\Xb)\}
$$
exists and is nonsingular. Without any loss of generality, we may assume that $r_0(\xb)\equiv~1$.

We also slightly generalize model \eqref{model1} by introducing a linear trend $\hb^T(\xb)\betab$; that is, we consider
\be \label{model2}
f(\xb)= \hb^T(\xb)\betab +Z_x \,,
\ee
where $\hb(\xb)=(h_0(\xb),\ldots,h_{p'}(\xb))^T$ is a vector of $p'+1$ known functions of $\xb$ and $\betab\in\mathds{R}^{p'+1}$ has the normal prior $\SN(\betabh^0,\ms^2 \Ab)$, non-informative so that we can replace $\Ab^{-1}$ by the null matrix $\0b$ in all calculations (the choice of $\betabh^0$ being then irrelevant). We assume that the matrix $\Ex_\mu\{\hb(\Xb)\hb^T(\Xb)\}$ is well defined. For reasons that will become clear below, we shall consider in particular the case where $\hb=\rb$.

The posterior mean and variance of $f(\xb)$, conditional on $\ms^2$ and $K$, are now, respectively,
\bea
  \hat \eta_n(\xb) &=& \hb^T(\xb)\betabh^n + \kb_n^T(\xb)\Kb_n^{-1}(\yb_n- \Hb_n\betabh^n)\,, \nonumber  \\
  \ms^2\rho_n^2(\xb)&=&\ms^2 \left\{K(\xb,\xb)-\kb_n^T(\xb)\Kb_n^{-1}\kb_n(\xb) \right. \nonumber \\
&&  \left. + \, [\hb(\xb)-\Hb_n^T\Kb_n^{-1}\kb_n(\xb)]^T
(\Hb_n^T\Kb_n^{-1}\Hb_n)^{-1}[\hb(\xb)-\Hb_n^T\Kb_n^{-1}\kb_n(\xb)]
\right\}\,, \label{MSE2}
\eea
where $\{\Hb_n(\xb)\}_{i,j}=h_j(\xb_i)$, $i=1,\ldots,n$, $j=0,\ldots,p'$, and
$$
\betabh^n = (\Hb_n^T\Kb_n^{-1}\Hb_n)^{-1}\Hb_n^T\Kb_n^{-1}\yb_n \,.
$$
The posterior mean and covariance matrix of $\Ib_\mu(f)$ are
\be
\widehat\Ib_n &=& \Bb(\mu)\betabh^n + \Pb_n(\mu)\Kb_n^{-1}(\yb_n- \Hb_n\betabh^n)\,, \label{Ibn1} \\
\ms^2\,\Vb_n &=& \ms^2 \left\{ \Ub(\mu) - \Pb_n(\mu)\Kb_n^{-1}\Pb_n^T(\mu) \right. \nonumber\\
&& \hspace{-1.5cm} \left. + \left[\Bb(\mu)-\Pb_n(\mu)\Kb_n^{-1}\Hb_n\right] (\Hb_n^T\Kb_n^{-1}\Hb_n)^{-1} \left[\Bb(\mu)-\Pb_n(\mu)\Kb_n^{-1}\Hb_n\right]^T \right\}\,, \label{sbn1}
\ee
where $\Bb(\mu)=\Ex_\mu\{\rb(\Xb)\hb^T(\Xb)\}$, $\Pb_n(\mu)=\Ex_\mu\{\rb(\Xb)\kb_n^T(\Xb)\}$ and
\be\label{U}
\Ub(\mu)=\Ex_\mu\{\rb(\Xb)\rb^T(\Xb')K(\Xb,\Xb')\}\,,
\ee
with $\Xb$ and $\Xb'$ i.i.d.\ $\sim$ $\mu$.

Consider now the special choice $\hb=\rb$ in \eqref{model2}. Following Section~\ref{S:BLUE-Kreduction}, we can write
$f(\xb)= \rb^T(\xb)\betab + {\mathcal{P}_r}Z_x + (\Id_{L^2}-{\mathcal{P}_r})Z_x$,
where ${\mathcal{P}_r}$ denotes the orthogonal projection of $L^2(\SX,\mu)$ onto the linear space spanned by $\rb(\cdot)$; that is, ${\mathcal{P}_r} g(\xb) = \rb^T(\xb) \Mb_r^{-1} \int_\SX \rb(\xb') g(\xb') \, \dd\mu(\xb')$ for all $g\in L^2(\SX,\mu)$. This gives
$$
{\mathcal{P}_r}Z_x = \rb^T(\xb) \Mb_r^{-1} \int_\SX \rb(\xb') Z_{x'} \, \dd\mu(\xb') \,.
$$
In absence of prior information on $\betab$ ($\Ab^{-1}=\0b$), the prior on the parameters $\betab'=\betab + \Mb_r^{-1} \int_\SX \rb(\xb') Z_{x'} \, \dd\mu(\xb')$ remains non-informative, and the covariance kernel of $\widetilde Z_x= (\Id_{L^2}-{\mathcal{P}_r})Z_x$ is
\bea
K_\mu(\xb,\xb')&=& K(\xb,\xb') - \ub_\mu^T(\xb)\Mb_r^{-1}\rb(\xb') - \rb^T(\xb)\Mb_r^{-1}\ub_\mu(\xb') \nonumber \\
&& + \, \rb^T(\xb)\Mb_r^{-1}\Ub(\mu)\Mb_r^{-1}\rb(\xb')\,, \label{tildeK2}
\eea
where $\Ub(\mu)$ is given by \eqref{U} and
$\ub_\mu(\xb)=\Ex_\mu\{\rb(\Xb)K(\Xb,\xb)\}\,, \ \xb\in\SX$.

Similarly to Remark~\ref{R:reducedK} (see \cite[Sect.~5.4]{GP-CSDA2016}), this kernel reduction does not modify predictions, and
direct calculation shows that $\Ex_\mu\{\rb(\Xb)\rb^T(\Xb')K_\mu(\Xb,\Xb')\}=\0b$ and $\Ex_\mu\{\rb(\Xb)\widetilde\kb_n^T(\Xb)\}=\0b$.
We thus obtain the following property, where we denote by $\Rb_n$ the $n\times (p+1)$ matrix $\{\Rb_n(\xb)\}_{i,j}=r_j(\xb_i)$, $i=1,\ldots,n$, $j=0,\ldots,p$.

\begin{lem}\label{L:L2}
When $K$ is SPD, $\widehat \Ib_n$ given by \eqref{Ibn1} satisfies
\bea \label{Ibn2}
\widehat \Ib_n = \Mb_r (\Rb_n^T\widetilde\Kb_n^{-1}\Rb_n)^{-1}(\Rb_n^T\widetilde\Kb_n^{-1}\yb_n) \,,
\eea
and the posterior covariance matrix \eqref{sbn1} satisfies
\be
\Vb_n = \Mb_r (\Rb_n^T\widetilde\Kb_n^{-1}\Rb_n)^{-1} \Mb_r \,. \label{sbn2}
\ee
\end{lem}

To ensure a precise estimation of $\Ib_\mu(f)$, we may select a design $\Xb_n$ that minimizes $\SJ(\Vb_n)$, with $\SJ(\cdot)$ a Loewner increasing function defined on the set of symmetric non-negative define matrices. Typical choices are $\SJ(\Vb_n)=\det(\Vb_n)$ (D-optimality) and $\SJ(\Vb_n)=\tr(\Vb_n)$ (A-optimality).
Greedy minimization of $\SJ(\Vb_n)$ corresponds to Sequential Bayesian Quadrature, see Section~\ref{S:any-time}.
Using \eqref{sbn2} and formulae for the inversion of a block matrix, we obtain the following expressions for $\det(\Vb_{n+1})$ and $\tr(\Vb_{n+1})$:
\bea
\det(\Vb_{n+1}) &=& \det(\Vb_n)\, \frac{K_\mu(\xb,\xb)-\widetilde\kb_n^T(\xb)\widetilde\Kb_n^{-1}\widetilde\kb_n(\xb)}{\widetilde\rho_n^2(\xb)} \,, \\
\tr(\Vb_{n+1}) &=& \tr(\Vb_n) \\
&& \hspace{-2.2cm} - \, \frac{[\rb(\xb)-\Rb_n^T\widetilde\Kb_n^{-1}\widetilde\kb_n(\xb)]^T
(\Rb_n^T\widetilde\Kb_n^{-1}\Rb_n)^{-1}\Mb_r^2 (\Rb_n^T\widetilde\Kb_n^{-1}\Rb_n)^{-1}[\rb(\xb)-\Rb_n^T\widetilde\Kb_n^{-1}\widetilde\kb_n(\xb)]}{\widetilde\rho_n^2(\xb)} \,,
\eea
with
$$
\widetilde\rho_n^2(\xb) = \left[K_\mu(\xb,\xb)-\widetilde\kb_n^T(\xb)\widetilde\Kb_n^{-1}\widetilde\kb_n(\xb)+\frac{(1-\widetilde\kb_n^T(\xb)\widetilde\Kb_n^{-1}\1b_n)^2}{\1b_n^T\widetilde\Kb_n^{-1}\1b_n} \right]\,.
$$
When $p=0$ ($\rb(\xb)\equiv 1$), $\Vb_n=s_n^2$ in \eqref{sn2} and $\det(\Vb_{n+1})=\tr(\Vb_{n+1})=s_{n+1}^2$ given by \eqref{sn+1^2}.

\section{Karhunen-Lo\`eve decomposition and Bayesian c-optimal design}\label{S:KL}

As shown in Section~\ref{S:BLUE-Kreduction}, estimation of $I_\mu(f)$ is equivalent to estimation of $\beta_0'$ in the model \eqref{model1p}. Following \cite{GP-CSDA2016}, we may the consider the Karhunen-Lo\`eve decomposition of the Gaussian RF $\widetilde Z_x$, which has zero mean and covariance $\Exx\{\widetilde Z_x\widetilde Z_{x'}\}=\ms^2\,K_\mu(\xb,\xb')$ and is such that $\int_\mu Z_x\, \dd\mu(\xb)=0$ a.s.; see also \cite{Fedorov96}. We thus write
$f(\xb) = \beta_0 + \sum_{k\geq 1} \beta_k \varphi_k(\xb) + \widetilde Z_{x,0}$,
where the $\varphi_k$ are orthonormal (in $L^2(\SX,\mu)$) eigenfunctions of the spectral decomposition of the integral operator $T_\mu$ defined on $L^2(\SX,\mu)$ by
$T_\mu[f](\xb)= \int_\SX K_\mu(\xb,\xb')f(\xb')\, \dd\mu(\xb')$, $f \in L^2(\SX,\mu)$, $\xb\in\SX$,
with associated eigenvalues $\Lambda_k$ (strictly positive and assumed to be ordered by decreasing values), and where the $\beta_k$, $k\geq 1$, are independent random variables distributed $\SN(0,\ms^2\Lambda_k)$. Here $\widetilde Z_{x,0}$ is a zero-mean Gaussian RF with covariance $\Exx\{\widetilde Z_{x,0}\widetilde Z_{x',0}\}=\ms^2\,\left[K_\mu(\xb,\xb')-\sum_{k\geq 1} \Lambda_k \varphi_k(\xb)\varphi_k(\xb')\right]$. Spectral truncation at $k=M-1$ yields the approximated model
\be\label{BLM1}
f(\xb) \approx \beta_0' + \sum_{k=1}^{M-1} \beta_k \varphi_k(\xb) + \widetilde Z_{x,1} \,,
\ee
where $\widetilde Z_{x,1}$ has zero mean and satisfies
$$
\Exx\{\widetilde Z_{x,1}\widetilde Z_{x',1}\}=\ms^2 \left[K_\mu(\xb,\xb')-\sum_{k=1}^{M-1} \Lambda_k \varphi_k(\xb)\varphi_k(\xb')\right] \,.
$$
The eigenfunctions $\varphi_k$ satisfy $\int_\SX \varphi_k(\xb)\,\dd\mu(\xb)=0$ for all $k$. When evaluating $f$ at $\xb$ outside the support of $\mu$ we replace all $\varphi_k$ in \eqref{BLM1} by their canonical extensions $\varphi_k'$ defined by
$\varphi_k'(\xb)=(1/\Lambda_k) T_\mu[\varphi_k](\xb)$, for all $k\geq 1$ and all $\xb$ in~$\SX$; they also satisfy $\int_\SX \varphi_k'(\xb)\,\dd\mu(\xb)=0$ for all $k$. Note that the $\varphi_k'$ are defined over the whole~$\SX$, with $\varphi_k'=\varphi_k$ $\mu$-almost everywhere, which is important when $\mu$ is only supported on a subset of~$\SX$.

The errors $\widetilde Z_{x,1}$ in \eqref{BLM1} are correlated. To facilitate the construction of experimental designs ensuring a precise estimation of $\beta_0'$, an additional approximation is introduced which yields a Bayesian Linear Model (BLM) with $M$ parameters $\betab=(\beta_0',\beta_1,\ldots,\beta_{M-1})^T$ and uncorrelated errors
\be\label{BLM2}
f(\xb) \approx \beta_0' + \sum_{k=1}^{M-1} \beta_k \varphi_k(\xb) + \widetilde Z'_{x,1} \,.
\ee
Here $\widetilde Z'_{x,1}$ has zero mean and satisfies $\Exx\{\widetilde Z_{x,1}\widetilde Z_{x',1}\}=\ms^2 v(\xb)$ if $\xb'=\xb$ and is zero otherwise, with
$$
v(\xb) =\left[K_\mu(\xb,\xb)-\sum_{k=1}^{M-1} \Lambda_k \varphi_k(\xb)\varphi_k(\xb)\right] \,.
$$
After observation of $\yb_n=[f(\xb_1),\ldots,f(\xb_n)]^T$, $\betab$ has the posterior normal distribution
$\SN(\hat\betab^n,\ms^2\,\Mb_B^{-1}(\Xb_n))$, with
$$
\Mb_B(\Xb_n)= \Psib_n^T \Sigmab_n^{-1} \Psib_n + \Lambdab^{-1} \,,
$$
where $\{\Psib_n\}_{i,0}=1$ and $\{\Psib_n\}_{i,k}=\varphi_k(\xb_i)$, $k=1,\ldots,M-1$, $i=1,\ldots,n$, $\Sigmab_n=\diag\{v(\xb_i),\, i=1,\ldots,n\}$, $\Lambdab = \diag\{\infty, \Lambda_1,\ldots,\Lambda_{M-1}\}$, $\hat\betab^n=\Mb_B^{-1}(\Xb_n)\Psib_n^T \Sigmab_n^{-1}\yb_n$. According to the BLM \eqref{BLM2}, an optimal design for the estimation of $I_\mu(f)$, or equivalently for the estimation of $\beta_0'$, is a Bayesian c-optimal design minimizing $\{\Mb_B^{-1}(\Xb_n)\}_{1,1}$; see \cite{Pilz83}. The design criterion is different when one is interested in approximation rather than integration: as shown in \cite{GP-CSDA2016}, a design minimizing the integrated MSPE in this model is an A-optimal design minimizing $\tr[\Mb_B^{-1}(\Xb_n)]$.

Approximate design theory can also be used. There, an optimal design $\xi^*$ is a probability measure on~$\SX$ that minimizes $\{\Mb_{B,m}^{-1}(\xi)\}_{1,1}$ with respect to $\xi\in\SM^+(1)$, where
\be\label{MB}
\Mb_{B,m}(\xi) = \int_\SX \frac{1}{v(\xb)} \, \psib(\xb)\psib^T(\xb) \dd\xi(\xb) + \Lambdab^{-1}/m\,, \quad m>0\,,
\ee
with $\psib(\xb)=(1,\varphi_1(\xb),\ldots,\varphi_M(\xb))^T$. The scalar $m$ defines a projected number of observations, in the sense that $\Mb_{B,n}(\xi_{n,e})=(1/n) \Mb_B(\Xb_n)$ when $\xi_{n,e}$ is the empirical measure associated with $\Xb_n$. The determination of $\xi^*$ forms a convex optimization problem \cite{Pilz83} for which many algorithms are available; see, e.g., \cite[Chap.~9]{PP2013}. Note that $\xi^*$ depends both on $m$ and $M$. When $M$ is small enough, $\xi^*$ has generally few support points; those points can be used for the construction of an exact design $\Xb_n$ and an extraction procedure is described in \cite{GP-CSDA2016, P-RESS2018}. Choosing $m$ and $M$ of the same order of magnitude as the projected number $n$ of observations is recommended. The construction of $\xi^*$ can also be sequential, using a vertex-direction algorithm similar to those in Section~\ref{S:any-time}. In general, the eigenfunctions $\varphi_k$ and eigenvalues $\Lambda_k$ corresponding to a given kernel $K$ (here $K_\mu$) are unknown and a numerical construction is required. It may rely on the substitution of a quadrature approximation for $\mu$ in the definition of the integral operator $T_\mu$; the construction is much facilitated when a tensor-product kernel can be used: the $\varphi_k$ and $\Lambda_k$ are then direct products of one-dimensional eigenfunctions and eigenvalues. In that case, the BLM \eqref{BLM2} also allows for straightforward estimation of the so-called Sobol' indices used in sensitivity analysis; see \cite{P-RESS2018}.

\begin{figure}[ht!]
\begin{center}
\includegraphics[width=.49\linewidth]{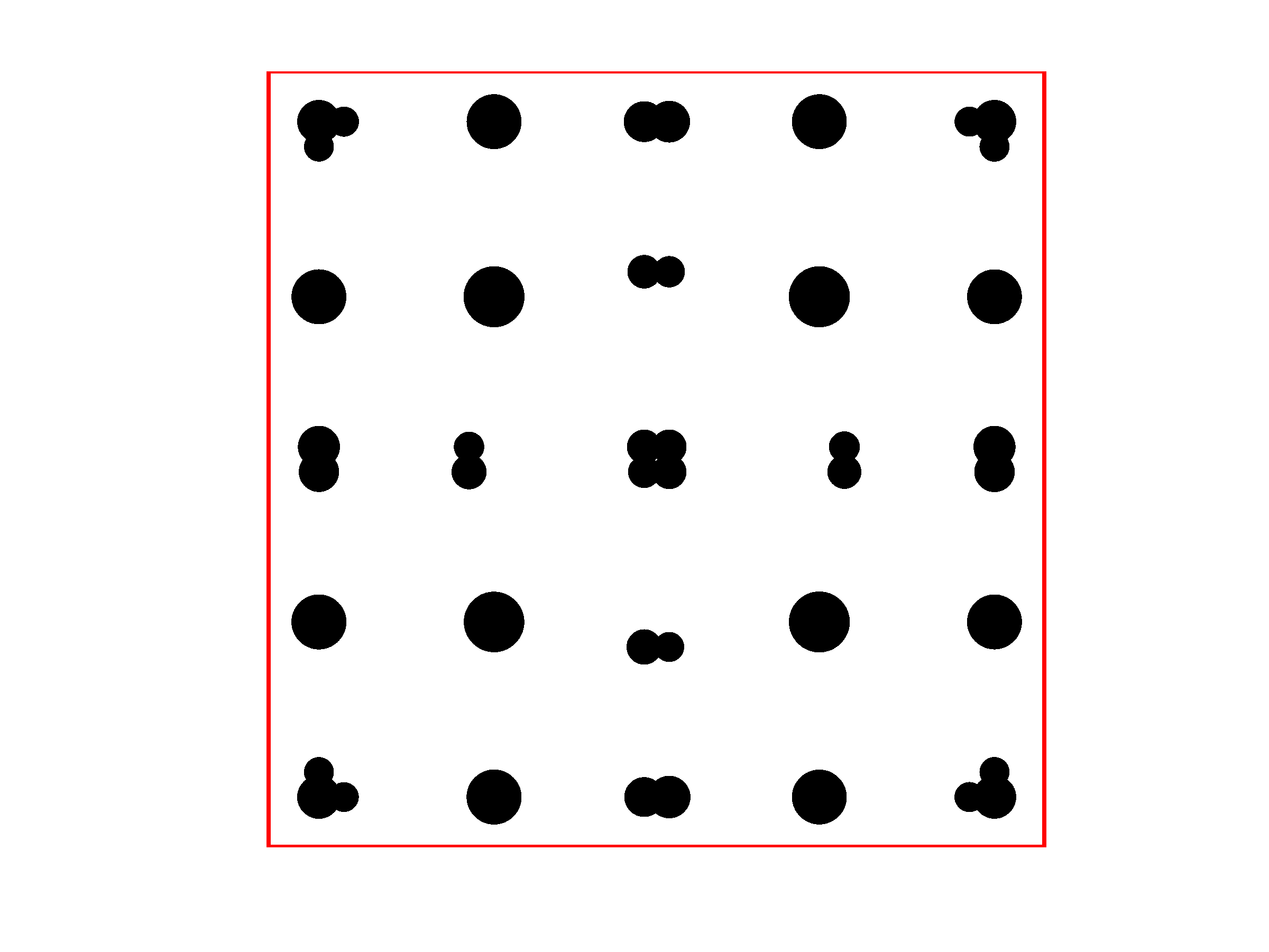} \includegraphics[width=.49\linewidth]{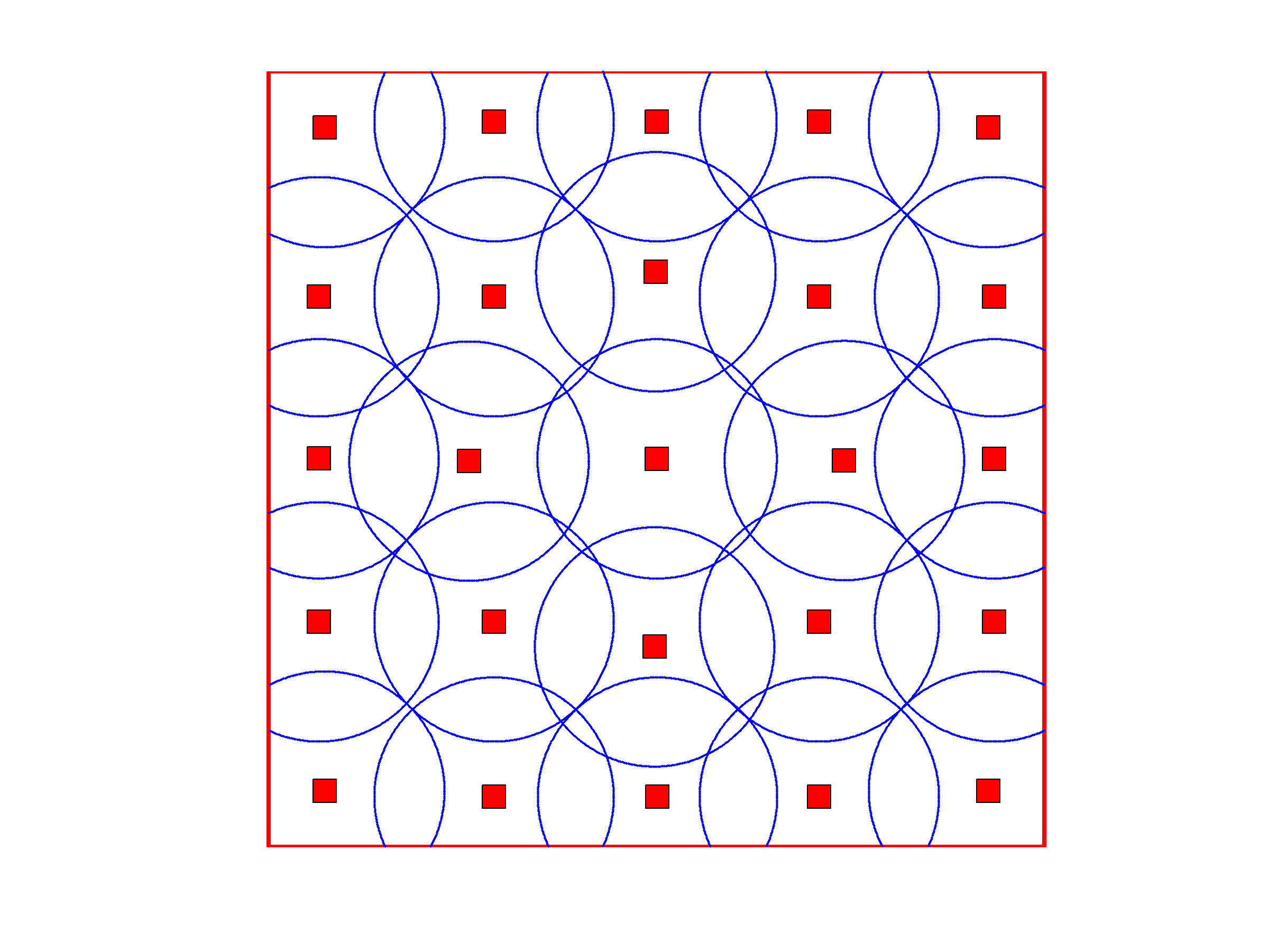}
\end{center}
\caption{\small Left: design measure $\xi^*$ minimizing $\{\Mb_{B,m}^{-1}(\xi)\}_{1,1}$ (44 support points with weights proportional to disk areas); Right: exact design $\Xb_{25}$ extracted from $\xi^*$, the circles have radius $\CR(\Xb_n)$.}
\label{F:designs-theta2_m15_M15}
\end{figure}

\noindent{\bf Example.}
We present a simple illustrative example where $\SX=[0,1]^2$, $\mu$ is uniform on~$\SX$ and $K$ is the tensor product of Mat\'ern 3/2 kernels \eqref{K32} with $\mt=2$. We consider designs supported on $\SX_\mO$ formed by a $32 \times 32$ regular grid in~$\SX$; $m=M=15$ in \eqref{MB}. Thanks to the the tensor structure, numerical approximations of eigenfunctions and eigenvalues can be obtained through one-dimensional quadrature approximations; we use here the uniform distribution on the 100 points $0,1/99,2/99,\ldots,1$; see \cite{P-RESS2018}.
The optimal design measure $\xi^*$ is constructed with a vertex-exchange algorithm \cite{Bohning86} with optimal step size, initialized at the uniform design on $\SX_\mO$. Iterations are stopped at iteration $k$ when
$$
\min_{\xb\in\SX_\mO} \frac{\mp \left(\{\Mb_{B,m}^{-1}[(1-\ma)\xi^{(k)}+\ma\delta_\xb]\}_{1,1}\right)}{\mp\ma}\bigg|_{\ma=0}>-10^{-6} \,.
$$
Figure~\ref{F:designs-theta2_m15_M15}-Left shows the corresponding measure $\xi^{(k)}$: there are 44 support points, their weights are proportional to the disk areas shown on the figure. Figure~\ref{F:designs-theta2_m15_M15}-right presents the exact design $\Xb_n$ extracted, with $n=25$ points; the circles have radius equal to $\CR(\Xb_n)$ and illustrate the good covering of~$\SX$ by $\Xb_{25}$. In general, the adjustment of $m$ and $M$ to obtain a design $\Xb_n$ with good space-filling properties for a prescribed value of $n$ remains a difficult task and requires further investigations.
\fin

\vsp
{\bf Acknowledgements. }
The first author thanks the Isaac Newton Institute for Mathematical Sciences, Cambridge, for support and hospitality during the programme UQ for inverse problems in complex systems where work on this paper was partly undertaken. This work was supported by EPSRC grant no EP/K032208/1.
We gratefully acknowledge S\'ebastien Da Veiga (Safran, Paris) who draw our attention to the machine learning literature on kernel herding and to connection with minimization of $L_2$ discrepancy.


\bibliographystyle{plain}
\bibliography{xampl,test}

\end{document}